\documentclass[12pt]{article}
\font\smallit=cmti10

\usepackage{amssymb,latexsym,epsfig}

\makeatletter

\renewcommand\section{\@startsection {section}{1}{\z@}
{-30pt \@plus -1ex \@minus -.2ex}
{2.3ex \@plus.2ex}
{\normalfont\normalsize\bfseries}}

\renewcommand\subsection{\@startsection{subsection}{2}{\z@}
{-3.25ex\@plus -1ex \@minus -.2ex}
{1.5ex \@plus .2ex}
{\normalfont\normalsize\bfseries}}

\renewcommand{\@seccntformat}[1]{\csname the#1\endcsname. }

\makeatother

\begin{document}

\begin{center}
{\bf THE 3x+1 PERIODICITY CONJECTURE IN $\mathbb{R}$}
\vskip 20pt
{\bf Josefina L\'opez}\footnote{E-mail: josefinapedro@hotmail.com}\\
{\smallit Camataqu\'i, Sud Cinti, Chuquisaca, Bolivia}\footnote{There is no institution sponsoring us.
It is the name of the little village in Bolivia where we live until 2019.} \\
\vskip 10pt
{\bf Peter Stoll}\footnote{Dr. phil. nat., University of Bern, Switzerland, (1942-2019).} \\
{\smallit Camataqu\'i, Sud Cinti, Chuquisaca, Bolivia}\\
\end{center}
\vskip 30pt

\centerline{\bf Abstract}

\noindent
The $3x+1$ map $T$ is defined on the 2-adic integers $\mathbb{Z}_2$ by $T(x)=x/2$ for even $x$ and $T(x)=(3x+1)/2$ for odd $x$. It is still unproved that under iteration of $T$ the trajectory of any rational 2-adic integer is eventually cyclic. A $2$-adic integer is rational if and only if its representation with $1$'s and $0$'s is eventually periodic. We prove that the $3x+1$ conjugacy $\Phi$ maps aperiodic $x\in\mathbb{Z}_2$ onto aperiodic $2$-adic integers provided that $\underline{\lim}\; (\frac{h}{\ell})_{\ell=1}^{\infty}>\frac{\ln(2)}{\ln(3)}$ where $h$ is the number of $1$'s in the first $\ell$ digits of $x$ with the following constraint: if there is a rational $2$-adic integer with a non-cyclic trajectory, then necessarily $\underline{\lim}\; (\frac{h}{\ell})_{\ell=1}^{\infty}=\frac{\ln(2)}{\ln(3)}$. We study $\Phi$ as an infinite series in $\mathbb{R}$ and obtain negative irrational numbers for which we compute their aperiodic $2$-adic expansion. We find prominent behaviors of the orbit of $x$ taking Sturmian words as parity vector. We also found amazing results of the terms of $\Phi$ in $\mathbb{R}$. We define the $\ell$'th iterate of $T$ for $\ell\rightarrow \infty$ in the ring of $3$-adic integers and obtain positive irrational numbers for which we compute their aperiodic $3$-adic expansion. 

\pagestyle{myheadings}

\thispagestyle{empty}
\baselineskip=15pt
\vskip 30pt

\newtheorem{thm}{Theorem} 
\newtheorem{defin}{Definition} 
\newtheorem{lem}{Lemma} 
\newtheorem{exa}{Example} 
\newtheorem{prop}{Proposition}
\newtheorem{cor}{Corollary} 
\newtheorem{esq}{Scheme}
\newtheorem{res}{Result}


\section{Introduction Part (I)}
Let $\mathbb{Z}_2$ denote the ring of 2-adic integers. Each $x\in\mathbb{Z}_2$ can be expressed uniquely as an infinite word $x=x_0x_1x_2\cdots$ of 1's and 0's. The $x_k$ are the digits of $x$, written from left to right. Let $0\leq d_0<d_1<d_2<\cdots$ be a finite or infinite sequence of non-negative integers defined by $d_i:=k$ whenever $x_k=1$ for a 2-adic integer $x=x_0x_1\cdots x_k\cdots$. Then $x$ can be written as the finite or infinite sum $x=2^{d_0}+2^{d_1}+2^{d_2}+\cdots$. Under 2-adic evaluation we have, for instance 
\begin{displaymath}
111111\cdots=2^0+2^1+2^2+\cdots=\frac{1}{1-2}=-1,
\end{displaymath}
\begin{displaymath}
1101010\cdots=2^0+2^1+2^3+2^5+\cdots=2^0+\frac{2^1}{1-2^2}=\frac{1}{3}.
\end{displaymath}

The $3x+1$ map $T$ is defined on the 2-adic integers $\mathbb{Z}_2$ by $T(x)=x/2$ for even $x$ and $T(x)=(3x+1)/2$ for odd $x$.\footnote{even: $x_0=0$; odd: $x_0=1$.} The \textit{trajectory} $\mathcal{T}(x)$ is the infinite sequence
\begin{displaymath}
\mathcal{T}(x):=(T^0(x)=x,T^1(x),T^2(x),\ldots,T^k(x),\ldots )
\end{displaymath}
where $T^k(x)$ denotes the $k$'th iterate of $T$. The set
\begin{displaymath}
\mathcal{O}(x):=\{T^k(x):k\in\mathbb{N}_0\}
\end{displaymath}
is called the \textit{orbit} of $x$. The orbit is finite if and only if the trajectory is eventually cyclic: $T^i(x)=T^j(x)$ for some $i\neq j$. If $T^k(x)=x$ for some $k\geq 1$, then the trajectory is purely cyclic, and $\mathcal{O}(x)$ contains the elements of a finite \textit{cycle}. If $\mathcal{O}(x)$ is an infinite set, the trajectory of the 2-adic integer $x$ is called \textit{divergent}.

Let $\mathbb{Q}_{odd}$ denote the ring of rational numbers having odd denominators in reduced fraction form. This ring is isomorphic to the subring $\mathbb{Q}_2\subset\mathbb{Z}_2$ of eventually periodic 2-adic integers. Thus $\mathbb{Q}_{odd}$ is the ring of \textit{rational 2-adic} integers. The question ab to whether there is a rational 2-adic integer with a divergent trajectory is part of the \textit{periodicity conjecture} raised in 1985 by Jeffrey C. Lagarias (\cite{Lag:1985}). Following (\cite{Lag:1985}), the trajectory of $x$ evaluated modulo 2 is called the \textit{parity vector} of $x$:
\begin{displaymath}
v(x)=(T^k(x)\bmod{2})_{k=0}^\infty. 
\end{displaymath}
The parity vector is an infinite word over the alphabet $\{0,1\}$. If $x\in\mathbb{Z}_2$, $v(x)$ is interpreted as a 2-adic integer as well, namely
\begin{displaymath}
v(x)=v_0v_1\cdots v_k\cdots = \sum_{k=0}^{\infty}(T^k(x) \,\bmod{2})\cdot2^k,
\end{displaymath}
and $x\mapsto v(x)$ defines a map $\mathbb{Z}_2\rightarrow\mathbb{Z}_2$. This map is called $Q_\infty$ in (\cite{Lag:1985}), where it is proved that $Q_\infty$ is bijective and continuous under the 2-adic metric (\cite{Lag:1985}, Theorem L). In 1994 Daniel J. Bernstein found the inverse map $x\leftarrow v(x)$, which he called $\Phi$ (\cite{Ber:1994}). Map $Q_\infty$ will be later referred to as $\Phi^{-1}$. The explicit formulas for homeomorphism $\Phi:\mathbb{Z}_2\rightarrow\mathbb{Z}_2$ are as follows:
\begin{equation}\label{PHI}
\Phi(v)=\Phi(2^{d_0}+2^{d_1}+2^{d_2}+\cdots)=-\sum_{i=0}^{\infty}\frac{1}{3^{i+1}}2^{d_i}\\[-12pt]
\end{equation}
\begin{equation}\label{PHI-}
\Phi^{-1}(x)=\sum_{k=0}^{\infty}(T^k(x) \,\bmod{2})\cdot2^k
\end{equation} 

The periodicity conjecture states that $\Phi^{-1}(Q_2)=Q_2$. A proof of this conjecture must have two parts: $\Phi^{-1}(Q_2)\supset Q_2$ and $\Phi^{-1}(Q_2)\subset Q_2$. The first inclusion states that any eventually periodic $y\in\mathbb{Z}_2$ is parity vector of some eventually periodic $x\in\mathbb{Z}_2$. Since $\Phi^{-1}(Q_2)\supset Q_2$ is equivalent to $Q_2\supset\Phi(Q_2)$, the proof is easily done by computing the corresponding geometric series (\ref{PHI}) under 2-adic evaluation. However, the second inclusion which states that the parity vector of any eventually periodic $x\in\mathbb{Z}_2$ is eventually periodic is still unproved. Hence we have the unsolved problem: Is it true that for any \textit{aperiodic} word $v$ there is no $x\in\mathbb{Q}_{odd}$, having $v$ as parity vector? 

In 1996 D. J. Bernstein and J. C. Lagarias proved that $\Phi$ is a \textit{conjugacy map}.  
The 2-adic shift map $S$ is defined on the 2-adic integers $\mathbb{Z}_2$ by $S(x)=x/2$ for even $x$ and $S(x)=(x-1)/2$ for odd $x$.
$T$ and $S$ are conjugates (\cite{BerLag:conjugacy}). There exists a unique homeomorphism $\Phi:\mathbb{Z}_2\to\mathbb{Z}_2$ (the $3x+1$ conjugacy map) with $\Phi(0)=0$ and
\begin{equation}\label{CO}
\Phi\circ S \circ\Phi^{-1}=T.\\[-5pt]
\end{equation}
\noindent Map $\Phi$ is a 2-adic isometry (\cite{BerLag:conjugacy}):
\begin{equation}\label{ISO}
|\Phi(x)-\Phi(y)|_2=|x-y|_2 \qquad\textrm{for all } x,y\in \mathbb{Z}_2.\\[-10pt]
\end{equation}
\noindent
Moreover (\cite{BerLag:conjugacy}),
\begin{equation}\label{MOD}
\Phi(x)\equiv x\pmod{2} \qquad\textrm{for all } x\in \mathbb{Z}_2.
\end{equation}
Since $\Phi\circ S=T\circ \Phi$, we have the following commutative scheme:
\begin{esq}\label{esquema1}
\begin{center}
\begin{tabular}{rrccccl}
$v=$ & $v_0v_1v_2v_3\cdots$ &  & $\longrightarrow$ &  & $\Phi(v)$ & $:=\zeta$ \\  
 & $S \downarrow$ &  & $\Phi$ &  & $\downarrow T$ &  \\ 
& $v_1v_2v_3\cdots$ &  & $\longrightarrow$ &  & $T(\zeta)$ &  \\ 
 & $S \downarrow$ &  & $\Phi$ &  & $\downarrow T$ &  \\ 
 & $v_2v_3\cdots$ &  & $\longrightarrow$ &  & $T^2(\zeta)$ &  \\ 
& $\downarrow$ &  & $\Phi$ &  & $\downarrow$ &  \\
\end{tabular}
\end{center}
\end{esq}

Let $v$ be any aperiodic word. Then we have $T^i(\zeta)\neq T^j(\zeta)$ for all $i\neq j$ since $\Phi$ is injective. Hence $\mathcal{O}(\zeta)\in\mathbb{Z}_2$ is an infinite set, and the trajectory of $\zeta$ is divergent. If $\zeta\in\mathbb{Q}_{odd}$, then the infinite set $\mathcal{O}(\zeta)\subset\mathbb{Q}_{odd}\subset\mathbb{R}$ cannot have accumulation points by (\cite{Monks:Yaz}, Lemma 3.2). Throughout this paper we look for accumulation points in the set $\mathcal{O}(\zeta)$. If we can find one, then $\zeta\notin\mathbb{Q}_{odd}$. 

If $\zeta\in\mathbb{Q}_{odd}$, then the trajectory $\mathcal{T}(\zeta)$ is strictly positive for $\zeta>0$ and eventually nonnegative for $-1<\zeta\leq 0$. However, the trajectory for $\zeta\leq -1$ is either strictly negative with all its elements being $\leq -1$, or eventually nonnegative with elements going beyond $-1$ by (\cite{Monks:Yaz}, Lemma 3.4).

Scheme \ref{esquema1} allows a coherent projection of the $3x+1$ periodicity conjecture from $\mathbb{Z}_2$ down to $\mathbb{R}$. To do so, we define a real-valued clone of $T$ as follows. Let $u=v_0v_1\cdots v_i\cdots v_{\ell-1}$ ($\ell\geq 1$) be a prefix of an aperiodic word $v$. The function $T_u$ is defined for all $x\in\mathbb{R}$ and computed by following the instruction $u$; that is $T_{v_i}:=(3x+1)/2$ for $v_i=1$, $T_{v_i}:=x/2$ for $v_i=0$, and $T_u(x):=T_{v_{\ell-1}}\circ\cdots\circ T_{v_i}\circ\cdots T_{v_1}\circ T_{v_0}(x)$.\footnote{The function $T_u(x)$ is the affine function $\phi:=\mathbb{R}\rightarrow\mathbb{R} $ of L. Halbeisen and N. Hungerb\H{u}hler in \cite{Halb:1997}} Obviously, $T_u(x)$ ignores the parity of $x$. The \textit{pseudo trajectory} of $x\in\mathbb{R}$ relative to $v$ is defined by $\mathcal{T}_v(x):=(x, T_{v_0}(x), T_{v_0v_1}(x), T_{v_0v_1v_2}(x),\ldots)$.

If the series (\ref{PHI}) has a limit $\Phi_{\mathbb{R}}(v):=\zeta$ in $\mathbb{R}$, then we have the following commutative scheme:
\begin{esq}\label{esquema2}
\begin{center}
\begin{tabular}{rrccccl}
$v=$ & $v_0v_1v_2v_3\cdots$ &  & $\longrightarrow$ &  & $\Phi_{\mathbb{R}}(v)$ & $:=\zeta$ \\  
 & $S \downarrow$ &  & $\Phi_{\mathbb{R}}$ &  & $\downarrow T_{v_0}$ &  \\ 
& $v_1v_2v_3\cdots$ &  & $\longrightarrow$ &  & $T_{v_0}(\zeta)$ &  \\ 
 & $S \downarrow$ &  & $\Phi_{\mathbb{R}}$ &  & $\downarrow T_{v_1}$ &  \\ 
 & $v_2v_3\cdots$ &  & $\longrightarrow$ &  & $T_{v_1}\circ T_{v_0}(\zeta)$ & $=T_{v_0v_1}(\zeta)$\\ 
& $\downarrow$ &  & $\Phi_{\mathbb{R}}$ &  & $\downarrow T_{v_2}$ &  \\
\end{tabular}
\end{center}
\end{esq}

Let $v=v_0v_1v_2v_3\cdots$ be an aperiodic word having the limit $\zeta=\Phi_{\mathbb{R}}(v)$. By Scheme \ref{esquema2} the pseudo trajectory $\mathcal{T}_v(\zeta)=(\zeta,T_{v_0}(\zeta),T_{v_0v_1}(\zeta),T_{v_0v_1v_2}(\zeta),\ldots)$ and the trajectory $\mathcal{T}(\zeta)$ defined by $\mathcal{T}(\zeta):= (\Phi_{\mathbb{R}}(v_0v_1v_2\cdots),\Phi_{\mathbb{R}}(v_1v_2v_3\cdots),\Phi_{\mathbb{R}}(v_2v_3v_4\cdots),\ldots)$ are the same sequence, and its terms are negative real numbers less than $-1$ (Lemma \ref{limarriba}). 

For $\zeta\in\mathbb{Q}_{odd}$ in $\mathbb{R}$ Scheme \ref{esquema2} does the same as Scheme \ref{esquema1} in $\mathbb{Z}_2$ since $v$ is the parity vector of $\zeta$. For $\zeta\notin\mathbb{Q}_{odd}$ it is an extension for a  value whose parity is undefined. For an aperiodic $v$ the elements of the pseudo trajectory $\mathcal{T}_v(\zeta)$ are all different from each other, so that the orbit $\mathcal{O}(\zeta)\subset (-\infty,-1)$ is an infinite set (Lemma \ref{infinito}). If this set has accumulation points, then $\zeta\notin\mathbb{Q}_{odd}$. 

We had originally thought that such accumulation points must exist for all aperiodic words having a limit $\Phi_{\mathbb{R}}(v)$. However, Example \ref{acht} is a counterexample. We furthermore thought that with a detailed knowledge of the structure of $v$ it would be relatively simple to find out some accumulation point. We were wrong again. Example \ref{On} is such a case where we cannot guarantee a (probably existing) accumulation point. Our hope that things in $\mathbb{R}$ would be much easier than in $\mathbb{Z}_2$ where a sequence of eventually periodic 2-adic integers can have an eventually periodic 2-adic integer as 2-adic accumulation point ---for instance, the sequence $\big(\frac{3+2^n}{7}\big)_{n=1,2,3,\cdots}$ converges to $3/7$)--- was an illusion.  

Here it is convenient to use a function introduced by L. Halbeisen and N. Hungerb\H{u}hler in \cite{Halb:1997}. Let $\mathbb{U^*}$ be the set of all \textit{finite} words $u$ over the alphabet $\{0,1\}$. Function $\varphi:\mathbb{U^*} \to \mathbb{N}_0$ is defined recursively by
\begin{eqnarray*}
\varphi(\varepsilon)&=&0; \\
\varphi(u0)&=&\varphi(u);\\
\varphi(u1)&=&3\varphi(u)+2^{\ell(u)}.
\end{eqnarray*}
where $\varepsilon$ is the empty word and $\ell(u)$ the length of $u$. Let $h(u)$ be the \textit{height} of $u$; that is, the number of 1's in $u$. Using the pointer notation $d_i$, we have (\cite{Halb:1997})
\begin{equation}\label{fi}
\varphi(u)=3^{h(u)}\sum_{i=0}^{h(u)-1} \frac{1}{3^{1+i}} 2^{d_i}.
\end{equation}
Further (\cite{Halb:1997}), for all $u,u'\in\mathbb{U^*}$, 
\begin{equation}\label{UV}
\varphi(uu')= 3^{h(u')}\varphi(u)+2^{\ell(u)} \varphi(u').
\end{equation} 

It can be proved that $\varphi$ is injective. Thus $\varphi(u)$ encodes the structure of $u$ with a non-negative integer. We define $\mathbb{U}:=\mathbb{U^*}-\varepsilon$. Then $\varphi(u)$ is not divisible by $3$ for all $u\in\mathbb{U}$; it is divisible by $2^n$ only if the first $n$ digits of $u\in\mathbb{U}$ are 0's.

Let $v=v_0v_1\cdots$ be an aperiodic word. In Section \ref{proof}, we associate to any prefix $u\in\mathbb{U}$ of $v$ with length $\ell$ and height $h$ the following three magnitudes: 
\begin{equation}\label{tres}
\Phi(u)=-\frac{\varphi(u)}{3^h},\qquad T_u(0)=\frac{\varphi(u)}{2^\ell},\qquad C(u)=\frac{\varphi(u)}{2^\ell-3^h},
\end{equation}
where $C(u)$ is a well known formula (written with $\varphi$) for finite cycles of length $\ell$; height $h$, $\Phi(v)$ is the sum of the first $h$ terms in the infinite series (\ref{PHI}), and $T_u(0)$ is the sum of the first $h$ non-zero terms in the finite series
\begin{equation}\label{1985}
T_u(0)=\sum_{i=0}^{h-1}\frac{3^{h-1-i}2^{d_i}}{2^{\ell}}=\sum_{i=0}^{\ell-1}v_i\cdot\frac{3^{v_i+v_{i+1}+\cdots +v_{\ell-1}-1}}{2^{\ell-i}}
\end{equation}
which corresponds to $\varrho_k(n)$ of relation (2.6) in \cite{Lag:1985}.

It is easily seen that we have the following \textit{symbolic} relations:
\begin{equation}\label{rel}
\Phi(u)=C(u)\bigg(1-\frac{2^\ell}{3^h}\bigg),\qquad T_u(0)=C(u)\bigg(1-\frac{3^h}{2^\ell}\bigg),\qquad T_u(0)=-\frac{3^h}{2^\ell}\Phi(u).
\end{equation}

The first relation is useful in $\mathbb{Z}_2$, because $\Phi(u)$, $C(u)$ and $2^\ell/3^h$ are 2-adic integers; the second one is useful in the ring $\mathbb{Z}_3$, because $T_u(0)$, $C(u)$ and $3^h/2^\ell$ are 3-adic integers; and the third one involves 2-adic or 3-adic fractions. However, the three relations are well-defined in $\mathbb{R}$ since $\Phi(u)$, $T_u(0)$ and $C(u)$ are rational numbers. Obviously, 
\begin{equation}\label{modulo}
C_{\mathbb{R}}(u)\equiv \Phi_{\mathbb{R}}(u)\pmod{2^\ell}\quad\textrm{and}\quad C_{\mathbb{R}}(u)\equiv T_u(0)\pmod{3^h}.
\end{equation}

\begin{exa}\label{null}
Suppose that $u=100101$ is the prefix of an aperiodic word $v$. Then $u$ is a prefix of the parity vector of $\Phi_{\mathbb{R}}(u)= -\frac{1}{3}-\frac{2^3}{3^2}-\frac{2^5}{3^3}=-\frac{65}{27}$, $u^{\infty}=uuu\cdots$ is the purely periodic parity vector of $C_{\mathbb{R}}(u)=\frac{65}{37}$, and $T_u(0)= \frac{3^2}{2^6}+\frac{3}{2^3}+\frac{1}{2}=\frac{65}{64}$ gives the output when instruction $u$ is applied on $0$. Indeed, we have the trajectories 

\noindent $\mathcal{T}(-\frac{65}{27})=\big(-\frac{65}{27},-\frac{28}{9},-\frac{14}{9},-\frac{7}{9},-\frac{2}{3},-\frac{1}{3},0,0,\ldots\big)$, and $\mathcal{T}(-\frac{65}{27})\bmod{2}=10010100\cdots$.\\$\mathcal{T}(\frac{65}{37})=\big(\frac{65}{37}, \frac{116}{37}, \frac{58}{37}, \frac{29}{37}, \frac{62}{37}, \frac{31}{37}, \;\;\frac{65}{37}, \frac{116}{37},\;\ldots\big)$, and $\mathcal{T}(\frac{65}{37})\bmod{2}=100101\; 100101\cdots$. Further, $\mathcal{T}_u(0)=\big(0, \frac{1}{2}, \frac{1}{4}, \frac{1}{8}, \frac{11}{16}, \frac{11}{32}, \frac{65}{64}\big)$. Finally,

\noindent $-\frac{65}{27}\equiv 45\pmod{64}$  yields the first integer of the trajectory $(45, 68, 34, 17, 26, 13, 20,\ldots)$, so that the parity vector of $45$ has $100101$ as prefix, and $T_u(0)=\frac{65}{64}\equiv\frac{65}{37}\equiv 20\pmod{27}$ yields the seventh integer, that is, the output of $100101$ when $45$ is its input.
\end{exa}

It was stated first by J. C. Lagarias (\cite{Lag:1985}, Relation 2.31) that the existence of an integer $n\in\mathbb{Z}$ with divergent trajectory to $\pm\infty$ implies
\begin{equation}\label{mayor}
\underline{\lim}\; \bigg(\frac{h}{\ell}\bigg)\geq\frac{\ln(2)}{\ln(3)}\\
\end{equation}
where $h$ is the number of 1's in the first $\ell$ digits of the parity vector $v(n)$. In 2004 it was proved that the same constraint holds for rational 2-adic integers as well (Monks, Yazinski \cite{Monks:Yaz}, Theorem 2.7 b): if the trajectory of $\zeta\in\mathbb{Q}_{odd}$ is divergent (that is, if the orbit is an infinite set), then it holds (\ref{mayor}). The contrapositive of this statement tells us that
\begin{equation}\label{menor}
\underline{\lim}\; \bigg(\frac{h}{\ell}\bigg)<\frac{\ln(2)}{\ln(3)}\\
\end{equation}
implies $\zeta\notin\mathbb{Q}_{odd}$, or $\mathcal{O}(\zeta)$ is finite. However, if $\Phi(v)=\zeta$ for an aperiodic word $v$, then $\mathcal{O}(\zeta)$ is an infinite set. Thus (\ref{menor}) implies $\zeta\notin\mathbb{Q}_{odd}$.

Our main result is a constraint stronger than (\ref{mayor}).\footnote{S. Eliahou (\cite{Eliahou:1993}) Observation: the inverse proportion of odd elements in a cycle is very close  to $\log _2(3)$, in a precise sense.}
\begin{thm}\label{aperiodic}
If the trajectory of some $\zeta\in\mathbb{Q}_{odd}$ is divergent (that is, if the orbit is an infinite set), then it holds that
\begin{equation}\label{equalR}
\underline{\lim}\; \bigg(\frac{h}{\ell}\bigg)=\frac{\ln(2)}{\ln(3)},\\
\end{equation}
where $h$ is the number of $1$'s in the first $\ell$ digits of the parity vector $v$ of $\zeta$.
\end{thm}

The proof of Theorem \ref{aperiodic} is in Section \ref{main proof}. It involves results of a former paper (\cite{Lopez:Lop}) about the $3x+1$ conjugacy map over \textit{Sturmian words}. Sturmian words with irrational slope $0<\alpha<1$ are the simplest aperiodic words because of their similarity with purely periodic words. Specifically, Sturmian words with irrational slope $1>\alpha\geq \ln(2)/\ln(3)$ satisfy the constraint (\ref{mayor}) in contrast with the binary expansion of classical irrationals like $\pi$, $e$, $\sqrt{2}$, or $\log(2)$ which seem to be \textit{normal} numbers such that ${\lim}_{\ell\rightarrow\infty}\;(h/\ell)=1/2$. 

The special Sturmian words $1c_{\alpha}$ and $0c_{\alpha}$ are defined in (\cite{Lot:2002}) by 
\begin{displaymath}
1c_\alpha:=\lceil(j+1)\alpha\rceil-\lceil j\alpha\rceil, \quad 0c_\alpha:=\lfloor(j+1)\alpha\rfloor-\lfloor j\alpha\rfloor \quad\textrm{for }j=0,1,2,\ldots\quad (0<\alpha<1).
\end{displaymath}  

\begin{figure}[htbp]
\begin{center}
\epsfxsize=3.5in
\epsfbox{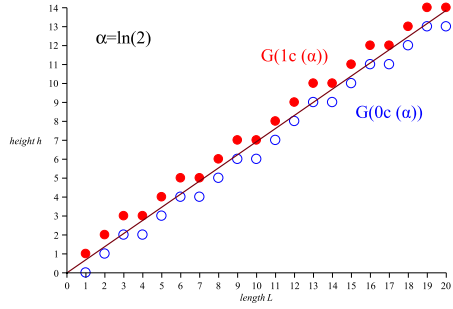}
\end{center}
\caption{Geometric representation of $1c_{\alpha}$ and $0c_{\alpha}$}
\label{intro}
\end{figure}

In Figure \ref{intro} we show the geometric representation of these words. They are aperiodic if and only if $0<\alpha<1$ is an irrational number. For instance, if $\alpha=\ln(2)$, we have $1c_{\alpha}=1110110110\cdots$. The number of $1$'s in the first $\ell$ digits (called \textit{height} $h$) is represented as a function of the length $\ell$ of the prefix. Thus $G(1c_{\;\ln(2)}):=(\ell,h)_{\ell=1}^\infty=\big((1,1),(2,2),(3,3),(4,3),(5,4),(6,5),(7,5),\ldots\big)$. The \textit{slope} $\big(h/\ell\big)_{\ell=1,2,\ldots}$ converges to $\ln(2)$.
In (\cite{Lopez:Lop}, Section 4) we discussed the devil's staircases $F$ and $F^*$ which show the irrationality of $\Phi_{\mathbb{R}}(1c_{\alpha})$ (for $1>\alpha > \ln(2)/\ln(3)$) and $\Phi_{\mathbb{R}}^*(1c_{\alpha})$ (for $0<\alpha < \ln(2)/\ln(3)$). The evaluation of $\Phi(1c_{\alpha})$ in $\mathbb{R}$ gave us the hint on how to prove the main result (\cite{Lopez:Lop}, Theorem 1, Corollary 2) in the 2-adic world. Otherwise we would not have been aware of the function \textit{gap} and the rather strange function $\Phi^*$ dual to $\Phi$. So far as well, but how can we relate the irrational real number $\Phi_{\mathbb{R}}(1c_{\alpha})$ with the $2$-adic ring $\mathbb{Z}_2$?

\section{Introduction Part (II)}

It is well known that the ring of $p$-adic integers is the \textit{inverse limit} of the rings $\mathbb{Z}/p^n\mathbb{Z}$ and therefore, a $p$-adic integer is a sequence $(x_n)_{n=1}^{\infty}$ of non-negative integers such that $x_n\in\mathbb{Z}/p^n\mathbb{Z}$ and if $n\leq m$, then $x_n\equiv x_m \pmod{p^n}$.

Let $u^{({\ell)}}$ be a prefix with length $\ell$ of a given aperiodic word $v$, and suppose that the limit $\Phi_{\mathbb{R}}(v)$ of the series (\ref{PHI}) exists. The inverse limit of $\Phi_{\mathbb{R}}(u^{({\ell)}})$ for a fixed $\ell$ is the sequence
\begin{displaymath}
\lim_{\longleftarrow}\Phi_{\mathbb{R}}(u^{({\ell)}}):=\bigg(\Phi_{\mathbb{R}}\big(u^{({\ell)}}\big)\bmod{2^n}\bigg)_{n=1}^{\infty}\;.
\end{displaymath}
Since $\Phi_{\mathbb{R}}\big(u^{({\ell)}}\big)$ are the partial sums of $\Phi_{\mathbb{R}}(v)$, we have that
\begin{displaymath}
\Bigg(\bigg(\Phi_{\mathbb{R}}\big(u^{({\ell)}}\big)\bmod{2^n}\bigg)_{n=1}^{\ell}\Bigg)_{\ell=1}^{\infty}
\end{displaymath}
is a sequence $\big(X_{\ell}\big)_{\ell=1}^{\infty}$ of terms $X_{\ell}$ which are sequences of length $\ell$ such that $X_{\ell}$ is prefix of $X_{\ell+1}$ for $\ell=1,2,3,\dots$. For $\ell\rightarrow\infty$ we obtain an infinite sequence $X_{\infty}$, which is the inverse limit of $\Phi_{\mathbb{R}}(v)$. The elements of $X_{\infty}$ coincide with the ultimate integer of each $X_{\ell}$. Therefore,
the inverse limit of $\Phi_{\mathbb{R}}(v)$ is the sequence
\begin{equation}\label{IL}
\lim_{\longleftarrow}\Phi_{\mathbb{R}}(v)=\bigg(\Phi_{\mathbb{R}}\big(u^{({\ell)}}\big)\bmod{2^{\ell}}\bigg)_{\ell=1}^{\infty}\;,\quad\textrm{ and we have}\quad\lim_{\longleftarrow}\Phi_{\mathbb{R}}(v)=\Phi(v).
\end{equation}
If $\Phi_{\mathbb{R}}(v)\in\mathbb{Q}_{odd}$, then $\lim_{\leftarrow}\Phi_{\mathbb{R}}(v)\in\mathbb{Q}_2$. If $\Phi_{\mathbb{R}}(v)$ is an irrational number, then the inverse limit $\lim_{\leftarrow}\Phi_{\mathbb{R}}(v)$ is a $2$-adic integer which is aperiodic, because the irrational $\Phi_{\mathbb{R}}(v)$ does not admit an eventually periodic $2$-adic expansion. 

\begin{exa}\label{ejln2}
Let $\alpha=\ln(2)$. \\Then $1c_{\alpha}=1110110110111011011011\cdots$. Further, $\Phi(1c_{\alpha})=1110010101011101000010\cdots$ by (\cite{Lopez:Lop}, Example 9) and $\Phi_{\mathbb{R}}(1c_{\alpha})=-2.5970181822928929618\cdots$ by (\cite{Lopez:Lop}, Corollary 2 in $\mathbb{R}$).

\noindent $\big(\Phi_{\mathbb{R}}\big(u^{({\ell)}}\big)\big)_{\ell=1}^{\infty}= (-\frac{1}{3}, -\frac{5}{9}, -\frac{19}{27}, -\frac{19}{27}, -\frac{73}{81}, -\frac{251}{243}, -\frac{251}{243}, -\frac{881}{729}, -\frac{2899}{2187}, -\frac{2899}{2187}, -\frac{9721}{6561}, -\frac{31211}{19683}, \;\ldots)$.

\noindent $\big(\Phi_{\mathbb{R}}\big(u^{({\ell)}}\big)\bmod{2^{\ell}}\big)_{\ell=1}^{\infty}=(1, 3, 7, 7, 7, 39, 39, 167, 167, 679, 679, 2727, \ldots)$.

\noindent Then $\lim_{\leftarrow}\Phi_{\mathbb{R}}(1c_{\alpha})=111001010101\cdots$, and we have $\lim_{\leftarrow}\Phi_{\mathbb{R}}(1c_{\alpha})=\Phi(1c_{\alpha})$.

\noindent Moreover, we can go back to the parity vector $1c_{\alpha}$ with the finite trajectory $\mathcal{T}(2727)=(2727, 4091, 6137, 9206, 4603, 6905, 10358, 5179, 7769, 11654, 5827, 8741)$ reduced modulo $2$.
\end{exa}

If limit $\Phi_{\mathbb{R}}(v)$ does not exist, then the series (\ref{PHI}) diverges to $-\infty$, because the sequence of partial sums $\big(\Phi_{\mathbb{R}}(u^{(\ell)})\big)_{\ell=1}^{\infty}$ is a monotone decreasing sequence of negative rational numbers. Since the $2$-adic $\Phi(v)$ always exists for all $v\in\mathbb{Z}_2$, we have
\begin{equation}\label{IL2}
\Phi(v)=\bigg(\Phi_{\mathbb{R}}\big(u^{({\ell)}}\big)\bmod{2^{\ell}}\bigg)_{\ell=1}^{\infty}.
\end{equation}
From (\ref{IL}) and (\ref{IL2}) we conclude that the existence of $\Phi_{\mathbb{R}}(v)$ implies $\Phi_{\mathbb{R}}(v)\notin \mathbb{Q}\setminus\mathbb{Q}_{odd}$, because numbers of $\mathbb{Q}\setminus\mathbb{Q}_{odd}$ have no $2$-adic expansion. Hence, if $\Phi_{\mathbb{R}}(v)$  exists and $\Phi_{\mathbb{R}}(v)\notin\mathbb{Q}_{odd}$, then $\Phi_{\mathbb{R}}(v)$ is an irrational real number with an aperiodic $2$-adic expansion given by (\ref{IL}). 

Example \ref{ejln2} shows that the irrational real numbers of the form $\Phi_{\mathbb{R}}(v)$ have $2$-adic representations, i.e., the irrational $-2.5970181822928929618\cdots$ is represented by the aperiodic $2$-adic integer $1110010101011101000010\cdots$. All these irrationals are negative and less than $-1$. If the limit $\Phi_{\mathbb{R}}(v)$ exists, then $\Phi_{\mathbb{R}}(v)$ and $C_{\mathbb{R}}(v):=\lim_{u\rightarrow v}\big(C_{\mathbb{R}}(u)\big)$ are the same number by Lemma \ref{Fi-C}. Therefore, it is no surprise that the new sequence\\
\noindent $\big(C_{\mathbb{R}}\big(u^{({\ell)}}\big)\big)_{\ell=1}^{\infty}= (-1, -1, -1, -\frac{19}{11}, -\frac{73}{49}, -\frac{251}{179}, -\frac{251}{115}, -\frac{881}{473}, -\frac{2899}{1675}, -\frac{2899}{1163}, -\frac{9721}{4513}, -\frac{31211}{15587}, \;\ldots)$\\
\noindent yields the same inverse limit as in Example \ref{ejln2}:\\
\noindent $\big(C_{\mathbb{R}}\big(u^{({\ell)}}\big)\bmod{2^{\ell}}\big)_{\ell=1}^{\infty}=(1, 3, 7, 7, 7, 39, 39, 167, 167, 679, 679, 2727, \ldots)$.

It may be that $\lim_{u\rightarrow v}\big(C_{\mathbb{R}}(u)\big)$ exists only for selected prefixes $u$. For such sequential accumulation points we write  $C_{\mathbb{R}}(v,\mathcal{L}):=\lim\;\big(C_{\mathbb{R}}(u^{({\ell)}})\big)_{\ell\in\mathcal{L}}$, where $\ell$ takes progressively the values of a strictly monotone increasing sequence $\mathcal{L}$ of positive integers. In Example \ref{ejgolden} we show that there are irrational real numbers of the form $C_{\mathbb{R}}(v,\mathcal{L})$ which are positive.

\begin{exa}\label{ejgolden}
Let $\alpha=\frac{2}{1+\sqrt5}$ and $\big(\frac{p_k}{q_k}\big)_{k=0}^{\infty}=\big(0,1,\frac{1}{2},\frac{2}{3},\frac{3}{5},\frac{5}{8},\frac{8}{13},\frac{13}{21},\frac{21}{34},\frac{34}{55},\ldots\big)$ its convergents.  

\noindent Then $1c_{\alpha}=1101101011011010110101\cdots$. Further, $\Phi(1c_{\alpha})=1101111011100100011001\cdots$ by (\cite{Lopez:Lop}, Example 7) and $\Phi^*_{\mathbb{R}}(1c_{\alpha})=10.370127141747049162\cdots$ by (\cite{Lopez:Lop}, Example 3).

\noindent Sequence (16) in \cite{Lopez:Lop} converges very rapidly to $\Phi^*_{\mathbb{R}}(1c_{\alpha})$. Subsequence \footnote{We use $w_k$ instead of $v_k$ for  Christoffel words.}\\
\noindent $\big(\frac{2\varphi(w_k)-3^{p_k-1}}{2(2^{q_k}-3^{p_k})}\big)_{k \;\mathrm{odd}} =(-\frac{1}{2},-\frac{7}{2},\frac{557}{26},\frac{11078989}{1005658},\frac{4605785884670507258361037943}{444520808633569958378391605},\ldots)$ converges very fast to $\Phi^*_{\mathbb{R}}(1c_{\alpha})$. Since  $\frac{2\varphi(w_k)-3^{p_k-1}}{2(2^{q_k}-3^{p_k})}=\frac{\varphi(w_k)}{2^{q_k}-3^{p_k}}+g'(\frac{p_k}{q_k})$ and $g'(\frac{p_k}{q_k})\rightarrow 0$, the sequence $\big(\frac{\varphi(w_k)}{2^{q_k}-3^{p_k}}\big)_{k \;\mathrm{odd}}$ converges very fast to $\Phi^*_{\mathbb{R}}(1c_{\alpha})$ and is a subsequence of $\big(C_{\mathbb{R}}\big(u^{({\ell)}}\big)\big)_{\ell=1}^{\infty}=\big(\frac{\varphi(u^{({\ell)}})}{2^{\ell}-3^{h(\ell)}}\big)_{\ell=1}^{\infty}$=\\
\noindent $\big(-1,-1,-5,-\frac{23}{11},-\frac{85}{49},-5,-\frac{319}{115},\frac{319}{13},-\frac{1213}{217},-\frac{4151}{139},-\frac{853}{145},-\frac{47599}{11491},-\frac{47599}{3299},-\frac{159181}{26281},\frac{319}{13},\ldots\big)$.\\ 
\noindent Thus the subsequence $\big(C_{\mathbb{R}}\big(u^{({\ell)}}\big)\big)_{\ell\in\mathcal{L}}$ for $\mathcal{L}:=\big(q_k\big)_{k\; odd}$ converges to $\Phi^*_{\mathbb{R}}(1c_{\alpha})$ (e.g., for $k=11$ we have $C_{\mathbb{R}}\big(u^{({144)}}\big)\approx 10.396$).\footnote{The subsequence $\big(C_{\mathbb{R}}\big(u^{({\ell)}}\big)\big)_{\ell\in\mathcal{L}}$ for $\mathcal{L}:=\big(q_k\big)_{k\; \mathrm{even}}$ converges to $3\Phi^*_{\mathbb{R}}(1c_{\alpha})+\frac{1}{2}$. The $w_k$ always end with $0$, but if $k$ is even we replace this ultimate $0$ by $1$ to get a prefix of $1c_{\alpha}$. For example, we have $u^{(5)}=11011$, $\frac{\varphi(11011)}{2^5-3^4}=-\frac{85}{49}$, but $w_4=11010$, $\frac{\varphi(11010)}{2^5-3^3}=\frac{23}{5}$.}
\end{exa}

Let $\big(C_{\mathbb{R}}\big(u^{({\ell)}}\big)\big)_{\ell\in\mathcal{L}}$ be the sequence of Example \ref{ejgolden} where $\mathcal{L}:=\big(q_k\big)_{k\; odd}$. We calculate the inverse limit of its terms $\frac{\varphi(u^{(q_k)})}{2^{q_k}-3^{p_k}}$ for $k$ odd in $\mathbb{Z}_2$ and in $\mathbb{Z}_3$.

In $\mathbb{Z}_2$ we have\\
\noindent $\scriptstyle {\lim_{\leftarrow}C_{\mathbb{R}}(u^{(1)})=(C_{\mathbb{R}}(u^{(1)}) \bmod {2^n})_{n=1}^{\infty}=(-1 \bmod {2^n})_{n=1}^{\infty}=}\\
\scriptstyle {(1, 3, \;\;7, 15, 31, 63, 127, 255, 511, 1023, 2047, 4095, 8191, 16383, 32767, 65535, 131071, 262143, 524287, 1048575, 2097151, 4194303, \ldots)}$;

\noindent $\scriptstyle {\lim_{\leftarrow}C_{\mathbb{R}}(u^{(3)})=(C_{\mathbb{R}}(u^{(3)}) \bmod {2^n})_{n=1}^{\infty}=(-5 \bmod {2^n})_{n=1}^{\infty}=}\\
\scriptstyle {(1, 3, \;\;3, 11, 27, 59, 123, \;\;251, 507, 1019, 2043, 4091, 8187, 16379, 32763, 65531, 131067, 262139, 524283, 1048571,2097147, 4194299, \ldots)}$;
 
\noindent $\scriptstyle {\lim_{\leftarrow}C_{\mathbb{R}}(u^{(8)})=(C_{\mathbb{R}}(u^{(8)}) \bmod {2^n})_{n=1}^{\infty}=(\frac{319}{13} \bmod {2^n})_{n=1}^{\infty}=}\\
\scriptstyle {(1, 3, \;\;3, 11, 27, 59, 123, \;\;123, 379, 891, 1915, 1915, 1915, 10107, 10107, 10107, 10107, 141179, 403323, 403323, \;\;1451899, 3549051, \ldots)}$;

\noindent $\scriptstyle {\lim_{\leftarrow}C_{\mathbb{R}}(u^{(21)})=(C_{\mathbb{R}}(u^{(21)}) \bmod {2^n})_{n=1}^{\infty}=(\frac{5805215}{502829} \bmod {2^n})_{n=1}^{\infty}=}\\
\scriptstyle {(1, 3, \;\;3, 11, 27, 59, 123, \;\;123, 379, 891, 1915, 1915, 1915, 10107, 10107, 10107, 10107, 141179, 403323, 403323, \;\;403323, 2500475, \ldots)}$.

In $\mathbb{Z}_3$ we have\\
\noindent $\scriptstyle {\lim_{\leftarrow}C_{\mathbb{R}}(u^{(1)})=(C_{\mathbb{R}}(u^{(1)}) \bmod {3^n})_{n=1}^{\infty}=(-1 \bmod {3^n})_{n=1}^{\infty}=}\\
\scriptstyle {(2, 8, 26, 80, 242, 728, 2186, 6560, 19682, 59048, 177146, 531440, 1594322, 4782968, 14348906, 43046720, 129140162, 387420488, \ldots)}$;

\noindent $\scriptstyle {\lim_{\leftarrow}C_{\mathbb{R}}(u^{(3)})=(C_{\mathbb{R}}(u^{(3)}) \bmod {3^n})_{n=1}^{\infty}=(-5 \bmod {3^n})_{n=1}^{\infty}=}\\
\scriptstyle {(1, \;\;4, 22, 76, 238, 724, 2182, 6556, 19678, 59044, 177142, 531436, 1594318, 4782964, 14348902, 43046716, 129140158, 387420484, \ldots)}$;
 
\noindent $\scriptstyle {\lim_{\leftarrow}C_{\mathbb{R}}(u^{(8)})=(C_{\mathbb{R}}(u^{(8)}) \bmod {3^n})_{n=1}^{\infty}=(\frac{319}{13} \bmod {3^n})_{n=1}^{\infty}=}\\
\scriptstyle {(1, \;\;1, 10, 37, \;\;118, 361, 361, 2548, 9109, 9109, 68158, 245305, 245305, 1839628, 6622597, 6622597, 49669318, 178809481, \ldots)}$;

\noindent $\scriptstyle {\lim_{\leftarrow}C_{\mathbb{R}}(u^{(21)})=(C_{\mathbb{R}}(u^{(21)}) \bmod {3^n})_{n=1}^{\infty}=(\frac{5805215}{502829} \bmod {3^n})_{n=1}^{\infty}=}\\
\scriptstyle {(1, \;\;1, 10, 37, \;\;199, 442, 442, 4816, 11377, 11377, 129475, 306622, \;\;306622, 3495268, 8278237, 8278237, 94371679, 94371679,  \ldots)}$;

\noindent $\scriptstyle {\lim_{\leftarrow}C_{\mathbb{R}}(u^{(55)})=(C_{\mathbb{R}}(u^{(55)}) \bmod {3^n})_{n=1}^{\infty}=(\frac{204377029899901229}{19351615319297399} \bmod {3^n})_{n=1}^{\infty}=}\\
\scriptstyle {(1, \;\;1, 10, 37, \;\;199, 442, 442, 4816, 11377, 11377, 129475, 306622, \;\;1369504, 1369504, 10935442, 25284349, 68331070, 326611396, \ldots)}$.

The computation above shows that for all $k\geq 3$ odd
\begin{displaymath}
\big(X_{q_k})_{k}:=\Bigg(\bigg(C_{\mathbb{R}}\big(u^{(q_k)}\big)\bmod{2^n}\bigg)_{n=1}^{q_k+q_{k+1}-1}\Bigg)_{k=3,5,\ldots}\big(Y_{q_k})_k:=\Bigg(\bigg(C_{\mathbb{R}}\big(u^{(q_k)}\big)\bmod{3^n}\bigg)_{n=1}^{p_k-1}\Bigg)_{k=3,5,\ldots}
\end{displaymath}
are sequences of finite sequences $X_{q_k}$, $Y_{q_k}$ such that $X_{q_k}$ is a prefix of $X_{q_{k+2}}$ and $Y_{q_k}$ is a prefix of $Y_{q_{k+2}}$. This can be proved expanding $\varphi(u^{(q_k)})$ with formula (\ref{fi})
\begin{displaymath}
C_{\mathbb{R}}\big(u^{(q_k)}\big)=\frac{\varphi(u^{(q_k)})}{2^{q_k}-3^{p_k}}=\frac{\sum_{i=0}^{p_k-1}3^{p_k-1-i}\;   2^{d_i}}{2^{q_k}-3^{p_k}}.
\end{displaymath}

\noindent For $k\rightarrow\infty$ we obtain the inverse limits $X_{\infty}$ in $\mathbb{Z}_2$ and $Y_{\infty}$ in $\mathbb{Z}_3$. Thus\\
\noindent $X_{\infty}=1101111011100100011001111010001\cdots$, $Y_{\infty}=101121021021202112100201211012\cdots$.

We can verify that $X_{\infty}$ and $\Phi(1c_{\alpha})$ are the same $2$-adic integer. However, by Lemma \ref{necessary}, the limit $\Phi_{\mathbb{R}}(1c_{\alpha})$ does not exist since $\frac{2}{1+\sqrt5}<\frac{\ln(2)}{\ln(3)}$. Things are quite different in $\mathbb{Z}_3$ where $Y_{\infty}$ is the inverse limit of $C_{\mathbb{R}}(1c_{\alpha},\mathcal{L})$, that is,
$\lim_\leftarrow C_{\mathbb{R}}(1c_{\alpha},\mathcal{L})
=\lim\;\big(C_{\mathbb{R}}(u^{({\ell)}})\big)_{\ell\in\mathcal{L}}=\Phi^*_{\mathbb{R}}(1c_{\alpha})$ for $\mathcal{L}=\big(q_k\big)_{k>1,\; \mathrm{odd}}$. Therefore, $Y_{\infty}=\lim_{\leftarrow}\Phi^*_{\mathbb{R}}(1c_{\alpha})$ is an aperiodic $3$-adic integer, because the irrational $\Phi^*_{\mathbb{R}}(1c_{\alpha})$ does not admit an eventually periodic $3$-adic expansion. Hence the irrational $10.370127141747049162\cdots$ is represented by the aperiodic $3$-adic integer $1011210210212021121002012110122\cdots$.

If $\lim\;\big(C_{\mathbb{R}}(u^{({\ell)}})\big)_{\ell\in\mathcal{L}}$ exists as in Example \ref{ejgolden}, then $\lim\;(T_{u^{(\ell)}}(0))_{\ell\in\mathcal{L}}$ exists as well by Lemma \ref{CT}, since $\lim_{\ell\rightarrow \infty}(\frac{h}{\ell})=\frac{2}{1+\sqrt5}<\frac{\ln(2)}{\ln(3)}$ and $\lim\;\big(C_{\mathbb{R}}(u^{({\ell)}})\big)_{\ell\in\mathcal{L}}=\lim\;(T_{u^{(\ell)}}(0))_{\ell\in\mathcal{L}}$. Therefore, it is no surprise that the new sequence\\
\noindent $\big(T_{u^{(\ell)}}(0)\big)_{\ell\in\mathcal{L}}=(\frac{5}{2^3},\frac{319}{2^8},\frac{5805215}{2^{21}},\frac{204377029899901229}{2^{55}},\ldots)$ yields the same inverse limit in $\mathbb{Z}_3$:\\
\noindent $(1,1,10,37,199,442,442,4816,11377,11377,129475,306622,1369504,1369504,\ldots)$.

In $\mathbb{Z}_2$ it holds that $\big|\Phi(u)-C(u)\big|_2=\Big| \frac{-\varphi(u)\cdot 2^\ell}{3^h(2^\ell-3^h)} \Big|_2\leq \frac{1}{2^\ell}$ where $u\in\mathbb{U}$ is a prefix of $v$ with length $\ell$. Hence $\Phi(u)$ and $C(u)$ coincide at least in the first $\ell$ digits. But, $\Phi(u)$ and $\Phi(v)$ coincide also in the first $\ell$ digits since $\big|\Phi(v)-\Phi(u)\big|_2=\big|v-u\big|_2$. It follows that $C(u)$ and $\Phi(v)$ coincide in the first $\ell$ digits. Since for any $v\in\mathbb{Z}_2$ exists $\Phi(v)$ we have
\begin{equation}\label{ILC}
\Phi(v)=\bigg(C_{\mathbb{R}}\big(u^{({\ell)}}\big)\bmod{2^{\ell}}\bigg)_{\ell=1}^{\infty}, \textrm{ or }\;\;\Phi(v)=\lim_{\ell\rightarrow \infty}\;C(u^{(\ell)})\;\textrm{ with }C(u^{(\ell)})\in\mathbb{Z}_2\;.
\end{equation}
The latter relation can be used to approximate $\Phi(v)$ for a given $v$. By choosing a prefix $u$ with length $\ell$, by calculating the rational $C(u)$ and converting it into a 2-adic integer $x$, then $\Phi(v)$ and $x$ will coincide in the first $\ell$ digits. The magnitude $\Phi_{\mathbb{R}}(u)$ evaluated in $\mathbb{Z}_2$ will yield the same $\ell$ digits. This will work in all cases, even when $|\Phi_{\mathbb{R}}(u)-C_{\mathbb{R}}(u)|$ does not converge, or when the limits $\Phi_{\mathbb{R}}(v)$, $C_{\mathbb{R}}(v)$ do not exist. For instance, such a case is $1c_{\alpha}$ if $\alpha=\frac{\ln(2)}{\ln(3)}$. We obtain $\Phi(1c_{\alpha})=110111111101101001111101100100111101010100000010100000001\cdots$, although the limits $\Phi_{\mathbb{R}}(v)$, $C_{\mathbb{R}}(v)$ do not exist.

In $\mathbb{Z}_3$ it holds that $\big|T_u(0)-C(u)\big|_3=\Big| \frac{-\varphi(u)\cdot 3^h}{2^\ell(2^\ell-3^h)} \Big|_3= \frac{1}{3^h}$ where $u\in\mathbb{U}$ is a prefix of $v$ with length $\ell$ and height $h$. Hence $T_u(0)$ and $C(u)$ coincide exactly in the first $h$ $3$-adic digits. However, $T_{u1}(0)=\frac{3T_u(0)+1}{2}$ and $T_{u0}(0)=\frac{T_u(0)}{2}$, so that with any additional bit we have a division by $2$ in $\mathbb{Z}_3$ altering the first $h$ digits, and $T_v(0)$ is undefined for any infinite word $v$ without suffix $0^{\infty}$. If $u^{(\ell)}$ is a prefix of $v$ with length $\ell$ and height $h(\ell)$, then we have
\begin{equation}\label{ILCT}
\bigg(T_{u^{(\ell)}}(0) \bmod{3^n}\bigg)_{n=1}^{h(\ell)}=\bigg(C_{\mathbb{R}}(u^{(\ell)})\bmod{3^n}\bigg)_{n=1}^{h(\ell)}
\end{equation}
for the first $h(\ell)$ integers of the inverse limit in $\mathbb{Z}_3$. For instance, in the case $1c_{\alpha}$ with $\alpha=\frac{\ln(2)}{\ln(3)}$ we obtain a sequence of $3$-adic output for increasing $\ell$ in Table \ref{tab:3adic}.

\begin{table}[htb]
\centering
\begin{tabular}{llllll}
$\ell$&$1c_{\alpha}$&$T_{u^{(\ell)}}(0)$&$C_{\mathbb{R}}(u^{(\ell)})$&$\bmod\;{3^n}$&$3$-adic\\
\hline
1&1&$\frac{1}{2}$&$-1$&$(2)$&2\\
2&11&$\frac{5}{4}$&$-1$&$(2,8)$&22\\
3&110&$\frac{5}{8}$&$-5$&$(1,4)$&11\\
4&1101&$\frac{23}{16}$&$-\frac{23}{11}$&$(2,2,20)$&202\\
5&11011&$\frac{85}{32}$&$-\frac{85}{49}$&$(2,8,17,71)$&2212\\
6&110110&$\frac{85}{64}$&$-5$&$(1,4,22,76)$&1122\\
7&1101101&$\frac{319}{128}$&$-\frac{319}{115}$&$(2,2,20,74,236)$&20222\\
8&11011011&$\frac{1085}{256}$&$-\frac{1085}{473}$&$(2,8,17,71,233,719)$&221222\\
9&110110110&$\frac{1085}{512}$&$-5$&$(1,4,22,76,238,724)$&112222\\
10&1101101101&$\frac{3767}{1024}$&$-\frac{3767}{1163}$&$(2,2,20,74,236,722,2180)$&2022222\\
11&11011011010&$\frac{3767}{2048}$&$-\frac{3767}{139}$&$(1,1,10,37,118,361,1090)$&1011111\\
\end{tabular}
\caption{$3$-adic output}
\label{tab:3adic}
\end{table}
Note that in Table \ref{tab:3adic} the first digit minus $1$ of each $3$-adic integer gives the digits of the parity vector $1c_{\alpha}$. This is true if and only if the parity vector does not contain two or more consecutive $0$'s.

The infinite sequence of $3$-adic integers in Table \ref{tab:3adic} is not convergent in $\mathbb{Z}_3$. However, we have convergence in $\mathbb{Z}_3$ for the subsequence $\mathcal{L}=\big(q_k\big)_{k>1  \;\mathrm{odd}}$ and also for $\mathcal{L}'=\big(q_k\big)_{k>1  \;\mathrm{even}}$ being $\big(\frac{p_k}{q_k}\big)_{k=0}^{\infty}=(0,1,\frac{1}{2},\frac{2}{3},\frac{5}{8},\frac{12}{19},\frac{41}{65},\frac{53}{84},\frac{306}{485},\frac{665}{1054},\ldots)$ the convergents of $\frac{\ln(2)}{\ln(3)}$. Since $1c_{\alpha}$ is a Sturmian word, we get the inverse limits with
\begin{displaymath}
\Bigg(\bigg(C_{\mathbb{R}}\big(u^{(q_k)}\big)\bmod{3^n}\bigg)_{n=1}^{p_k-1}\Bigg)_{k=3,5,7,\ldots}\textrm{ and }\;\Bigg(\bigg(C_{\mathbb{R}}\big(u^{(q_k)}\big)\bmod{3^n}\bigg)_{n=1}^{p_k+1}\Bigg)_{k=2,4,6,\ldots}.
\end{displaymath}
\noindent The corresponding computation yields\\
\noindent $\lim\;\big(C(u^{(\ell)})\big)_{\ell\in\mathcal{L}}=101111011011011022120100121111011121211122021111122112\cdots\;$ and\\
\noindent
$\lim\;\big(C(u^{(\ell)})\big)_{\ell\in\mathcal{L}'}=221222212212212210101221120000022220010001120000001100\cdots$,\\
\noindent where $C(u^{(\ell)})\in\mathbb{Z}_3\;$. 

\noindent These two limits are possible $3$-adic outputs for the given aperiodic word $v:=1c_{\ln(2)/\ln(3)}$ equipped with an appropriate $\mathcal{L}$ (respectively $\mathcal{L}'$), while the $2$-adic $\Phi(v)$ is its unique input. In the same way using $q_k$ we look for possible outputs for the words $S(v)=v_1v_2v_3\cdots, S^2(v)=v_2v_3\cdots,S^3(v)=v_3v_4\cdots$, etc.; each of them with two possible outputs, namely,\\ 
\noindent $20020200200200200\cdots$ and $\;11222220112220112\cdots\;$ for $\;S(v)$;\\
\noindent $22112121121121121\cdots$ and $\;20222222102222102\cdots\;$ for $\;S^2(v)$;\\
\noindent $11201220220220220\cdots$ and $\;10111111212222001\cdots\;$ for $\;S^3(v)$.\\
\noindent Since $(\Phi(v_0v_1v_2\cdots),\Phi(v_1v_2v_3\cdots),\Phi(v_2v_3v_4\cdots),\ldots)$ is the trajectory of $v$, these new outputs are also outputs of $v$. It follows that there are infinitely many $3$-adic outputs for the given aperiodic word $v$ equipped with a set of appropriated $\mathcal{L}$, while the $2$-adic $\Phi(v)$ is its unique input.

We propose the following model which relates $\mathbb{Z}_2$, $\mathbb{Z}_3$ and $\mathbb{R}$ as a challenge for further research. Let $v$ be an aperiodic word and $\mathcal{L}$ a subsequence (proper or improper) of the positive integers $(1,2,3,\dots)$ such that the $3$-adic limit $T_{v,\mathcal{L}}(0)=\lim\;\big(C(u^{(\ell)})\big)_{\ell\in\mathcal{L}}$ exists. We conjecture that if such an $\mathcal{L}$ always exists, then there are infinitely many of them. We interpret the pair $(v,\mathcal{L})$ as a machine with a parity vector $v$ and an instruction $\mathcal{L}$. The input is given in $Z_2$ by the unique $\Phi(v)$, and the outputs are given in $Z_3$ by $T_{v,\mathcal{L}}(0)$. In $\mathbb{R}$ we have access to $\Phi(v)$ and $T_{v,\mathcal{L}}(0)$ by means of
\begin{equation}\label{machine}
\Phi(v)=\bigg(C_{\mathbb{R}}\big(u^{({\ell)}}\big)\bmod{2^{\ell}}\bigg)_{\ell=1}^{\infty}\;\textrm{ and }\;T_{v,\mathcal{L}}(0)=\Bigg(\bigg(C_{\mathbb{R}}(u^{(\ell)})\bmod{3^n}\bigg)_{n=1}^{h(\ell)}\Bigg)_{\ell\in\mathcal{L}}\;.
\end{equation}
$\Phi(v)$ is an eventually periodic word with $1$'s and $0$'s if and only if $T_{v,\mathcal{L}}(0)$ is eventually periodic with digits $\{0,1,2\}$. So, the $2$-adic $\Phi(v)$ is aperiodic if and only if the $3$-adic $T_{v,\mathcal{L}}(0)$ is aperiodic. If $\Phi_{\mathbb{R}}(v)$ exists, then $C_{\mathbb{R}}(v,\mathcal{L})=\Phi_{\mathbb{R}}(v)$ is a negative irrational number which admits an aperiodic $2$-adic expansion. If $\Phi_{\mathbb{R}}(v)$ does not exist, but $C_{\mathbb{R}}(v,\mathcal{L})$ does, then $C_{\mathbb{R}}(v,\mathcal{L})=T_{v,\mathcal{L}}(0)$ is a non-negative real number. If this number is irrational, then it admits an aperiodic $3$-adic expansion.

\section{Further Results and Examples}
The following results are enumerated in a different order from their proof in the text. They form part of a chain of $58$ Lemmas which begins in Section \ref{proof}. In each result we indicate with brackets where the proof can be found.
 
\begin{figure}[htbp]
\begin{center}
\epsfxsize=6in
\epsfbox{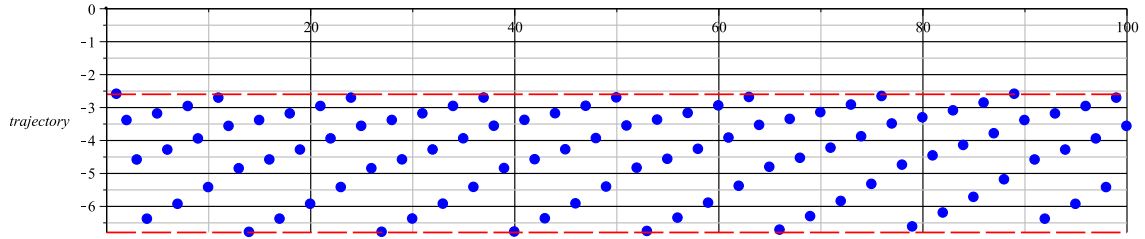}
\end{center}
\caption{Trajectory of $\Phi_{\mathbb{R}}(1c_{\;\ln(2)})\approx -2.5970181822928929618$}
\label{abov}
\end{figure}

\begin{res}\label{eins} (Figure \ref{abov}).
If $\alpha$ is irrational and $1>\alpha>\ln(2)/\ln(3)$, then $\Phi_{\mathbb{R}}(1c_{\alpha})$ is the limit upper and $\Phi_{\mathbb{R}}(0c_{\alpha})=3\Phi_{\mathbb{R}}(1c_{\alpha})+1$ is the lower limit of the trajectory of $\Phi_{\mathbb{R}}(1c_{\alpha})$. Both limits are accumulation points of the orbit $\mathcal{O}(\Phi_{\mathbb{R}}(1c_{\alpha}))$. Hence $\Phi_{\mathbb{R}}(1c_{\alpha})\notin\mathbb{Q}_{odd}$. Moreover, \textit{all} points of the trajectory are accumulation points of $\mathcal{O}(\Phi_{\mathbb{R}}(1c_{\alpha}))$. The number $\Phi^*_{\mathbb{R}}(1c_{\alpha})$ does not exist $(=+\infty)$. \hfill \textup{[Lemma \ref{uplo1lim}]}$\blacktriangleleft$
\end{res}

\begin{figure}[htbp]
\begin{center}
\epsfxsize=6in
\epsfbox{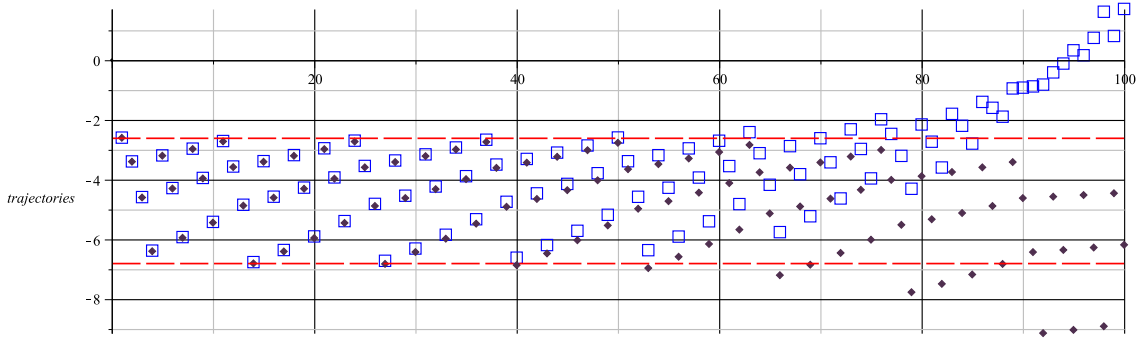}
\end{center}
\caption{Trajectories $\mathcal{T}_{1c_{\;\ln(2)}}(-2.593)$, $\mathcal{T}_{1c_{\;\ln(2)}}(-2.599)$ for $\alpha=\ln(2)\approx 0.693147$}
\label{abov2}
\end{figure}

\begin{figure}[htbp]
\begin{center}
\epsfxsize=6in
\epsfbox{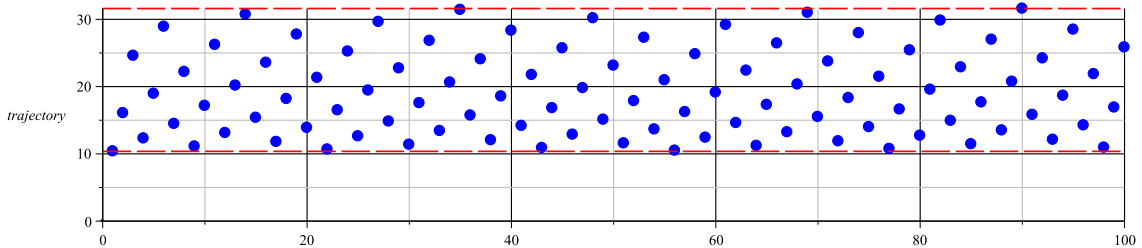}
\end{center}
\caption{Trajectory of $\Phi^*_{\mathbb{R}}(1c_{2/(1+\sqrt 5)})\approx 10.370127141747049162$}
\label{belo}
\end{figure}

\begin{figure}[htbp]
\begin{center}
\epsfxsize=6in
\epsfbox{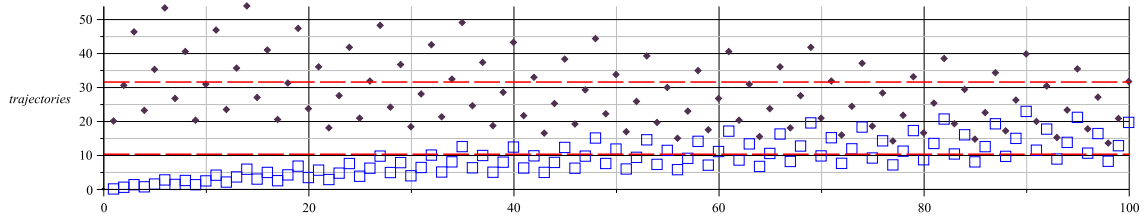}
\end{center}
\caption{Trajectories $\mathcal{T}_{1c_{\alpha}}(20)$ and $\mathcal{T}_{1c_{\alpha}}(0)$ for $\alpha=2/(1+\sqrt 5)$}
\label{belo2}
\end{figure}

What happens if we choose a real number $x\neq\Phi_{\mathbb{R}}(1c_{\alpha})$?  Figure \ref{abov2} illustrates the situation. 
For $x<\Phi_{\mathbb{R}}(1c_{\alpha})$ we have $\mathcal{T}_{1c_\alpha}(x)\rightarrow -\infty$  without an upper limit of the pseudo trajectory, and for $x>\Phi_{\mathbb{R}}(1c_{\alpha})$ we have $\mathcal{T}_{1c_\alpha}(x)\rightarrow \infty$  without a lower limit (Lemmas \ref{x+y} and \ref{2/3}). But for $x=\Phi_{\mathbb{R}}(1c_{\alpha})$ we have by definition $\mathcal{T}_{1c_\alpha}(x)=\mathcal{T}(x)$ with the corresponding upper and lower limits. Thus $\Phi_{\mathbb{R}}(1c_{\alpha})$ must be something special. Indeed, $\Phi_{\mathbb{R}}(1c_{\alpha})$ is the real number that fits as input for the aperiodic $1c_\alpha$. We say that $1c_\alpha$ is the parity vector of the irrational $\Phi_{\mathbb{R}}(1c_{\alpha})$.
\begin{res}\label{zwei} (Figures \ref{belo} and \ref{belo2}).
If $\alpha$ is irrational and $0<\alpha<\ln(2)/\ln(3)$, then $\Phi_{\mathbb{R}}(1c_{\alpha})$ does not exist $(=-\infty)$. But then $\Phi^*_{\mathbb{R}}(1c_{\alpha})$ is the lower limit and $\Phi^*_{\mathbb{R}}(0c_{\alpha})=3\Phi^*_{\mathbb{R}}(1c_{\alpha})+1/2$ is the upper limit of the pseudo trajectory $\mathcal{T}_{1c_\alpha}(x)$ of \textit{any} given $x\in\mathbb{R}$. Both limits are accumulation points of the orbit $\mathcal{O}(x)$. Therefore, no $\zeta\in\mathbb{Q}_{odd}$ has $1c_{\alpha}$ as parity vector.\\ 
\noindent If we choose $x=\Phi^*_{\mathbb{R}}(1c_{\alpha})$, then the whole trajectory lies within the semi-open interval
$[ \Phi^*_{\mathbb{R}}(1c_{\alpha}),3\Phi^*_{\mathbb{R}}(1c_{\alpha})+1/2 )$.
\hfill \textup{[Lemmas \ref{Tu0}, \ref{limits}]}$\blacktriangleleft$
\end{res}

\begin{figure}[htbp]
\begin{center}
\epsfxsize=5in
\epsfbox{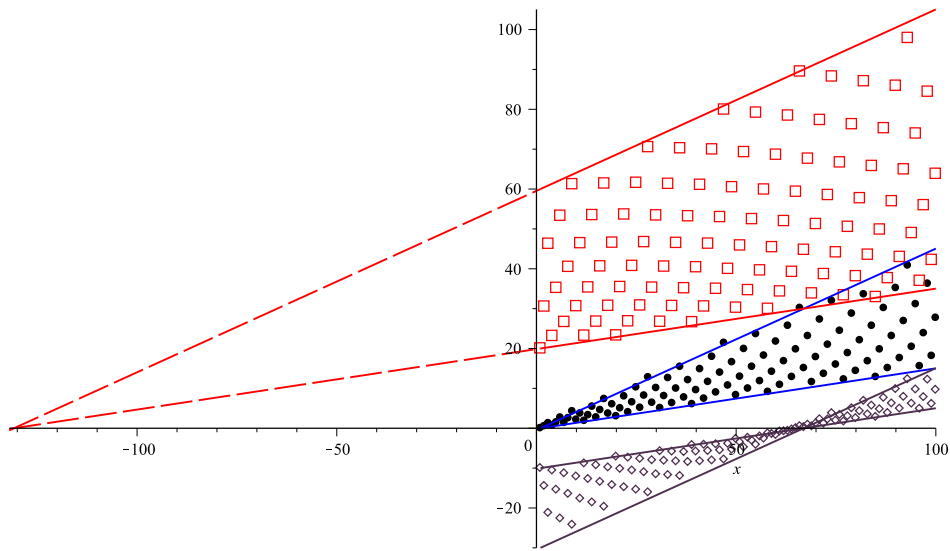}
\end{center}
\caption{Trajectories $\mathcal{T}_{1c_{\alpha}}(20)$, $\mathcal{T}_{1c_{\alpha}}(0)$ and $\mathcal{T}_{1c_{
\alpha}}(-10)$ for $\alpha=\ln(2)/\ln(3)$}
\label{con}
\end{figure}
\begin{res}\label{drei} (Figure \ref{con}).
If $\alpha=\ln(2)/\ln(3)$, then $\Phi_{\mathbb{R}}(1c_{\alpha})$ $(=-\infty)$ and $\Phi^*_{\mathbb{R}}(1c_{\alpha})$ $(=+\infty)$ do not exist. The pseudo trajectory $\mathcal{T}_{1c_{\alpha}}(0)$ lies entirely in a (slightly enlarged) cone bounded by two lines with slopes $\frac{1}{6\ln(3)}$ and $\frac{1}{2\ln(3)}$. The pseudo trajectory $\mathcal{T}_{1c_{\alpha}}(x)$ of any other $x\in\mathbb{R}$ (in the figure $x=20$ and $x=-10$) lies in a horizontally transferred cone with the same slopes. Since there are only finitely many points below any horizontal line, orbits $\mathcal{O}(0)$ and $\mathcal{O}(x)$ have no accumulation point.\hfill \textup{[Section \ref{pseudo}]}$\blacktriangleleft$
\end{res}

\begin{res}\label{1}
Let $v$ be an aperiodic word.
\begin{displaymath}
\textrm{If}\quad\underline{\lim}\;\bigg(\frac{h}{\ell}\bigg)_{\ell=1}^\infty>\frac{\ln(2)}{\ln(3)}\;,\quad\textrm{then limit }\Phi_{\mathbb{R}}(v)<-1\;\textrm{exists},\; \textrm{ and }\; \Phi_{\mathbb{R}}(v)\notin\mathbb{Q}_{odd}\textrm{ holds }.\\[-20pt]
\end{displaymath}
\end{res}
\hfill \textup{[Lemmas \ref{limarriba}, \ref{orbitarriba}]}$\blacktriangleleft$

\noindent Example \ref{vier}  below illustrates this result. For any real number $\alpha\in (\frac{\ln(2)
}{\ln(3)},1)$ we construct aperiodic words built up by factors $0^n1^m$ to increase $n$ 
such that $\lim_{\ell\rightarrow\infty}\big(\frac{h}{\ell}\big)=\alpha$ using 
the pseudo code: 
\begin{tabular}{ll}
1.&input $top:=$ number of factors $0^n1^m$, parameter $\alpha$\\
2.&$\ell:=0$, $h:=0$, $v:=$ empty list\\ 
3.&for $j=1, 2, 3,\ldots$ to $top$ do\\
4.&$n:=j$\\ 
5.&$m:=\lfloor\frac{(\ell+n)\alpha-h}{1-\alpha}\rfloor$\\
6.&redefine $v:=v0^n1^m$\\ 
7.&redefine $\ell:=\ell+n+m$, $h:=h+m$\\
8.&end loop\\
9.&output $v$\\
\end{tabular}

\noindent If the slope is $\frac{h}{\ell}$ after having executed the loop several times, then we adjust the slope in the next turn by $\frac{h+m}{\ell+n+m}\approx\alpha$ and hence $m\approx\frac{(\ell+n)\alpha-h}{1-\alpha}$.

\begin{figure}[htbp]
\begin{center}
\epsfxsize=6.5in
\epsfbox{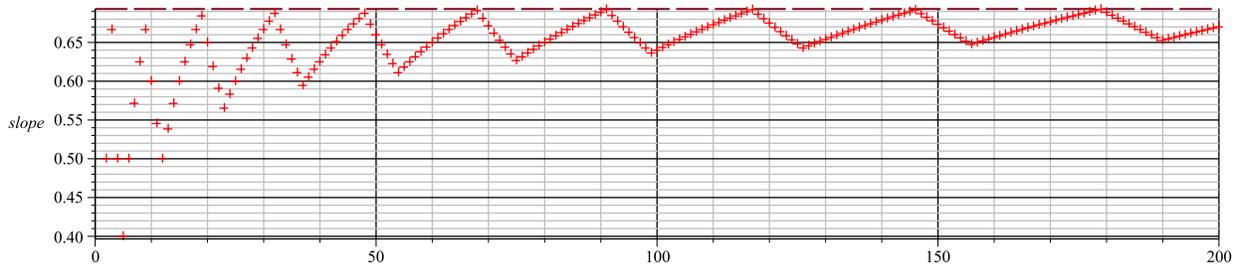}
\end{center}
\caption{Slope $(\frac{h}{\ell})_{\ell=2}^{200},\;\alpha=\ln(2)\approx 0.693147$}
\label{ln2}
\end{figure}

\begin{figure}[htbp]
\begin{center}
\epsfxsize=5in
\epsfbox{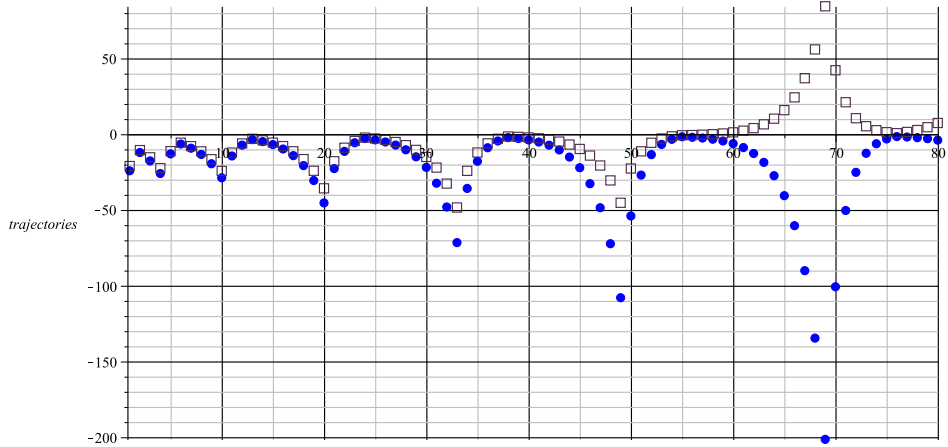}
\end{center}
\caption{Trajectory of $\Phi_{\mathbb{R}}(v)\approx -24.17309$ and $\mathcal{T}_v(-21)$}
\label{stand}
\end{figure}

\begin{exa}\label{vier} (Figures \ref{ln2} and \ref{stand}).
Let $1>\alpha> \frac{\ln(2)}{\ln(3)}$, $v=\big(0^n1^m\big)_{n=1,2,3,\ldots}$ and $m=\lfloor\frac{(\ell+n)\alpha-h}{1-\alpha}\rfloor$.
For example, if $\alpha=\ln(2)$, then $v=01^2\; 0^21^4\; 0^31^7\; 0^41^9\; 0^51^{11}\; 0^61^{14}\; 0^71^{16}\; 0^81^{18}\; 0^91^{20}\; 0^{10}1^{23}\cdots$, $\Phi_{\mathbb{R}}(v)\approx -24.17309$, $\lim_{\ell\rightarrow\infty}\big(\frac{h}{\ell}\big)=\alpha$. Since $\underline{\lim}\big(\frac{h}{\ell}\big)=\alpha>\frac{\ln(2)}{\ln(3)}$, orbit $\mathcal{O}(\Phi_{\mathbb{R}}(v))$ has accumulation points. Thus $\Phi_{\mathbb{R}}(v)\notin\mathbb{Q}_{odd}$.
\end{exa}


If $\underline{\lim}\;\big(\frac{h}{\ell}\big)_{\ell=1}^\infty>\frac{\ln(2)}{\ln(3)}$, then limit $\Phi_{\mathbb{R}}(v)$ exists by Result \ref{1}. If $\underline{\lim}\;\big(\frac{h}{\ell}\big)_{\ell=1}^\infty=\frac{\ln(2)}{\ln(3)}$, then limit $\Phi_{\mathbb{R}}(v)$ may not exist, although condition $\big(\frac{2^{\ell}}{3^{h}}\big)_{\ell=1}^{\infty}\rightarrow 0$ is fulfilled. We have the following criterion.

\begin{res}\label{final}
Let $v$ be an aperiodic word such that 
$\underline{\lim}\; \big(\frac{h}{\ell}\big)_{\ell=1}^{\infty}=\frac{\ln(2)}{\ln(3)}=\alpha$ and $1c_v$ its associated word (see Definition \ref{def1cv}, Section \ref{unocv}).

\noindent If $\big( \frac{2^{\ell}}{3^{\lceil\ell\alpha\rceil+n_{\ell}}}\big)_{\ell=1}^{\infty}\rightarrow 0$, where $n_{\ell}:=h(\ell)-\lceil \ell\alpha\rceil$ is the relative height of $v$ over the Sturmian word $1c_{\alpha}$ at position $\ell$, then $\ell_j:=\max\{\ell\;|\;n_{\ell}=j\}$ is well defined for $j=0,1,2,\ldots$, and the following holds: 
\begin{displaymath}
\textrm{If }\quad\overline{\lim}\;\bigg(\frac{\ell_{j+2}-\ell_{j+1}}{\ell_{j+1}-\ell_j}\bigg)_{j=1}^{\infty}< 3\quad\textrm{, then limits } \Phi_{\mathbb{R}}(1c_v)\leq\Phi_{\mathbb{R}}(v)\;\textrm{exist}.
\end{displaymath}
\begin{displaymath}
\textrm{If }\quad\overline{\lim}\;\bigg(\frac{\ell_{j+2}-\ell_{j+1}}{\ell_{j+1}-\ell_j}\bigg)_{j=1}^{\infty}> 3\quad\textrm{, then limit } \Phi_{\mathbb{R}}(1c_v) \textrm{ does not exist.} 
\end{displaymath}
\noindent If $\big( \frac{2^{\ell}}{3^{\lceil\ell\alpha\rceil+n_{\ell}}}\big)_{\ell=1}^{\infty}$ is not convergent to $0$, then limits 
$\Phi_{\mathbb{R}}(1c_v)$, $\Phi_{\mathbb{R}}(v)$ do not exist.
\end{res}
\hfill \textup{[Lemmas \ref{max}, \ref{suff}]}$\blacktriangleleft$

In Example \ref{vier} the number of consecutive $0$'s increases by $1$. How fast can the number of consecutive $0$'s grow such that $\Phi_{\mathbb{R}}(v)$ will exists? The following Result \ref{funf} gives the answer. 

\begin{res}\label{funf}
Let $1>\alpha> \frac{\ln(2)}{\ln(3)}$, $v=\big(0^{f(n)}1^m\big)_{n=1,2,3,\ldots}$, $m=\lfloor\frac{(\ell+f(n))\alpha-h}{1-\alpha}\rfloor$, and let $f(n)\in\emph{\textbf{O}}\Big(\bigg(\frac{\alpha(\frac{\ln(3)}{\ln(2)}-1)}{1-\alpha}\bigg)^n\Big)$ be an increasing function with range in $\mathbb{N}$. Then limit $\Phi_{\mathbb{R}}(v)$ exists. 
\end{res}
\hfill \textup{[Section \ref{geom}]}$\blacktriangleleft$

\noindent If function $f(n)$ is growing faster than the exponential function in Result \ref{funf}, then $\Phi_{\mathbb{R}}(v)$ diverges to $-\infty$. Example \ref{bigO} illustrates Result \ref{funf}. 

\begin{figure}[htbp]
\begin{center}
\epsfxsize=6in
\epsfbox{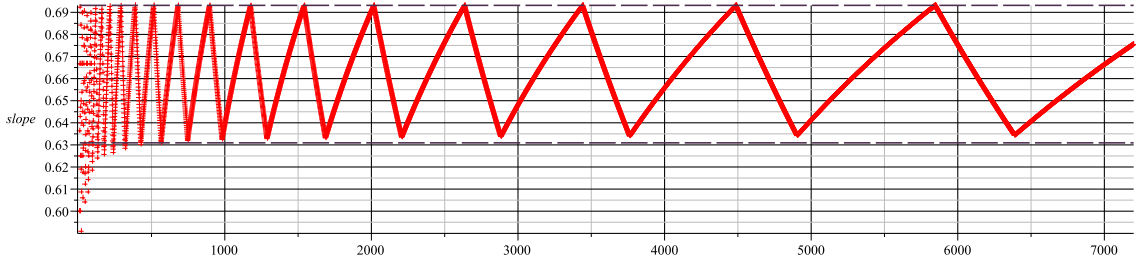}
\end{center}
\caption{Slope $(\frac{h}{\ell})_{\ell=10}^{7200}$, $\overline{\lim}\;\big(\frac{h}{\ell}\big)=\ln(2)\approx 0.693$, $\underline{\lim}\;\big(\frac{h}{\ell}\big)\approx 0.634 > \frac{\ln(2)}{\ln(3)}\approx 0.6309$}
\label{standO}
\end{figure}

\begin{figure}[htbp]
\begin{center}
\epsfxsize=5in
\epsfbox{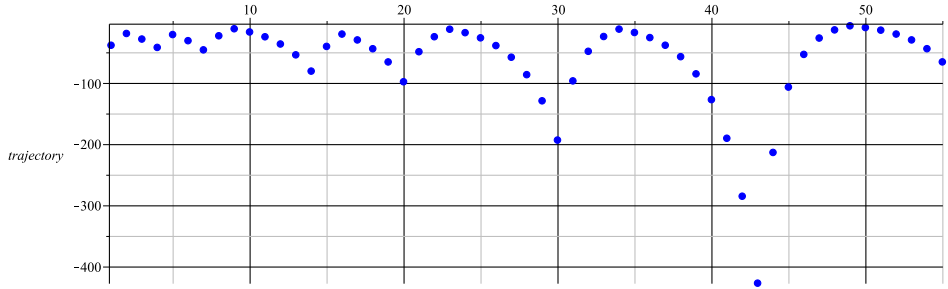}
\end{center}
\caption{Trajectory of $\Phi_{\mathbb{R}}(v)\approx -38.43026$}
\label{standOtraj}
\end{figure}

\begin{exa}\label{bigO} (Figures \ref{standO} and \ref{standOtraj}).
Let $\alpha=\ln(2)$, then we have the constraint $f(n)\in\emph{\textbf{O}}\big(\big(1.321366\ldots\big)^n\big)$. Polynomials of any degree with positive integers as coefficients satisfy this constraint, but, for  instance, $2^n$ does not. We choose the exponential function $f(n)=\lfloor 1.3^n\rfloor$, replacing $n:=j$ by $n:=\lfloor 1.3^j \rfloor$ in the pseudo code. Then we have\\
$v=01^2\; 01^2\; 0^21^5\; 0^21^4\; 0^31^{7}\; 0^41^{9}\; 0^61^{13}\; 0^81^{18}\; 0^{10}1^{23}\; 0^{13}1^{29}\; 0^{17}1^{39}\cdots$, $\Phi_{\mathbb{R}}(v)\approx -38.43026$, $\overline{\lim}\;\big(\frac{h}{\ell}\big)=\alpha=\ln(2)$. Since $\underline{\lim}\;\big(\frac{h}{\ell}\big)\approx 0.634 > \frac{\ln(2)}{\ln(3)}$, orbit $\mathcal{O}(\Phi_{\mathbb{R}}(v))$ has accumulation points (being $-1$ one of them). Thus $\Phi_{\mathbb{R}}(v)\notin\mathbb{Q}_{odd}$ by Result \ref{1}.
\end{exa}


\begin{figure}[htbp]
\begin{center}
\epsfxsize=6.5in
\epsfbox{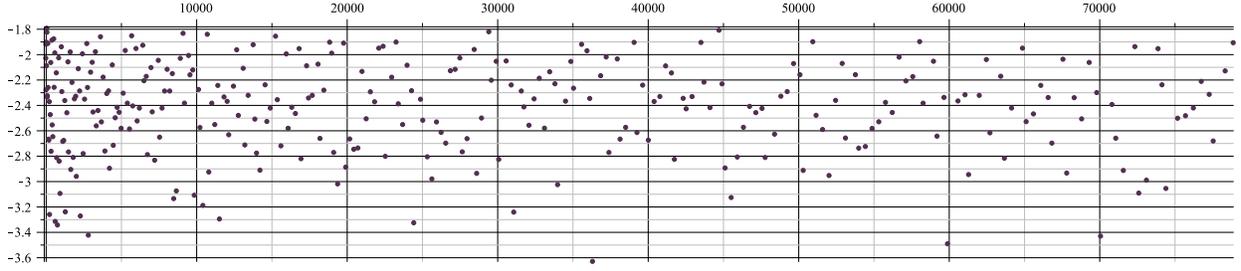}
\end{center}
\caption{Sub-sequence $(\ell= 9,17,25,38,52,74,\ldots)$ of the Trajectory of $\Phi_{\mathbb{R}}(v)\approx -2.29230$}
\label{starts}
\end{figure}

\begin{figure}[htbp]
\begin{center}
\epsfxsize=6in
\epsfbox{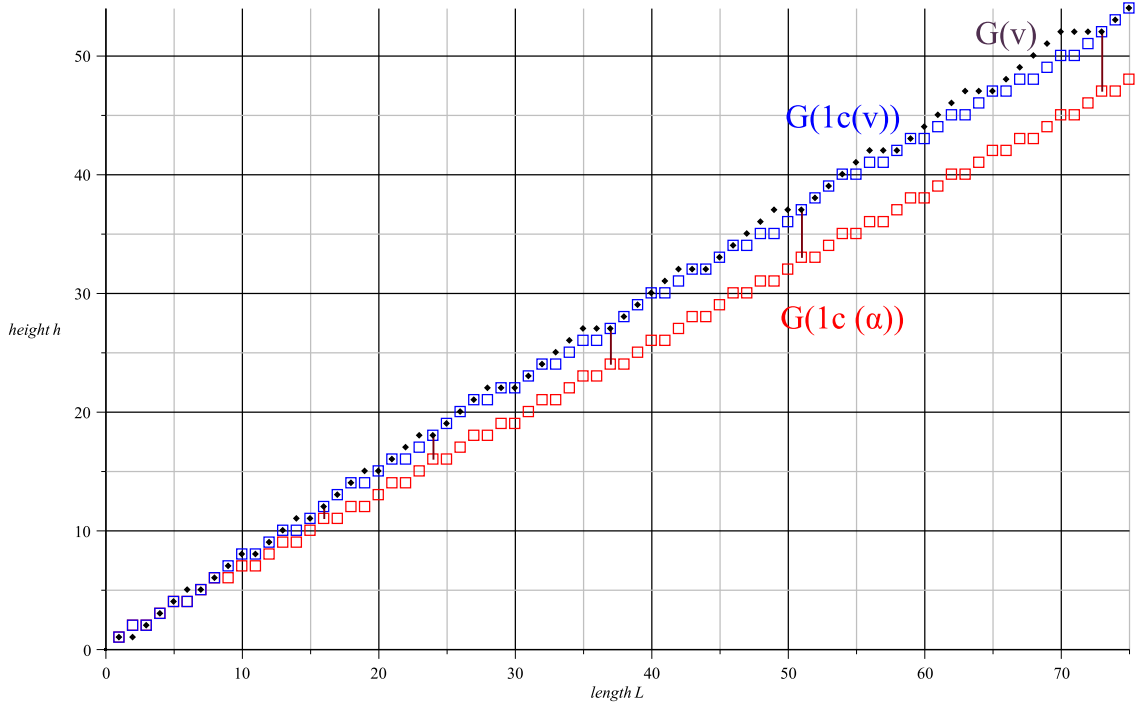}
\end{center}
\caption{Words $v$, $1c_v$ and $1c\;_{\ln(2)/\ln(3)}$ of Examples \ref{On} and \ref{acht}.}
\label{pal}
\end{figure}

\begin{exa}\label{On} (Figure \ref{starts}, also Figure \ref{pal}).
We choose $f(n)=n-\lfloor\frac{5}{6}n\rfloor$, replacing $v:=v0^n1^m$ by $v:=v1^m0^{f(n)}$ and $\ell:=\ell+n+m$ by $\ell:=\ell+f(n)+m$ in the pseudo code. Then\\ 
$v=10\; 1^40\; 1^30\; 1^30\; 1^{4}0\; 1^{3}0\; 1^{4}0^2\; (1^{5}0^2)^5\;  1^{5}0^{3}\; (1^{7}0^{3})^5\; 1^{7}0^{4}\; 1^{8}0^{4}\; 1^{9}0^{4}\; 1^{8}0^{4}\; 1^{9}0^{4}\; 1^{8}0^{4}\; 1^{9}0^{5}\cdots$, $\Phi_{\mathbb{R}}(v)\approx -2.29230$ and $\lim_{\ell\rightarrow\infty}\big(\frac{h}{\ell}\big)=\frac{\ln(2)}{\ln(3)}$. Figure \ref{starts} suggests that infinitely many points of trajectory $\mathcal{T}(\Phi_{\mathbb{R}}(v))$ fall into range $(2\Phi_{\mathbb{R}}(v)+1,-1)$. It seems that $\Phi_{\mathbb{R}}(v)$ is an accumulation point. We have no proof. However, there is strong evidence that $\Phi_{\mathbb{R}}(v)$ is an irrational number (Section \ref{end}, at the end).
\end{exa}

\begin{res}\label{sieben} 
Let $\alpha=\frac{\ln(2)}{\ln(3)}$, $v=\big(1^m0^{f(n)}\big)_{n=1,2,3,\ldots}$, $m=\lfloor\frac{(\ell+f(n))\alpha-h}{1-\alpha}\rfloor$, and let $f(n)\in\emph{\textbf{O}}\big(n\big)$ be an increasing function with range in $\mathbb{N}$. Then limit $\Phi_{\mathbb{R}}(v)$ exists.\hfill \textup{[Section \ref{end}]}$\blacktriangleleft$
\end{res}

\begin{figure}[htbp]
\begin{center}
\epsfxsize=6in
\epsfbox{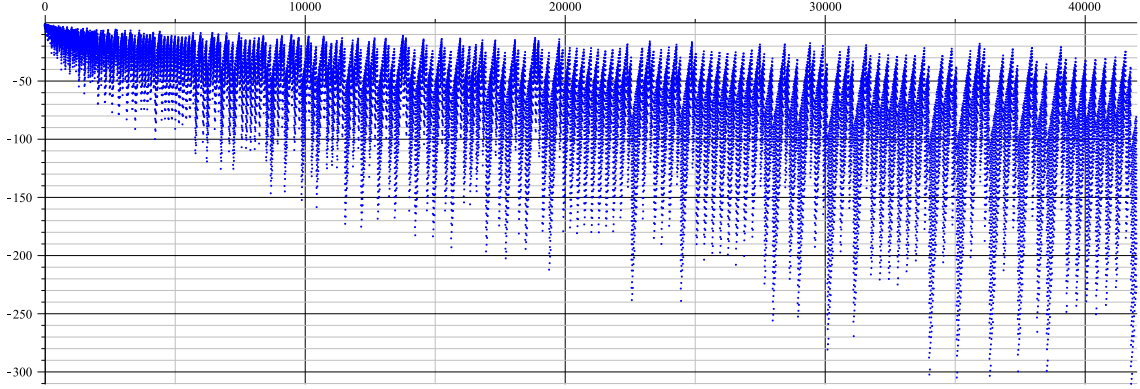}
\end{center}
\caption{Trajectory of $\Phi_{\mathbb{R}}(1c_v)\approx -2.33103$}
\label{pei}
\end{figure}

\begin{exa}\label{acht} 
Let $1c_v$ be the associated word to the $v$ of Example \ref{On}. This word coincides piecewise with $1c_{\alpha}$ but has 1's instead of 0's at certain positions which depend on $v$ as follows:

\noindent Let $G(v)=\big(\ell,h\big)_{\ell=1}^\infty$, $G(1c_v)=\big(\ell,h'\big)_{\ell=1}^\infty$ and $G(1c_{\alpha})=\big(\ell,\lceil\ell\alpha\rceil\big)_{\ell=1}^\infty$ be the geometric representations. At positions $\big(1+\max\{\ell\;|\;h=\lceil\ell\alpha\rceil+j\}\big)_{j=0}^{\infty}=(9,17,25,38,52,74,\ldots)$ of $G(v)$, we define $h':=\lceil\ell\alpha\rceil+j+1$ for the height of $G(1c_v)$, maintaining at all other positions $G(1c_v)$ piecewise parallel to $G(1c_{\alpha})$ (see Figure \ref{pal}). In this case we have\\  
$1c_v=1^20\;1^20\;1^40\;1^20\;1^40\;1^20\;1^50\;10\;1^20\;1^20\;1^40\;1^20\;1^20\;10\;1^50\;10\;1^20\;1^20\;1^20\;10\;1^20\;1^40\;\cdots$,\\
$\Phi_{\mathbb{R}}(1c_v)\approx -2.33103$ and $\lim_{\ell \rightarrow\infty}\big(\frac{h}{\ell}\big)=\frac{\ln(2)}{\ln(3)}$. Trajectory $\mathcal{T}(\Phi_{\mathbb{R}}(1c_v))$ diverges to $-\infty$, by Lemma \ref{ellj}, and there is no accumulation point (see Figure \ref{pei}).
\end{exa}

\begin{res}\label{same traj}
Let $v$ be an aperiodic word having the limit $\Phi_{\mathbb{R}}(v):=\zeta$. Then $C_{\mathbb{R}}(v)=\zeta$ also holds. If $u\rightarrow v$ (a prefix of $v$ such that $\ell(u)\rightarrow\infty$), then the trajectory of $\;C_{\mathbb{R}}(u)$ and the trajectory of $\;\Phi_{\mathbb{R}}(u)$ converge pointwise to the same limit trajectory $\mathcal{T}_v(\zeta)$.\hfill \textup{[Section \ref{infinitecyc}]}$\blacktriangleleft$ 
\end{res}
Since for any non-empty prefix $u$ of $v$ it holds that $T^{\ell}(C_{\mathbb{R}}(u))=T_u(C_{\mathbb{R}}(u))=C_{\mathbb{R}}(u)$ and for all prefixes, but finite many prefixes $\Phi_{\mathbb{R}}(u^{\infty})=C_{\mathbb{R}}(u)$, we arrive to the following (perhaps controversial) interpretation.

\begin{res}
For $u\rightarrow v$, the trajectories of the $C_{\mathbb{R}}(u)$'s converge to the limit trajectory $\mathcal{T}_v(\zeta)$, and instead of the word $u^{\infty}$, arises the purely periodic, transfinite ordinal labeled word $v^{\omega}=v_0v_1\cdots\; v_{\omega}v_{\omega+1}\cdots\; v_{2\omega}v_{2\omega+1}\cdots\; v_{3\omega}v_{3\omega+1}\cdots \; v_{n\omega}v_{n\omega+1}\cdots$ arises, having the aperiodic $v$ as period. 
\end{res}
\noindent With this result we may define for all positive integers $n$
\begin{displaymath}
\Phi_{\mathbb{R}}(v_{n\omega}v_{n\omega+1}\cdots):=\Phi_{\mathbb{R}}(v)=T_{v^n}(\Phi_{\mathbb{R}}(v)),
\end{displaymath}
regarding the transfinite trajectory $\mathcal{T}_{v^{\omega}}(\zeta)$ as the trajectory of $\zeta=\Phi_{\mathbb{R}}(v)$ being infinitely repeated many times.

\begin{res}\label{a}
Let $\alpha=\frac{\ln(2)}{\ln(3)}$. The arithmetic mean $\mu$ and the geometric $\mu_g$ mean of the terms in the series $-\Phi_{\mathbb{R}}(1c_{\alpha})=\sum_{i=1}^{\infty}t_i$ are given by 
\begin{displaymath}
\mu:=\lim_{m\rightarrow\infty}\frac{1}{m}\sum_{i=1}^{m}t_i=\frac{1}{6\ln(2)}\approx 0.240449\; ,\qquad \mu_g:=\lim_{m\rightarrow\infty}\bigg(\prod_{i=1}^{m}t_i\bigg)^{\frac{1}{m}}=\frac{1}{6}\sqrt{2}\approx 0.23570.
\end{displaymath}
\end{res}
\hfill \textup{[Section \ref{promedio}, Lemmas \ref{mu},\ref {mug}]}$\blacktriangleleft$
\begin{res}
Function $\vartheta: \mathbb{R}\rightarrow \{x\in\mathbb{R}\;|\;x\geq 0\}$ defined by $\vartheta(x):=\frac{1}{x\ln(2)}$ if $x\in [\frac{1}{6},\frac{1}{3}]$ and $\vartheta(x)=0$; otherwise, it is a \textit{probability density function}. Indeed, $\int_{-\infty}^{\infty}\vartheta(x)\;dx=1$, and $\mathcal{D}(x)=\int_{-\infty}^{x}\vartheta(t)\;dt=\frac{\ln(x)+\ln(6)}{\ln(2)}$, as a continuous and monotone increasing function, is the \textit{cumulative distribution function} of a random variable $X$ with density $\vartheta(x)$.
\end{res}
\noindent For instance, the probability of finding a term of $-\Phi_{\mathbb{R}}(1c_{\alpha})$ in the interval $[1/6,1/4]$, $\alpha=\ln(2)/\ln(3)$, is given by $\Pr[1/6\leq X\leq 1/4]=\int_{1/6}^{1/4}\frac{1}{t\ln(2)}\;dt\approx 0.58496$. The \textit{expected value} of a random variable $X$ is defined by $E[X]=\int_{-\infty}^{\infty}t\vartheta(t)\;dt=\frac{1}{6\ln(2)}$, thus $E[X]=\mu$. The solution for the equation $\int_{-\infty}^{x}\vartheta(t)\;dt=\frac{1}{2}$ yields the \textit{median} $M[X]$ of the distribution. Hence $M[X]=\frac{1}{6}\sqrt{2}$. Thus we may expect that about ``half'' of the terms lie above and the ``other half'' below the median $M[X]=\mu_g$. 
\section{Basic facts in a new form}\label{proof}

\begin{lem}\label{out}
Let $u=v_0v_1v_2\cdots v_{\ell-1}\in \mathbb{U}$ be a prefix of the parity vector $v=v_0v_1v_2\cdots$ for some $x\in\mathbb{Z}_2$. Then the $\ell$'th iterate of $T$ is given by 
\begin{equation}
T^\ell(x)=\frac{\varphi(u)+3^hx}{2^\ell}\label{ecout}
\end{equation}
where $\ell$ is the length, and $h=\sum_{j=0}^{\ell-1}v_j$ is the height of $u$.
\end{lem}
\noindent{\it Proof}. The formula is a 2-adic version of the relations (2.3), (2.4), (2.5), and (2.6) in \cite{Lag:1985}. Therein $\varrho_k(n)$ corresponds to ${\varphi(u)}$/$2^\ell$ and $\lambda_k(n)$ to $3^h$/$2^\ell$. \hfill $\Box$

The $\ell$'th iterate $T^\ell(x)$ is the $(\ell+1)$'th element of the trajectory $\mathcal{T}(x)$ since we have $v_\ell=T^\ell(x)\bmod 2$. The relation (\ref{ecout}) is a functional equation in $x$ with the unique solution $x=\Phi(v)$ in $\mathbb{Z}_2$. The term ${\varphi(u)}$/$2^\ell$ is a 2-adic fraction, unless the $\ell$ digits of $u$ are 0's.

\begin{lem}\label{Rout}
Let $u=v_0v_1v_2\cdots v_{\ell-1}\in \mathbb{U}$ be the prefix of any infinite word $v=v_0v_1v_2v_3\cdots$ over the alphabet $\{0,1\}$, and $x\in\mathbb{R}$. Then the $(\ell+1)$'th element  of the pseudo-trajectory of $x$ relative to $v$ is given by 
\begin{equation}
T_u(x)=T_u(0)+\frac{3^h}{2^\ell} x\label{Tu}
\end{equation}
where $\ell$ is the length, and $h=\sum_{j=0}^{\ell-1}v_j$ is the height of $u$.
\end{lem}
\noindent{\it Proof}. Function $T_u(x)$ ignores the parity of $x$ but depends on $u$ in the same way as it does $T^\ell(x)$ in (\ref{ecout}). The term ${\varphi(u)}$/$2^\ell$ is a constant. Thus (\ref{Tu}) is the equation of a line with slope $3^h$/$2^\ell$. Therefore, ${\varphi(u)}$/$2^\ell$ = $T_u(0)$. \hfill $\Box$
  
\begin{lem}\label{x+y}
Let $u\in\mathbb{U}$ be a finite word with length $\ell$ and height $h$, and $x,y\in\mathbb{R}$. Then
\begin{equation}\label{x+y1}
T_u(x+y)=T_u(x)+\frac{3^h}{2^\ell}y
\end{equation}
\end{lem}
\noindent{\it Proof}. By Lemma \ref{Rout}, $T_u(x+y)=T_u(0)+\frac{3^h}{2^\ell}(x+y)=T_u(0)+\frac{3^h}{2^\ell}x+\frac{3^h}{2^\ell}y=T_u(x)+\frac{3^h}{2^\ell}y$.\hfill $\Box$

\begin{lem}\label{uw}
Let $u$, $w$ be finite words over the alphabet $\{0,1\}$. If $\ell(w)$ is the length and $h(w)$ the height of $w$, then 
\begin{equation}\label{uw1}
T_{uw}(0)=T_w(0)+\frac{3^{h(w)}}{2^{\ell(w)}}T_u(0).
\end{equation}
\end{lem}
\noindent{\it Proof}. $T_{uw}(0)=\frac{\varphi(uw)}{2^{\ell(u)+\ell(w)}}=\frac{3^{h(w)}\varphi(u)+2^{\ell(u)}\varphi(w)}{2^{\ell(w)}\cdot 2^{\ell(u)}}=\frac{3^{h(w)}\varphi(u)}{2^{\ell(w)}\cdot 2^{\ell(u)}}+\frac{\varphi(w)}{2^{\ell(w)}}$ by (\ref{UV}) and (\ref{tres}).\hfill $\Box$  

\begin{lem}\label{Tv0}
Let $v\neq 0^\infty$ be any infinite word over the alphabet $\{0,1\}$. The elements of the pseudo-trajectory $\mathcal{T}_v(0)$ relative to $v$ are all different from each other.\footnote{even in case of a periodic word $v$.}
\end{lem}
\noindent{\it Proof}. We assume that there is a (non-empty) prefix $u$ and a larger prefix $uw$ of $v$ with lengths $\ell(u)$ and $\ell(u)+\ell(w)$ such that $T_u(0)=T_{uw}(0)$. Then
\begin{displaymath}
\frac{\varphi(u)}{2^{\ell(u)}}=\frac{\varphi(uw)}{2^{\ell(u)+\ell(w)}}\quad \textrm{and}\quad \varphi(u)=\frac{\varphi(uw)}{2^{\ell(w)}}.
\end{displaymath}
Since $\varphi(u)$ is an integer, $uw$ has a prefix of at least $\ell(w)$ zeroes. If $\ell(w)\geq \ell(u)$, then $u=0^{\ell(u)}$ and $\varphi(u)=0$. So $\varphi(uw)=0$, and $uw=0^{\ell(u)+\ell(w)}$ is a prefix of $0$'s.

\noindent If $\ell(w)< \ell(u)$, then $u=0^{\ell(w)}u'$. So $\varphi(u)=2^{\ell(w)}\varphi(u')$ and $\varphi(uw)=2^{\ell(w)}\varphi(u'w)$. But $\varphi(u)=\varphi(uw)/2^{\ell(w)}$ yields 
\begin{displaymath}
\varphi(u')=\frac{\varphi(u'w)}{2^{\ell(w)}}.
\end{displaymath} 
The same arguments given above for $u$ can be repeated now for $u'$ with length $\ell(u')<\ell(u)$. This process comes to an end, because $u$ is a finite word. Hence $T_u(0)=T_{uw}(0)$ implies that all digits of $uw$ are 0's.\hfill $\Box$ 

\begin{lem}\label{in}
Let $u=v_0v_1\cdots v_{\ell-1}\in \mathbb{U}$ be a finite word with height $h=\sum_{j=0}^{\ell-1}v_j$, then we have
\begin{equation}
\Phi(u)=-\frac{\varphi(u)}{3^h}.\label{in1}
\end{equation}
\end{lem}
\noindent{\it Proof}. The word $u$ is the prefix of the eventually periodic word $v_0v_1\cdots v_{\ell-1}0^\infty$. Then apply (\ref{PHI}), and simplify with (\ref{fi}).\hfill $\Box$

\begin{lem}\label{in0}
Let $u\in \mathbb{U}$ be a finite word with length $\ell$ and height $h$. Then $T^\ell(\Phi(u))=0$.
\end{lem}
\noindent{\it Proof}. Substitute $x=\Phi(u)$ in (\ref{ecout}), and simplify with (\ref{in1}).\hfill $\Box$

\begin{lem}\label{id}
Let $u\in \mathbb{U}$ be a finite word with length $\ell$ and height $h$. Then we have the identity
\begin{equation}
T^\ell\bigg(\Phi(u)+\frac{2^\ell}{3^h} x\bigg)=x\quad (\textrm{for all }x\in\mathbb{Z}_2).\label{id1}
\end{equation}
\end{lem}

\noindent{\it Proof}. Let $v=v_0v_1v_2\cdots$ be the parity vector of $(\Phi(u)+\frac{2^\ell}{3^h} x)$ for some $x\in\mathbb{Z}_2$. Then we have
\begin{displaymath}
T(\Phi(u)+\frac{2^\ell}{3^h} x)=\left\{ \begin{array}{ll}
\frac{\Phi(u)}{2}+\frac{2^{\ell-1}}{3^h} x=T(\Phi(u))+\frac{2^{\ell-1}}{3^h} x &\textrm{if }v_0=0;\\[5pt]
\frac{3\Phi(u)+1}{2}+\frac{2^{\ell-1}}{3^{h-1}} x=T(\Phi(u))+\frac{2^{\ell-1}}{3^{h-1}} x &\textrm{if }v_0=1.
\end{array}\right.
\end{displaymath}
For $1\leq i\leq\ell$ define $h_i:=\sum_{j=0}^{i-1}v_j$ and $h_0:=0$. Then
\begin{displaymath}
T^i(\Phi(u)+\frac{2^\ell}{3^h} x)=\left\{ \begin{array}{ll}
\frac{1}{2}\cdot T^{i-1}(\Phi(u)+\frac{2^{\ell}}{3^h} x)=T^i(\Phi(u))+\frac{2^{\ell-i}}{3^{h-h_{i-1}}} x &\textrm{if }v_i=0;\\[5pt]
\frac{1}{2}\bigg(3\cdot T^{i-1}\big(\Phi(u)+\frac{2^{\ell}}{3^h} x\big)+1)\bigg)=T^i(\Phi(u))+\frac{2^{\ell-i}}{3^{h-h_{i-1}-1}} x &\textrm{if }v_i=1.
\end{array}\right.
\end{displaymath}
Let $u=u_0u_1\cdots u_{\ell-1}$ be the finite word with height $h=\sum_{j=0}^{\ell-1}u_j$. The word $u0^\infty$ is the parity vector of $\Phi(u)$. From $v_0=(\Phi(u)+\frac{2^\ell}{3^h} x)\bmod{2}=(\Phi(u))\bmod{2}$, and $\Phi(u)\equiv u\pmod{2}$ follows $v_0=u_0$. Furthermore, for $1\leq i<\ell$ we have
\begin{displaymath}
T^i(\Phi(u)+\frac{2^\ell}{3^h} x)\equiv T^i(\Phi(u))+\frac{2^{\ell-i}}{3^{h-(\ldots)}}x\equiv T^i(\Phi(u))\pmod{2}.
\end{displaymath}
Thus $v_{i-1}=u_{i-1}$, and $u$ is a prefix of $v$. Therefore, $u=v_0v_1\cdots v_{\ell-1}$ and $h=\sum_{j=0}^{\ell-1}v_j$.

\noindent For $i=\ell$ we have two possibilities. If $v_{\ell-1}=0$, then $h_{\ell-1}=\sum_{j=0}^{\ell-2}v_j=h$. Thus $h-h_{\ell-1}=0$. If $v_{\ell-1}=1$, then $h_{\ell-1}=h-1$. Thus $h-h_{\ell-1}-1=h-(h-1)-1=0$. Hence in both cases
\begin{displaymath}
T^\ell(\Phi(u)+\frac{2^\ell}{3^h} x)=T^\ell(\Phi(u))+\frac{2^0}{3^0}x=0+x=x.\\[-15pt]
\end{displaymath}
\hfill $\Box$

\begin{lem}\label{ur}
Let $u=v_0v_1\cdots v_{\ell-1}\in \mathbb{U}$ be a prefix with height $h=\sum_{j=0}^{\ell-1}v_j$ of any infinite word $v=ur$ over the alphabet $\{0,1\}$. Then it holds that
\begin{equation}
\Phi(ur)=\Phi(u)+\frac{2^\ell}{3^h}\Phi(r).\label{ur1}
\end{equation}
\end{lem}

\noindent{\it Proof}. The number of 1's in $u=u0^\infty$ is $h$. Thus $\Phi(u)=\Phi(2^{d_0}+2^{d_1}+\cdots +2^{d_{h-1}})$. We have ${d_{h-1}}\leq\ell-1$ since ${d_{h-1}}=\ell-1 \Leftrightarrow v_{\ell-1}=1$. Also ${d_{h}}\geq\ell$ since ${d_{h}}=\ell \Leftrightarrow v_{\ell}=1$. The exponent $d_h$ points to the first $1$ in $r=v_\ell v_{\ell+1}v_{\ell+2}\cdots$. Since $d_h-\ell\geq 0$, we can write
\begin{displaymath}
\Phi(ur)=-\frac{2^{d_0}}{3}-\cdots -\frac{2^{d_{h-1}}}{3^h}+\frac{2^\ell}{3^h}\bigg(-\frac{2^{d_h-\ell}}{3}-\frac{2^{d_{h+1}-\ell}}{3^2}-\cdots\bigg).\\[-15pt]
\end{displaymath}
\hfill $\Box$

\begin{lem}\label{infinitos}
Let $u=u_0u_1\cdots u_{\ell-1}\in \mathbb{U}$ be fixed. The set of all 2-adic integers having $u$ as prefix of their parity vector is given by
\begin{displaymath}
\bigg\{\Phi(u)+\frac{2^{\ell}}{3^h}x\;\bigg|\;x\in\mathbb{Z}_2\bigg\}.
\end{displaymath}
\end{lem}

\noindent{\it Proof}. Any 2-adic integer with $u$ as prefix has the form $u_0u_1\cdots u_{\ell-1}v_{\ell}v_{\ell+1}\cdots$. Then  $\Phi(u_0u_1\cdots u_{\ell-1}v_{\ell}v_{\ell+1}\cdots)=\Phi(u)+\frac{2^{\ell}}{3^h}\Phi(v_{\ell}v_{\ell+1}\cdots)$. We define $x:=\Phi(v_{\ell}v_{\ell+1}\cdots)$. Digits $v_{\ell}v_{\ell+1}\cdots$ are not fixed, so that they represent any of all possible $2$-adic integers $y=v_{\ell}v_{\ell+1}\cdots$. But $\Phi$ is bijective and, therefore, for every $x\in\mathbb{Z}_2$ an appropriated  $y=v_{\ell}v_{\ell+1}\cdots$ exists such that $x=\Phi(y)$.\hfill $\Box$ 

\begin{lem}\label{ciclo}
Let $u\in \mathbb{U}$ be the prefix with length $\ell$ and height $h$ of the parity vector $v$ of $x$.\\ 
\noindent In $\mathbb{Z}_2$ it holds that
\begin{displaymath}
T^\ell(x)=x \quad\textrm{if and only if}\quad x=\frac{\varphi(u)}{2^\ell-3^h}.
\end{displaymath}
\end{lem}

\noindent{\it Proof}. [only if]. By Lemma \ref{out}, $T^\ell(x)=x$ yields $x=\frac{\varphi(u)}{2^\ell-3^h}$.\\
\noindent [if]. By Lemma \ref{id}, we have
\begin{displaymath}
T^\ell\bigg(\Phi(u)+\frac{2^\ell}{3^h}\cdot \frac{\varphi(u)}{2^\ell-3^h}\bigg)=\frac{\varphi(u)}{2^\ell-3^h}.
\end{displaymath}
But the term inside the bracket is $\frac{\varphi(u)}{2^\ell-3^h}$.\\
\noindent In fact, by Lemma \ref{in} we have $\Phi(u)+\frac{2^\ell}{3^h}\cdot \frac{\varphi(u)}{2^\ell-3^h}=\frac{-\varphi(u)}{3^h}+\frac{2^\ell\varphi(u)}{3^h(2^\ell - 3^h)}=\frac{\varphi(u)}{2^\ell-3^h}$.\hfill $\Box$

\begin{lem}\label{Cu}
Let $u\in\mathbb{U}$ be a finite word with length $\ell$ and height $h$. The rational 2-adic number
\begin{equation}\label{Cu1}
C(u):=\frac{\varphi(u)}{2^\ell-3^h}
\end{equation}
generates a purely cyclic trajectory $\mathcal{T}(C(u))$ with parity vector $u^\infty=uuu\cdots$, and $\mathcal{O}(C(u))$ being a finite cycle.
\end{lem}
\noindent{\it Proof}. Clearly $C(u)\in\mathbb{Q}_{odd}$, and $T^\ell(C(u))=C(u)$ by Lemma \ref{ciclo}.\hfill $\Box$ 

\section{Auxiliary results}

Let $u\in\mathbb{U}$ be a prefix with length $\ell$ and height $h$ converging to the aperiodic word $v$, that is, $\ell\rightarrow\infty$. Throughout this paper we summarize this situation by using the economic notation $u\rightarrow v$.

\begin{lem}\label{2/3}
Let $u\rightarrow v$.
\begin{equation}\label{pend}
\underline{\lim}\;\bigg(\frac{h}{\ell}\bigg)>\frac{\ln(2)}{\ln(3)}\Longrightarrow \bigg(\frac{2^\ell}{3^h}\bigg)\rightarrow 0,\qquad 
\overline{\lim}\;\bigg(\frac{h}{\ell}\bigg)<\frac{\ln(2)}{\ln(3)}\Longrightarrow \bigg(\frac{3^h}{2^\ell}\bigg)\rightarrow 0.
\end{equation}
\end{lem}
\noindent{\it Proof}. [Left side]. There is a $\delta>0$ such that $\underline{\lim}\;\big(\frac{h}{\ell}\big)-\frac{\ln(2)}{\ln(3)}>\delta$. Then there is a $L\in\mathbb{N}$ such that $\frac{h}{\ell}\geq \frac{\ln(2)}{\ln(3)}+\delta$ for all $\ell>L$. So we have $\frac{2^\ell}{3^h}=\frac{e^{\ell\ln(2)}}{e^{h\ln(3)}}\leq e^{\ell\ln(2)-\ell \big(\frac{\ln(2)}{\ln(3)}+\delta\big)\ln(3)}=e^{-\ell\delta\ln(3)}$.
The last term converges to 0 for $\ell\rightarrow\infty$, since $\delta > 0$.

[Right side]. There is a $\delta>0$ such that $\frac{\ln(2)}{\ln(3)}-\overline{\lim}\;\big(\frac{h}{\ell}\big)>\delta$. Then there is a $L\in\mathbb{N}$ such that $\frac{h}{\ell}\leq\frac{\ln(2)}{\ln(3)}-\delta$ for all $\ell>L$. So we have $\frac{3^h}{2^\ell}=\frac{e^{h\ln(3)}}{e^{\ell\ln(2)}}\leq e^{\ell \big(\frac{\ln(2)}{\ln(3)}-\delta\big)\ln(3)-\ell\ln(2)}=e^{-\ell\delta\ln(3)}$.\hfill $\Box$ 

Relations (\ref{pend}) and (\ref{rel}) show that either $\Phi(u)$ and $C(u)$ are both negative, or that $T_u(0)$ and $C(u)$ are both positive for a sufficiently large $u$. Cases $\underline{\lim}=\frac{\ln(2)}{\ln(3)}$ or $\overline{\lim}=\frac{\ln(2)}{\ln(3)}$ need special attention.

\begin{lem}\label{ln2ln3}
Let $\alpha=\frac{\ln(2)}{\ln(3)}$. Then it holds for all $\ell\in\mathbb{N}$
\begin{equation}\label{ln2ln31}
\frac{1}{3}<\frac{2^\ell}{3^{\lceil\ell\alpha\rceil}}<1,\qquad \frac{1}{3}<\frac{3^{\lfloor\ell\alpha\rfloor}}{2^{\ell}}<1,\qquad 1<\frac{2^\ell}{3^{\lfloor\ell\alpha\rfloor}}<3,\qquad 1<\frac{3^{\lceil\ell\alpha\rceil}}{2^{\ell}}<3.
\end{equation}
\end{lem}
\noindent{\it Proof}. Let $\xi:=\lceil\ell\alpha\rceil-\ell\alpha$. We have the identity $\frac{2^\ell}{3^{\ell\frac{\ln(2)}{\ln(3)}}}=1$. Therefore, $\frac{2^\ell}{3^{\lceil\ell\alpha\rceil}}=\frac{2^\ell}{3^{\ell\alpha+\xi}}=\frac{1}{3^\xi}$ and $\frac{3^{\lfloor\ell\alpha\rfloor}}{2^{\ell}}=\frac{3^{\ell\alpha+\xi-1}}{2^\ell}=3^{\xi -1}$. Since $0<\xi<1$, both values lie in the open interval $(\frac{1}{3},1)$.\hfill $\Box$

Let $G(v):=(\ell,h)_{\ell=1}^\infty$ be the sequence of points representing geometrically the aperiodic word $v$ and $\alpha=\ln(2)/\ln(3)$. The coordinates of points in $G(1c_\alpha)$ are $(\ell,\lceil\ell\alpha\rceil)$ and those in $G(0c_\alpha)$ are $(\ell,\lfloor\ell\alpha\rfloor)$. Points of $G(v)$ above line $h=\ell\alpha$ have coordinates $(\ell,\lceil\ell\alpha\rceil+n_\ell)$ and those below line $(\ell,\lfloor\ell\alpha\rfloor-m_\ell)$, where $n_\ell$ and $m_\ell$ are non-negative integers depending on $\ell$ such that $0\leq n_\ell\leq \ell-\lceil\ell\alpha\rceil$ and $0\leq m_\ell\leq\lfloor\ell\alpha\rfloor$.
\begin{lem}\label{=}
Let $u\rightarrow v$ be such that
\begin{displaymath}
\lim_{\ell\rightarrow\infty}\bigg(\frac{h}{\ell}\bigg)=\frac{\ln(2)}{\ln(3)},\quad\textrm{where }h=\lceil\ell\alpha\rceil+n_{\ell}\textrm{ , or}\quad h=\lfloor\ell\alpha\rfloor-m_\ell \;.\\[-10pt]
\end{displaymath}
\begin{displaymath}
\noindent\textrm{Then},\quad\bigg(\frac{2^\ell}{3^h}\bigg)\rightarrow 0 \quad\textrm{if and only if the following three conditions are fullfilled:}
\end{displaymath}
\begin{equation}\label{above}
\frac{h}{\ell}>\frac{\ln(2)}{\ln(3)}\quad \textrm{for all but finitely many points of G(v)}, \quad\frac{n_\ell}{\ell}\rightarrow 0,\quad n_\ell\rightarrow \infty.
\end{equation}
\begin{displaymath}
\noindent\textrm{Furthermore},\quad\bigg(\frac{3^h}{2^\ell}\bigg)\rightarrow 0 \quad\textrm{if and only if the following three conditions are fulfilled:}
\end{displaymath}
\begin{equation}\label{below}
\frac{h}{\ell}<\frac{\ln(2)}{\ln(3)}\quad \textrm{for all but finitely many points of G(v)}, \quad\frac{m_\ell}{\ell}\rightarrow 0,\quad m_\ell\rightarrow \infty.
\end{equation}
\end{lem}
\noindent{\it Proof}. [if]. Let $\alpha=\ln(2)/\ln(3)$. If $h>\ell\alpha$, then $\frac{1}{3^{n_\ell+1}}<\frac{2^\ell}{3^{\lceil\ell\alpha\rceil+n_\ell}}<\frac{1}{3^{n_\ell}}$ by Lemma \ref{ln2ln3}. If $n_\ell\rightarrow \infty$, then $\frac{1}{3^{n_\ell+1}}\rightarrow 0$ and $\frac{1}{3^{n_\ell}}\rightarrow 0$. Thus $\frac{2^\ell}{3^{\lceil\ell\alpha\rceil+n_\ell}}\rightarrow 0$. Finally, $\frac{n_\ell}{\ell}\rightarrow 0$ if and only if $\frac{h}{\ell}=\frac{\lceil\ell\alpha\rceil+n_\ell}{\ell}=\frac{\lceil\ell\alpha\rceil}{\ell}+\frac{n_\ell}{\ell}$ converges to $\alpha$.

If $h<\ell\alpha$, then $\frac{1}{3^{m_{\ell+1}}}<\frac{3^{\lfloor\ell\alpha\rfloor-m_\ell}}{2^\ell}<\frac{1}{3^{m_\ell}}$ by Lemma \ref{ln2ln3}. If $m_\ell\rightarrow \infty$, then $\frac{1}{3^{m_\ell+1}}\rightarrow 0$ and $\frac{1}{3^{m_\ell}}\rightarrow 0$. Thus $\frac{3^{\lfloor\ell\alpha\rfloor-m_\ell}}{2^\ell}\rightarrow 0$. Finally, $\frac{m_\ell}{\ell}\rightarrow 0$ if and only if $\frac{h}{\ell}=\frac{\lfloor\ell\alpha\rfloor-m_\ell}{\ell}=\frac{\lfloor\ell\alpha\rfloor}{\ell}-\frac{m_\ell}{\ell}$ converges to $\alpha$.

[only if]. There are only three possibilities for the points of $G(v)$: $h>\ell\alpha$ for all, except for a finite number of them,  
$h<\ell\alpha$ for all, except for a finite many of them, or they change infinitely, often on the side of $h=\ell\alpha$.\footnote{More precisely, line $y=\alpha x$ for $x\in\mathbb{R}$.} In the latter case, many of them have coordinates $(\ell,\lceil\ell\alpha\rceil)\in G(1c_\alpha)$ infinitely, and many of them have coordinates $(\ell,\lfloor\ell\alpha\rfloor)\in G(0c_\alpha)$ infinitely. By Lemma \ref{ln2ln3} there are infinitely many terms larger than $1/3$ in the sequences $(2^\ell/3^h)_{\ell\rightarrow\infty}$ and $(3^h/2^\ell)_{\ell\rightarrow\infty}$. Finally, if $\{n_\ell\;|\;\ell\in\mathbb{N}\}$ has an upper bound, say $n_\ell<L$ for all $\ell\in\mathbb{N}$, then all terms of $(2^\ell/3^h)_{\ell\rightarrow\infty}$ are larger than $1/3^{L+1}$. A similar argument holds for $m_\ell$.\hfill $\Box$ 

For an aperiodic word $v$ let $\Theta$ be the set of those \textit{sub-sequences} of $\big(\frac{h}{\ell}\big)_{\ell=1}^\infty$ which converge to limit $\frac{\ln(2)}{\ln(3)}$. This set may be empty. 

\begin{lem}\label{geq}
Let $u\rightarrow v$ be such that $\underline{\lim}\;\big(\frac{h}{\ell}\big)=\frac{\ln(2)}{\ln(3)}$.
\begin{displaymath}
\textrm{It holds that } \bigg(\frac{2^\ell}{3^h}\bigg)_{\ell=1}^\infty\rightarrow 0\; \textrm{ if and only if the elements of }\Theta\; \textrm{satisfy the conditions } (\ref{above}).
\end{displaymath}
\end{lem}
\noindent{\it Proof}. Applying Lemma \ref{2/3} (left side) and Lemma \ref{=}.\hfill $\Box$ 

\begin{lem}\label{leq}
Let $u\rightarrow v$ be such that $\overline{\lim}\;\big(\frac{h}{\ell}\big)=\frac{\ln(2)}{\ln(3)}$.
\begin{displaymath}
\textrm{It holds that } \bigg(\frac{3^h}{2^\ell}\bigg)_{\ell=1}^\infty\rightarrow 0\; \textrm{ if and only if the elements of }\Theta\; \textrm{satisfy the conditions } (\ref{below}).
\end{displaymath}
\end{lem}
\noindent{\it Proof}. Applying Lemma \ref{2/3} (right side) and Lemma \ref{=}.\hfill $\Box$

\begin{lem}\label{-2/3}
Let $u\rightarrow v$.
\begin{equation}\label{-pend}
\bigg(\frac{2^\ell}{3^h}\bigg)_{\ell=1}^{\infty}\rightarrow 0\;\;\Longrightarrow \underline{\lim}\;\bigg(\frac{h}{\ell}\bigg)\geq\frac{\ln(2)}{\ln(3)},\qquad 
\bigg(\frac{3^h}{2^\ell}\bigg)_{\ell=1}^{\infty}\rightarrow 0\;\;\Longrightarrow \overline{\lim}\;\bigg(\frac{h}{\ell}\bigg)\leq\frac{\ln(2)}{\ln(3)}.
\end{equation}
\end{lem}
\noindent{\it Proof}. We prove the contrapositive of the statements.

\noindent [Left side]. If $\underline{\lim}\;\big(\frac{h}{\ell}\big)<\frac{\ln(2)}{\ln(3)}$, then there are infinitely many terms $(\frac{h}{\ell})$ smaller than $\frac{\ln(2)}{\ln(3)}$. But $\frac{h}{\ell}<\frac{\ln(2)}{\ln(3)}$ implies $\frac{2^\ell}{3^h}>1$. Thus $\big(\frac{2^\ell}{3^h}\big)_{\ell=1}^{\infty}$ does not converge to 0.

\noindent [Right side]. If $\overline{\lim}\;\big(\frac{h}{\ell}\big)>\frac{\ln(2)}{\ln(3)}$, then there are infinitely many terms $(\frac{h}{\ell})$ larger than $\frac{\ln(2)}{\ln(3)}$. But $\frac{h}{\ell}>\frac{\ln(2)}{\ln(3)}$ implies $\frac{3^h}{2^\ell}>1$.\hfill $\Box$

\section{Proof of Theorem 1}\label{main proof}

Let $u\in\mathbb{U}$ be a prefix with length $\ell$ and height $h$ converging to the aperiodic word $v$, that is, $u\rightarrow v$. The sum $\Phi_{\mathbb{R}}(u)$ is a \textit{partial sum} of the first $h$ terms of $\Phi(v)$ evaluated in $\mathbb{R}$. The sequence $\big(\Phi_{\mathbb{R}}(u)\big)_{u\rightarrow v}$ is monotone decreasing since $\Phi_{\mathbb{R}}(u0)=\Phi_{\mathbb{R}}(u)$ and $\Phi_{\mathbb{R}}(u1)=\Phi_{\mathbb{R}}(u)-\frac{2^\ell}{3^{h+1}}$, and its terms are negative rational 2-adic integers. Therefore, if it can be shown that this sequence has a lower bound then limit $\Phi_{\mathbb{R}}(v)$ exists.

\begin{lem}\label{necessary}
Let $u\rightarrow v$. 
\begin{displaymath}
\textrm{ If } \Phi_{\mathbb{R}}(v)\textrm{ exists, then}\quad\bigg(\frac{2^\ell}{3^h}\bigg)_{\ell=1}^{\infty}\rightarrow 0\quad\textrm{and therefore,}\quad\underline{\lim}\;\bigg(\frac{h}{\ell}\bigg)\geq\frac{\ln(2)}{\ln(3)}.\\[-10pt] 
\end{displaymath}
\end{lem}
\noindent{\it Proof}. Every prefix $u\in\mathbb{U}$ of length $\ell$ terminating in $1$ yields a term $\frac{2^{\ell-1}}{3^h}$. Height $h$ counts the number of terms in $\Phi_{\mathbb{R}}(u)$. Terms corresponding to consecutive $0$'s which precede the last $1$ do not contribute to the sum, but they increase by a factor 2. If limit $\Phi_{\mathbb{R}}(v)$ exists, then the terms of $\Phi_{\mathbb{R}}(u)$ must converge to $0$. Thus $\frac{2^{\ell-1}}{3^h}\rightarrow 0$ for $u\rightarrow v$. Hence $\big(\frac{2^\ell}{3^h}\big)_{\ell=1}^{\infty}\rightarrow 0$, and $\underline{\lim}\big(\frac{h}{\ell}\big)\geq\frac{\ln(2)}{\ln(3)}$ by Lemma \ref{-2/3}.\hfill $\Box$

\begin{lem}\label{costura}
Let $v=ur$ be an aperiodic word with prefix $u\in\mathbb{U}$.\\  
If $\Phi_{\mathbb{R}}(v)=\zeta$ exists, then $\Phi_{\mathbb{R}}(r)$ also exists, namely
\begin{equation}\label{costura1}
\Phi_{\mathbb{R}}(r)=T_u(\zeta).
\end{equation}
If $\Phi_{\mathbb{R}}(r)$ exists, then $\Phi_{\mathbb{R}}(v)$ also exists.
\end{lem}

\noindent{\it Proof}. Let $u=v_0v_1\cdots v_{\ell-1}$. By Scheme \ref{esquema2} we have\\
\noindent  $\Phi_{\mathbb{R}}(S^\ell(v))=T_{v_{\ell-1}}\circ \cdots\circ T_{v_1}\circ T_{v_0}(\zeta)=T_u(\zeta)$, and $S^\ell(v)=v_\ell v_{\ell+1}\cdots =r$.\\
\noindent If $\Phi_{\mathbb{R}}(r)$ exists, then $\Phi_{\mathbb{R}}(u)+\frac{2^\ell}{3^h}\Phi_{\mathbb{R}}(r)=\Phi_{\mathbb{R}}(ur)$ by Lemma \ref{ur}, since $\Phi_{\mathbb{R}}(u)\in\mathbb{Q}_{odd}$.\hfill $\Box$

\begin{lem}\label{01}
Let $v=u01r$ be an aperiodic word and $u\in\mathbb{U}^*$. If $\Phi_{\mathbb{R}}(u01r)$ exists, then
\begin{equation}\label{011}
\Phi_{\mathbb{R}}(u10r)=\Phi_{\mathbb{R}}(u01r)+\frac{2^{\ell (u)}}{3^{h(u)+1}}.
\end{equation}
\end{lem}
\noindent{\it Proof}. We have \\$\Phi_{\mathbb{R}}(u10r)=\Phi_{\mathbb{R}}(u)+\frac{2^{\ell(u)}}{3^{h(u)}}\big(-\frac{1}{3}+\frac{2^2}{3}\Phi_{\mathbb{R}}(r)\big)$; $\Phi_{\mathbb{R}}(u01r)=\Phi_{\mathbb{R}}(u)+\frac{2^{\ell(u)}}{3^{h(u)}}\big(-\frac{2}{3}+\frac{2^2}{3}\Phi_{\mathbb{R}}(r)\big)$.\hfill $\Box$

\begin{lem}\label{lifting}
Let $1c_\alpha$ be a Sturmian word with an irrational slope $0<\alpha<1$, and $v=v_0v_1\cdots$ an infinite word (aperiodic or periodic) with the geometric representation $G(v)=(\ell,h)_{\ell=1}^\infty$. If $\frac{h}{\ell}>\alpha$ holds for all points of $G(v)$, then $1c_\alpha$ can be mapped (lifted) onto $v$ replacing successively $(01)$'s by $(10)$'s in an appropriate order.  
\end{lem}
\noindent{\it Proof}.

\begin{figure}[htbp]
\begin{center}
\epsfxsize=3.5in
\epsfbox{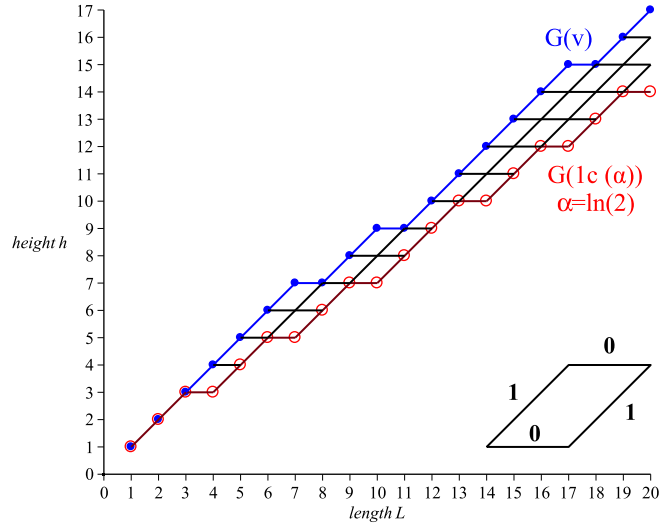}
\end{center}
\caption{Proof of Lemma \ref{lifting}}
\label{lift}
\end{figure}

By joining the points of $G(v)$ by a polygonal line $L_v$ and those of $G(1c_\alpha)$ by $L_{1c_\alpha}$, the area between $L_v$ and $L_{1c_\alpha}$ can be subdivided into small rhomboids, each of them with two horizontal sides representing 0's, and the other two sides representing 1's in a $45$ degree angle with the axis $\ell$. Each rhomboid represents a replacement $01\rightarrow 10$. If there are $2$ or more rhomboids at the same level, the replacements are done from right to left. For instance, $0001\rightarrow 0010\rightarrow 0100\rightarrow 1000$.

The number of rhomboids from one level to the next one increases or decreases at most $1$. The first replacement $01\rightarrow 10$ corresponds to the rightmost rhomboid at the bottom.\hfill $\Box$

\begin{lem}\label{limarriba}
Let $v$ be an aperiodic word.
\begin{displaymath}
\textrm{If}\quad\underline{\lim}\;\bigg(\frac{h}{\ell}\bigg)_{\ell=1}^\infty>\frac{\ln(2)}{\ln(3)}\;,\quad\textrm{then limit }\Phi_{\mathbb{R}}(v)<-1 \textrm{ exists}.
\end{displaymath}
\end{lem}
\noindent{\it Proof}. There is an irrational $\alpha$ such that $\underline{\lim}\;\big(\frac{h}{\ell}\big)_{\ell=1}^\infty>\alpha>\frac{\ln(2)}{\ln(3)}$. There is a $L\in\mathbb{N}$ such that $\frac{h}{\ell}>\alpha$ for all $\ell\geq L$. Let $L$ be the smallest one with this property. We compare $v$ with $1c_\alpha$. The limit $\Phi_{\mathbb{R}}(1c_\alpha)$ exists, since $\alpha>\frac{\ln(2)}{\ln(3)}$ (\cite{Lopez:Lop}, Section 4). Let $1c_\alpha=s_0s_1\cdots s_{L-1}s_Ls_{L+1}\cdots$ and $v=v_0v_1\cdots v_{L-1}v_Lv_{L+1}\cdots$. Since $L$ is minimal, we have $v_{L-1}=1$, $s_{L-1}=0$, and the prefix $s_0\cdots s_{L-2}$ and $u:=v_0\cdots v_{L-2}$ reach the same height $h=\lceil (L-1)\alpha\rceil=\sum_0^{L-2}s_i=\sum_0^{L-2}v_i$.

\noindent We apply Lemma \ref{lifting}. For all points of $G_L(v)=\big(\ell,h\big)_{\ell=L}^\infty$, $\frac{h}{\ell}>\alpha$ holds, thus $\frac{h}{\ell}\geq \frac{\lceil \ell\alpha\rceil}{\ell}$, so that these points do not fall below $G_L(1c_\alpha)$. Limit $\Phi_{\mathbb{R}}(us_{L-1}s_L\cdots)$ exists, because $\Phi_{\mathbb{R}}(s_{L-1}s_L\cdots)$ can be calculated by $\Phi_{\mathbb{R}}(s_{L-1}s_L\cdots)=T_{s_0s_1\cdots s_{L-2}}\big(\Phi_{\mathbb{R}}(1c_\alpha)\big)$ and therefore, $\Phi_{\mathbb{R}}(us_{L-1}s_L\cdots )=\Phi_{\mathbb{R}}(u)+\frac{2^{L-1}}{3^h}\;\Phi_{\mathbb{R}}(s_{L-1}s_L\cdots)$. Hence $us_{L-1}s_L\cdots$ can be lifted onto $uv_{L-1}v_L\cdots$, replacing (01)'s by (10)'s, first in $s_{L-1}s_L$. The lifting yields $-\infty<\Phi_{\mathbb{R}}(us_{L-1}s_L\cdots )<\Phi_{\mathbb{R}}(u10s_{L+2}\cdots )\leq\Phi_{\mathbb{R}}(uv_{L-1}v_L\cdots)=\Phi_{\mathbb{R}}(v)$ by Lemma \ref{01}. Finally, $v$ can be lifted onto $111\cdots$, thus $\Phi_{\mathbb{R}}(v)\leq -1$. But $v$ is aperiodic, so that $\Phi_{\mathbb{R}}(v)< -1$.\hfill $\Box$ 

\begin{lem}\label{infinito}
Let $v$ be an aperiodic word having the limit $\zeta=\Phi_{\mathbb{R}}(v)$. 
The elements of trajectory $\mathcal{T}(\zeta)$ are all different from each other, so that orbit $\mathcal{O}(\zeta)$ is an infinite set.
\end{lem}
\noindent{\it Proof}. We assume that two elements of $\mathcal{T}(\zeta)$ are 
equal, then $T_u(\zeta)=T_{u'}(\zeta)$ for different (non-empty) prefixes $u$ 
and $u'$ of $v$. But $T_u(0)+\frac{3^h}{2^\ell}\zeta=T_{u'}(0)+\frac{3^{h'}}{2
^{\ell '}}\zeta$ is a linear equation in $\zeta$ with rational coefficients. 
Thus $\zeta$ is a rational number. 

\noindent If $\zeta\notin\mathbb{Q}_{odd}$, then $\zeta$ is a fraction in its lowest terms with an odd numerator and an even denominator. If $2^n$ divides the denominator of $\zeta$, then $2^{n+1}$ divides the denominator of $T_{v_0}(\zeta)$, $2^{n+2}$ divides the denominator of $T_{v_0v_1}(\zeta)$, \ldots, introducing in each step an additional divisor $2$. However, it follows from $T_u(\zeta)=\frac{p}{q}$, $T_{u'}(\zeta)=\frac{p'}{q'}$, and $\frac{p}{q}=\frac{p'}{2^m q}$ ($m> 0$), a contradiction given that numerator $p'=2^m p$ is even.

\noindent If $\zeta\in\mathbb{Q}_{odd}$, then the parity of $\zeta$ is well-defined. Thus the function $T_{v_i}$ ($i\in\mathbb{N}_0$) does the same as the map $T$, and $v$ is the aperiodic parity vector of $\zeta$. However, we have $T^i(\zeta)\neq T^j(\zeta)$ for all $i\neq j$ by Scheme {\ref{esquema1}} in contradiction with the assumption that two elements of $\mathcal{T}(\zeta)$ are equal.\hfill $\Box$

 \begin{lem}\label{1cbound}
Let $\alpha$ be irrational and $1>\alpha>\frac{\ln(2)}{\ln(3)}$. The orbit of $\Phi_{\mathbb{R}}(1c_{\alpha})$ is bounded. $\Phi_{\mathbb{R}}(0c_{\alpha})=3\Phi_{\mathbb{R}}(1c_{\alpha})+1$ is a lower bound and $\Phi_{\mathbb{R}}(1c_{\alpha})$ is an upper bound.
\end{lem}
\noindent{\it Proof}. Let $\Phi_{\mathbb{R}}(s_Ls_{L+1}\cdots)$ be any element of the trajectory of $\Phi_{\mathbb{R}}(1c_{\alpha})=\Phi_{\mathbb{R}}(s_0s_1\cdots s_L\cdots)$. Let $G(1c_{\alpha})=\big(\ell,h\big)_{\ell=1}^{\infty}$ be the geometric representation of $1c_{\alpha}$. We have $h=\sum_{i=0}^{\ell-1}s_i$. Since $1c_{\alpha}=\big(\lceil(j+1)\alpha\rceil-\lceil j\alpha\rceil\big)_{j=0}^{\infty}=s_0s_1s_2\cdots s_j\cdots$, we have also $h=\lceil(j+1)\alpha\rceil$. Indeed, the telescopic sum $\lceil(j+1)\alpha\rceil=\sum_{i=0}^{j}(\lceil(i+1)\alpha\rceil-\lceil i\alpha\rceil )$ yields $\sum_{i=0}^{j}s_i=h$.

\noindent We represent $s_Ls_{L+1}\cdots$ geometrically as a word starting from the origin, that is, $s_Ls_{L+1}\cdots=s'_0s'_1s'_2\cdots$. Digits $s_Ls_{L+1}s_{L+2}\cdots$ are given by $\big\lceil (j+1)\alpha -\lceil L\alpha\rceil\big\rceil - \big\lceil j\alpha - \lceil L\alpha\rceil\big\rceil$ for $j=L,L+1,\ldots$. Substitution $k:=j-L$ yields digits $s'_0s'_1s'_2\cdots$, namely $\big\lceil(k+L+1)\alpha-\lceil L\alpha\rceil\big\rceil-\big\lceil(k+L)\alpha-\lceil L\alpha\rceil\big\rceil$ for $k=0, 1, 2\ldots$. Let $G(s'_0s'_1s'_2\cdots)=\big(\ell,h'\big)_{\ell=1}^{\infty}$ be the geometric representation. Then $h'=\big\lceil(k+L+1)\alpha-\lceil L\alpha\rceil\big\rceil$.

\noindent We have $h-h'=\lceil(j+1)\alpha\rceil-\big\lceil(k+L+1)\alpha-\lceil L\alpha\rceil\big\rceil$. Let $\delta:= L\alpha-\lceil L\alpha\rceil$. Thus $-1<\delta < 0$. Then $h-h'=\lceil (j+1)\alpha\rceil-\lceil(k+1)\alpha+\delta\rceil$. Since $j$ and $k$ take simultaneously  the same value, we have $h-h'=\lceil (k+1)\alpha\rceil-\lceil(k+1)\alpha+\delta\rceil$. 

\noindent The well known inequality $\lceil(k+1)\alpha\rceil+\lceil\delta\rceil-1 \leq \lceil(k+1)\alpha+\delta\rceil \leq\lceil(k+1)\alpha\rceil+\lceil\delta\rceil$ yields $-\lceil\delta\rceil+1 \geq \lceil (k+1)\alpha\rceil-\lceil(k+1)\alpha+\delta\rceil \geq \lceil\delta\rceil$. Hence $1\geq h-h'\geq 0$. It follows that $h\geq h'$ and $h'\geq h-1$. Height $h$ is for points of $G(1c_{\alpha})$ and $h-1$ for those of $G(0c_{\alpha})$. Hence $0c_{\alpha}$ can be lifted onto $s'_0s'_1s'_2\cdots$ and $s'_0s'_1s'_2\cdots$ onto $1c_{\alpha}$. Therefore, $\Phi_{\mathbb{R}}(0c_{\alpha})\leq \Phi_{\mathbb{R}}(s'_0s'_1s'_2\cdots)\leq \Phi_{\mathbb{R}}(1c_{\alpha})$, where $\Phi_{\mathbb{R}}(s'_0s'_1s'_2\cdots)=\Phi_{\mathbb{R}}(s_Ls_{L+1}\cdots)$.
By Lemma \ref{uplo1lim} $\Phi_{\mathbb{R}}(0c_{\alpha})=3\Phi_{\mathbb{R}}(1c_{\alpha})+1$.\hfill $\Box$ 

\begin{lem}\label{orbitarriba}
Let $v$ be an aperiodic word.
\begin{displaymath}
\textrm{If}\quad\underline{\lim}\;\bigg(\frac{h}{\ell}\bigg)_{\ell=1}^\infty>\frac{\ln(2)}{\ln(3)}\;,\;\textrm{ then the orbit of }\;\Phi_{\mathbb{R}}(v) \textrm{ has accumulation points}.
\end{displaymath}
\end{lem}
\noindent{\it Proof}. Limit $\Phi_{\mathbb{R}}(v)$ exists by Lemma \ref{limarriba}. Let $k_0<k_1<\cdots <k_i<\cdots$ be indexes in $v=v_0v_1\cdots v_{k_i}\cdots$ and $G_{k_i}(v)=\big(\ell,h_i\big)_{\ell=1}^\infty$ the sequence of points representing geometrically the suffix $v_{k_i}v_{k_i+1}v_{k_i+2}\cdots$, where $h_i:=\sum_{j=0}^{\ell-1}v_{k_i+j}$. $\underline{\lim}\;\big(\frac{h_i}{\ell}\big)_{\ell=1}^\infty =\underline{\lim}\;\big(\frac{h}{\ell}\big)_{\ell=1}^\infty   >\frac{\ln(2)}{\ln(3)}$ holds for the elements of $G_{k_i}(v)$ (\cite{Monks:Yaz}, Lemma 4.1). Limits $\Phi_{\mathbb{R}}(v_{k_i}v_{k_i+1}\cdots)$ are elements of the trajectory $\mathcal{T}(\Phi_{\mathbb{R}}(v))$. We apply the arguments given for $\Phi_{\mathbb{R}}(v)$ in the proof of Lemma \ref{limarriba}  on each of these $\Phi_{\mathbb{R}}(v_{k_i}v_{k_i+1}\cdots)$. 

\noindent Let $\alpha$ be an irrational such that $\underline{\lim}\;\big(\frac{h_i}{\ell}\big)_{\ell=1}^\infty>\alpha>\frac{\ln(2)}{\ln(3)}$ and $L_i$ is the smallest positive integer such that $\frac{h_i}{\ell}>\alpha$ for all $\ell\geq L_i$ in $G_{k_i}(v)$. Furthermore we define $u^{(i)}:=v_{k_i}v_{k_i+1}\cdots v_{k_i+L_i-1}$ and $1c_{\alpha}=s_0s_1\cdots s_{L_i-1}\cdots$. Then $\Phi_{\mathbb{R}}(u^{(i)}s_{L_i}s_{L_i+1}\cdots)\leq \Phi_{\mathbb{R}}(u^{(i)}v_{k_i+L_i}v_{k_i+L_i+1}\cdots)$ and therefore, $\Phi_{\mathbb{R}}(s_{L_i}s_{L_i+1}\cdots)\leq \Phi_{\mathbb{R}}(v_{k_i+L_i}v_{k_i+L_i+1}\cdots)$. Since $\Phi_{\mathbb{R}}(0c_{\alpha})\leq \Phi_{\mathbb{R}}(s_{L_i}s_{L_i+1}\cdots)$ by Lemma \ref{1cbound}, we conclude that for all $L_i$ bound $\;\Phi_{\mathbb{R}}(0c_{\alpha})\leq\Phi_{\mathbb{R}}(v_{k_i+L_i}v_{k_i+L_i+1}\cdots)<-1$ holds. Limits $\Phi_{\mathbb{R}}(v_{k_i+L_i}v_{k_i+L_i+1}\cdots)$ are elements of trajectory $\mathcal{T}(\Phi_{\mathbb{R}}(v))$. If we choose indexes $k_i$ so that $k_0+L_0<k_1$, $k_1+L_1<k_2$, $\ldots, k_i+L_i<k_{i+1},\ldots$, then $\Phi_{\mathbb{R}}(v_{k_i+L_i}v_{k_i+L_i+1}\cdots)$ are all different from each other by Lemma \ref{infinito}. Hence $\{\Phi_{\mathbb{R}}(v_{k_i+L_i}v_{k_i+L_i+1}\cdots)\;\big|\;i=0,1,\ldots\}$ is an infinite set with an upper and a lower bound, having at least one accumulation point.\hfill $\Box$

\noindent {\bf Proof of Theorem \ref{aperiodic}.} The word $v$ is aperiodic.

If $\underline{\lim}\;\big(\frac{h}{\ell}\big)_{\ell=1}^\infty>\frac{\ln(2)}{\ln(3)}$, then limit $\Phi_{\mathbb{R}}(v)$ exists by Lemma \ref{limarriba}, and the orbit of $\Phi_{\mathbb{R}}(v)$ has accumulation points by Lemma \ref{orbitarriba}, which implies $\Phi_{\mathbb{R}}(v)\notin\mathbb{Q}_{odd}$. From (\ref{IL}) and (\ref{IL2}) it follows that the $2$-adic $\Phi(v)$ is aperiodic.\\

If $\underline{\lim}\;\big(\frac{h}{\ell}\big)_{\ell=1}^\infty<\frac{\ln(2)}{\ln(3)}$, then $\Phi(v)\notin\mathbb{Q}_{odd}$, because a divergent trajectory of an $\Phi(v)\in\mathbb{Q}_{odd}$ implies $\underline{\lim}\;\big(\frac{h}{\ell}\big)_{\ell=1}^\infty\geq\frac{\ln(2)}{\ln(3)}$ (Monks, Yazinski \cite{Monks:Yaz}, Theorem 2.7 b). Thus $\Phi(v)$ is aperiodic. 

Hence $\underline{\lim}\;\big(\frac{h}{\ell}\big)_{\ell=1}^\infty=\frac{\ln(2)}{\ln(3)}$ is the only remaining possibility.\hfill $\Box$

\section{The word \textbf{$1c_v$} associated to \textbf{$v$}}\label{unocv}

The Sturmian word $1c_{\alpha}$ for $\alpha=\frac{\ln(2)}{\ln(3)}$ has exactly 3 different factors of length 2, namely 11, 10 and 01. There is no factor 00. When we replace many 0's by 1's infinitely in $1c_{\alpha}=1101101101011011010\cdots$ leaving many 0's unchanged infinitely, we obtain an infinite word $v'$ which coincides piecewise with $1c_{\alpha}$. The following Lemma shows that in the special case of $\underline{\lim}\;\big(\frac{h}{\ell}\big)_{\ell=1}^\infty=\frac{\ln(2)}{\ln(3)}$, limit $\Phi_{\mathbb{R}}(v')$ may not exist. 

\begin{lem}\label{piecewise}
Let $\alpha=\frac{\ln(2)}{\ln(3)}$. Let $1c_{\alpha}$ be written as a product of factors $u^{(0)}=s_{0}\cdots s_{\ell_{0}-2}0$, $u^{(1)}=s_{\ell_0}\cdots s_{\ell_{1}-2}0$, $u^{(2)}=s_{\ell_1}\cdots s_{\ell_{2}-2}0$, ... terminating all of them with 0, and having length $\ell_{j+1}-\ell_j$ and height $h_j$ for $j=0,1,2,\cdots$. Let $v'$ be the infinite word obtained by the substitution $\big(s_{\ell{j}-1}=0\big)_{j=0}^{\infty}\longrightarrow\big(s_{\ell{j}-1}=1\big)_{j=0}^{\infty}$. Then it holds that
\begin{displaymath}
\overline{\lim}\;\bigg(\frac{h_{j+1}}{h_j}\bigg)<\frac{3}{2}\Longrightarrow \Phi_{\mathbb{R}}(v') \textrm{ exists},\qquad 
\underline{\lim}\;\bigg(\frac{h_{j+1}}{h_{j}}\bigg)>6\Longrightarrow \Phi_{\mathbb{R}}(v')\rightarrow -\infty.\\[-10pt]
\end{displaymath}
\end{lem}
\noindent{\it Proof}. The terms of $-\Phi_{\mathbb{R}}(1c_{\alpha})$ are given by $\bigg(\frac{2^{\lfloor\frac{i}{\alpha}\rfloor}}{3^{1+i}}\bigg)_{i=0}^{\infty}:=\big(t_i\big)_{i=0}^{\infty}$ (\cite{Lopez:Lop}, Lemma 11). Function $[0,\infty)\rightarrow\mathbb{R}$ defined by $\tau(x):=\frac{2^{\lfloor\frac{x}{\alpha}\rfloor}}{3^{1+x}}$ has a range $(\frac{1}{6},\frac{1}{3}]$. The $t_i$'s defined for non-negative integers fall within this range and therefore, $\frac{1}{6}<t_i<\frac{1}{3}$ for $i\neq 0$. The divergent series $-\Phi_{\mathbb{R}}(1c_{\alpha})$ can be written formally as
$-\Phi_{\mathbb{R}}(1c_{\alpha})=\sum_{i=0}^{h_0-1}t_i+\sum_{i=h_0}^{h_1-1}t_i+\sum_{i=h_1}^{h_2-1}t_i+\cdots$. Introducing the arithmetic means $\mu_j:=\frac{1}{h_j}\sum_{i={h_{j-1}}}^{h_{j}-1}t_i$ for $j\neq 0$ and $\mu_0:=\frac{1}{h_0}\sum_{i=0}^{h_0-1}t_i$, we have formally $-\Phi_{\mathbb{R}}(1c_{\alpha})=h_0\mu_0+h_1\mu_1+h_2\mu_2+\cdots$.

\noindent Digit $s_{\ell{j}-1}=0$ is replaced by $s_{\ell{j}-1}=1$. Since $s_{\ell_j}=1$, there is an $i$ such that $d_i=\ell_j$ by equation (\ref{PHI}). Thus $t_i=\frac{2^{d_i}}{3^{1+i}}=\frac{2^{\ell_j}}{3^{1+i}}$. For the same $i$ we have now $d_i=\ell_j-1$ pointing to the new 1, and the corresponding term is $\frac{2^{\ell_j-1}}{3^{1+i}}=\frac{t_i}{2}:=t_j'$. For $s_{\ell_j}=1$ and all following 1's the $i$'s have been increased by 1, dividing the corresponding terms $t_i$ by 3. Hence $-\Phi_{\mathbb{R}}(v')=h_0\mu_0+t_0'+\frac{1}{3}h_1\mu_1+\frac{1}{3}t_1'+\frac{1}{9}h_2\mu_2+\frac{1}{9}t_2'+\cdots$ where $\frac{1}{12}<t_j'<\frac{1}{6}$ are the terms generated by the additional 1's. Since $\frac{3}{2}\frac{1}{12}<t_0'+\frac{1}{3}t_1'+\frac{1}{9}t_2'+\cdots <\frac{3}{2}\frac{1}{6}$, we have $-\Phi_{\mathbb{R}}(v')=\vartheta +\sum_{j=0}^{\infty}\frac{1}{3^j}h_j\mu_j$ where $\frac{1}{8}<\vartheta <\frac{1}{4}$.

\noindent We apply the ratio test and use the (rough) estimate $\frac{1}{6}h_j<\sum_{i=h_{j-1}}^{h_j-1}t_i<\frac{1}{3}h_j\;$  ($j\neq 0$). For the ratio we have $\frac{1}{3}\frac{h_{j+1}\mu_{j+1}}{h_j\mu_j}=\frac{1}{3}\frac{\sum_{i=h_{j}}^{h_{j+1}-1}t_i}{\sum_{i=h_{j-1}}^{h_j-1}t_i}$. Thus $\frac{1}{3}\cdot\frac{1}{2}\cdot\frac{h_{j+1}}{h_j} < \frac{1}{3}\cdot\frac{\sum_{i=h_{j}}^{h_{j+1}-1}t_i}{\sum_{i=h_{j-1}}^{h_j-1}t_i} < \frac{1}{3}\cdot 2\cdot\frac{h_{j+1}}{h_j}$. Hence, if
$\overline{\lim}\big(\frac{h_{j+1}}{h_j}\big)<\frac{3}{2}$ the series $-\Phi_{\mathbb{R}}(v')$ converges, and if $\underline{\lim}\big(\frac{h_{j+1}}{h_j}\big)>6$ the series diverges.\hfill $\Box$

\begin{lem}\label{max}
Let $v$ be an aperiodic word such that 
$\underline{\lim}\; \big(\frac{h}{\ell}\big)_{\ell=1}^{\infty}=\frac{\ln(2)}{\ln(3)}=\alpha$, and let $n_{\ell}=h(\ell)-\lceil \ell\alpha\rceil$ be the relative height of $v$ over the Sturmian word $1c_{\alpha}$ at position $\ell$ such that $\big( \frac{2^{\ell}}{3^{\lceil\ell\alpha\rceil+n_{\ell}}}\big)_{\ell=1}^{\infty}\rightarrow 0$.\\
For any non-negative integer $j$ $\;\max\{\ell\;|\;n_{\ell}=j\}$ exists.
\end{lem}
\noindent{\it Proof}. By applying Lemma \ref{geq}, the conditions (\ref{above}) are fulfilled. There is a positive integer $L$ such that $\frac{h}{\ell}>\frac{\lceil\ell\alpha\rceil}{\ell}$ for all $\ell\geq L$. Let $L$ be the smallest one with this property. Let $1c_\alpha=s_0s_1\cdots s_{L-1}s_Ls_{L+1}\cdots$ and $v=v_0v_1\cdots v_{L-1}v_Lv_{L+1}\cdots$. Since $L$ is minimal, we have $v_{L-1}=1$, $s_{L-1}=0$, and the prefix $s:=s_0\cdots s_{L-2}$ and $u:=v_0\cdots v_{L-2}$ reach the same height $\lceil (L-1)\alpha\rceil=\sum_0^{L-2}s_i=\sum_0^{L-2}v_i$. Thus $n_{L-1}=0$. From $n_{\ell}\rightarrow\infty$ it follows that $\max\{\ell\;|\;n_{\ell}=j\}$ exists for any non-negative integer $j$.\hfill $\Box$

Let $v=v_0v_1\cdots v_i\cdots$ be an aperiodic word such that 
$\underline{\lim}\; \big(\frac{h}{\ell}\big)_{\ell=1}^{\infty}=\frac{\ln(2)}{\ln(3)}$, $\big( \frac{2^{\ell}}{3^{\lceil\ell\alpha\rceil+n_{\ell}}}\big)_{\ell=1}^{\infty}\rightarrow 0$ and $1c_{\alpha}=s_0s_1\cdots s_i\cdots$ for $\alpha=\frac{\ln(2)}{\ln(3)}$. The index $i$ in $v_i$ is related with the index $\ell$ in $n_{\ell}$ by $i=\ell-1$. Since $i=\max\{\ell\;|\;n_{\ell}=j\}-1$ is the largest index where $n_{\ell}=j$ for some fixed $j$, we have $v_{i+1}=1$ and $s_{i+1}=0$. We replace these special 0's in $1c_{\alpha}$ by 1's applying Lemma \ref{max} and define the resulting word as follows:

\begin{defin}\label{def1cv}
Let $v=v_0v_1\cdots v_i\cdots$ be an aperiodic word such that $\underline{\lim}\; \big(\frac{h}{\ell}\big)_{\ell=1}^{\infty}=\frac{\ln(2)}{\ln(3)}=\alpha$ and $\big( \frac{2^{\ell}}{3^{h}}\big)_{\ell=1}^{\infty}\rightarrow 0$. Let $\ell_j:=\max\{\ell\;|\;n_{\ell}=j\}$ and $1c_{\alpha}=s_0s_1\cdots s_i\cdots$.\\
Word $1c_v:=s_0's_1'\cdots s_i'\cdots$ ---called the associated word to $v$--- is defined by $s_i':=1$ if $i\in\{\ell_j\;|\;j=0,1,2,\ldots\}$ otherwise, $s_i':=s_i$.
\end{defin}
Note that $s_0'=1$. Word $1c_v$ coincides piecewise with $1c_{\alpha}$, wherein infinitely many $0$'s have been replaced by $1$'s, depending on $v$, by means of Lemma \ref{max}.

\begin{lem}\label{iff}
Let $v$ be an aperiodic word such that 
$\underline{\lim}\; \big(\frac{h}{\ell}\big)_{\ell=1}^{\infty}=\frac{\ln(2)}{\ln(3)}=\alpha$, and let $n_{\ell}=h(\ell)-\lceil \ell\alpha\rceil$ be the relative height of $v$ over the Sturmian word $1c_{\alpha}$ at position $\ell$ such that $\big( \frac{2^{\ell}}{3^{\lceil\ell\alpha\rceil+n_{\ell}}}\big)_{\ell=1}^{\infty}\rightarrow 0$.\\
If $\Phi_{\mathbb{R}}(1c_v)$ exists, then limit $\Phi_{\mathbb{R}}(v)$ also exists.
\end{lem}
\noindent{\it Proof}. We write $1c_{\alpha}=us^{(0)}s^{(1)}\cdots s^{(j)}\cdots$ where $u=s_0s_1\cdots s_{\ell_{0}-1}$ is a prefix,\footnote{$\ell_0-1$ is the largest index $i$ such that $n_{\ell_0}=0$ (see Lemma \ref{max}).} and
$s^{(0)}=0s_{\ell_{0}+1}\cdots s_{\ell_{1}-1}$, $s^{(1)}=0s_{\ell_{1}+1}\cdots s_{\ell_{2}-1}$, $s^{(j)}=0s_{\ell_{j}+1}\cdots s_{\ell_{j+1}-1}$ are factors with length $\ell_{j+1}-\ell_j$, having $s_{\ell_j}=0$ as its first digit. These 0's will be replaced by 1's to get $1c_v$ according to $\ell_j:=\max\{\ell\;|\;n_{\ell}=j\}$. Consider the words \\
$v^{(j)}:=v_{\ell_j}v_{\ell_{j}+1}\cdots v_{\ell_{j+1}-1}=1v_{\ell_{j}+1}\cdots v_{\ell_{j+1}-1}$;$\;\quad s'^{(j)}=s_{\ell_j}'s_{\ell_j+1}'\cdots s_{\ell_{j+1}-1}'=1s_{\ell_j+1}\cdots s_{\ell_{j+1}-1}$.

\noindent Note that $v_{\ell_j}$ is the ($\ell_j+1$)'th digit of $v$. Thus $n_{\ell_j}=j$ and $n_{\ell_{j+1}}=j+1$. We have $n_{\ell}\geq j+1$ for $\ell_j+1\leq\ell\leq\ell_{j+1}$, because $n_{\ell_{j+1}}=j+1$ is the last $n_{\ell}$ with value $j+1$. Hence for the points of $G(v^{(j)}):=\big(\ell,h\big)_{\ell=1+\ell_j}^{\ell_{j+1}}$, $n_{\ell}\geq j+1$ and $n_{\ell_j+1}= n_{\ell_{j+1}}= j+1$ hold. 

\noindent Word $s'^{(j)}$ is a factor of $1c_v$. The 0's in $1c_{\alpha}$ with index $i\in\{\ell_0,\ell_1,\ldots,\ell_j\}$ have been replaced by 1's. With each of these $j+1$ substitutions the corresponding pieces $G_k(1c_{\alpha}):=\big(\ell,h''\big)_{\ell=1+\ell_k}^{\infty}$ for $k=0,1,\ldots ,j$ have been moved one unit upwards. Hence for the points of $G(s'^{(j)}):=\big(\ell,h'\big)_{\ell=1+\ell_j}^{\ell_{j+1}}$, $n'_{\ell}\geq 0$ holds. Since $n_{\ell}\geq j+1$, we have $n'_{\ell}=n_{\ell}-(j+1)\geq 0$. Thus $h-h'\geq 0$. Also, $G(v^{(j)})$ and $G(s'^{(j)})$ coincide in their first and their last points, so that they have the same length and the same height. It follows $\varphi\big(v^{(j)}\big)\leq \varphi\big(s'^{(j)}\big)$ (\cite{Halb:1997}, Lemma 4). Then $\Phi_{\mathbb{R}}(v^{(j)})\geq \Phi_{\mathbb{R}}(s'^{(j)})$ by Lemma \ref{in}.

\noindent We have $v=uv^{(0)}v^{(1)}\cdots v^{(j)}\cdots$ and $1c_v=us'^{(0)}s'^{(1)}\cdots s'^{(j)}\cdots$.
Points $G(v)=\big(\ell,h\big)_{\ell=1}^{\infty}$ do not fall below points $G(1c_v)=\big(\ell,h'\big)_{\ell=1}^{\infty}$. Sequences $G(v)$ and $G(1c_v)$ coincide at least in those points with $\ell\in \{\ell_0,\ell_0+1,\;\ell_1,\ell_1+1,\;\ell_2,\ell_2+1,\cdots\}$. Thus $1c_v$ can be lifted onto $v$. If $\Phi_{\mathbb{R}}(1c_v)$ exists, then $\Phi_{\mathbb{R}}(1c_v)\leq \Phi_{\mathbb{R}}(v)<-1$.\hfill $\Box$

\begin{lem}\label{Phi1cv}
Let $v$ be an aperiodic word such that 
$\underline{\lim}\; \big(\frac{h}{\ell}\big)_{\ell=1}^{\infty}=\frac{\ln(2)}{\ln(3)}=\alpha$, and let $n_{\ell}=h(\ell)-\lceil \ell\alpha\rceil$ be the relative height of $v$ over the Sturmian word $1c_{\alpha}$ at position $\ell$ such that $\big( \frac{2^{\ell}}{3^{\lceil\ell\alpha\rceil+n_{\ell}}}\big)_{\ell=1}^{\infty}\rightarrow 0$.\\
If $\Phi_{\mathbb{R}}(1c_v)$ exists, then
\begin{equation}\label{serie1cv}
\Phi_{\mathbb{R}}(1c_v)=\Phi_{\mathbb{R}}(u)+\frac{2^{\ell_0}}{3^{\lceil \ell_0\alpha\rceil}}\Phi_{\mathbb{R}}(s'^{(0)})+
\frac{2^{\ell_1}}{3^{\lceil\ell_1\alpha\rceil+1}}\Phi_{\mathbb{R}}(s'^{(1)})\cdots +\frac{2^{\ell_j}}{3^{\lceil\ell_j\alpha\rceil+j}}\Phi_{\mathbb{R}}(s'^{(j)})+\cdots,
\end{equation}
where $1c_v=us'^{(0)}s'^{(1)}\cdots s'^{(j)}\cdots$, $u=s_0s_1\cdots s_{L-1}\cdots s_{\ell_{0}-1}$, and $s'^{(j)}=1s_{\ell_j+1}\cdots s_{\ell_{j+1}-1}$.\\$L-1$ is the smallest index and $\ell_0-1$ the largest index $i$ in $1c_{\alpha}$ for $\alpha=\frac{\ln(2)}{\ln(3)}$ such that $n_L=n_{\ell_0}=0$.
If $\Phi_{\mathbb{R}}(v)$ exists, then
\begin{equation}\label{seriev}
\Phi_{\mathbb{R}}(v)=\Phi_{\mathbb{R}}(u)+\frac{2^{\ell_0}}{3^{\lceil \ell_0\alpha\rceil}}\Phi_{\mathbb{R}}(v^{(0)})+
\frac{2^{\ell_1}}{3^{\lceil\ell_1\alpha\rceil+1}}\Phi_{\mathbb{R}}(v^{(1)})\cdots +\frac{2^{\ell_j}}{3^{\lceil\ell_j\alpha\rceil+j}}\Phi_{\mathbb{R}}(v^{(j)})+\cdots,
\end{equation}
where $v^{(j)}=1v_{\ell_{j}+1}\cdots v_{\ell_{j+1}-1}$.
\end{lem}
\noindent{\it Proof}. Apply Lemma \ref{ur}.\hfill $\Box$

\begin{lem}\label{ellj}
Let $v$ be an aperiodic word such that 
$\underline{\lim}\; \big(\frac{h}{\ell}\big)_{\ell=1}^{\infty}=\frac{\ln(2)}{\ln(3)}=\alpha$, and let $n_{\ell}=h(\ell)-\lceil \ell\alpha\rceil$ be the relative height of $v$ over the Sturmian word $1c_{\alpha}$ at position $\ell$ such that $\big( \frac{2^{\ell}}{3^{\lceil\ell\alpha\rceil+n_{\ell}}}\big)_{\ell=1}^{\infty}\rightarrow 0$.\\
Lengths $\ell_{j+1}-\ell_j$ of words $s'^{(j)}$ in the associated word $1c_v$ diverge to $\infty$ for $j\rightarrow\infty$.
\end{lem}
\noindent{\it Proof}. It should be noted that words $s'^{(j)}= 1s_{\ell_j}\cdots s_{\ell_{j+1}-2}$ have a $\ell_{j+1}-\ell_j$ length. The sequence of slopes $\big(\frac{\lceil\ell_j\alpha\rceil+j}{\ell_j}\big)_{j=0}^{\infty}$  converges to $\alpha$. It follows that $\frac{\lceil\ell_{j+1}\alpha\rceil+j+1}{\ell_{j+1}}-\frac{\lceil\ell_j\alpha\rceil+j}{\ell_j}$ converges to $0$.

\noindent We have the following sum of slopes (adding numerator with numerator and denominator with denominator):
\begin{displaymath}
\frac{\lceil\ell_j\alpha\rceil+j}{\ell_j}\oplus\frac{\lceil\ell_{j+1}\alpha\rceil+j+1-\lceil\ell_j\alpha\rceil-j}{\ell_{j+1}-\ell_j} = \frac{\lceil\ell_{j+1}\alpha\rceil+j+1}{\ell_{j+1}}
\end{displaymath}
The second term is the slope of the line segment defined by points $(\ell_j,\lceil\ell_j\alpha\rceil+j)$ and $(\ell_{j+1},\lceil\ell_{j+1}\alpha\rceil+j+1)$ which together with $(0,0)$ form a triangle. Since the slopes of the other two sides of the triangle converge both to $\alpha$ and their difference converges to $0$, the slope of the line segment converges to $\alpha$ as well. Hence, the irrational $\alpha$ is the limit of the proper fractions $\frac{\lceil\ell_{j+1}\alpha\rceil+j+1-\lceil\ell_j\alpha\rceil-j}{\ell_{j+1}-\ell_j}$ and therefore, $\ell_{j+1}-\ell_j$ diverges to $\infty$ for $j\rightarrow\infty$.\hfill $\Box$

\section{Infinite cycles}\label{infinitecyc}

\begin{lem}\label{Fi-C}
Let $u\rightarrow v$.
\begin{displaymath}
\textrm{If}\quad \lim_{u\rightarrow v}\Phi_{\mathbb{R}}(u)=\Phi_{\mathbb{R}}(v),\quad\textrm{then}\quad \lim_{u\rightarrow v}C_{\mathbb{R}}(u)=\Phi_{\mathbb{R}}(v).\\[-15pt]
\end{displaymath}
\end{lem} 
\noindent{\it Proof}. In $\mathbb{R}$ we have $|\Phi(u)-C(u)|=\Big|\frac{-\varphi(u)2^\ell}{3^h(2^\ell-3^h)}\Big|=|\Phi(u)|\cdot\Big|\frac{2^\ell}{2^\ell-3^h}\Big|<|\Phi_{\mathbb{R}}(v)|\cdot\bigg|\frac{\frac{2^\ell}{3^h}}{\frac{2^\ell}{3^h}-1}\bigg|$. But $\frac{2^\ell}{3^h}\rightarrow 0$, thus $|\Phi(u)-C(u)|\rightarrow 0$.\hfill $\Box$ 

If $\Phi_{\mathbb{R}}(v)$ exists, then $|\Phi_{\mathbb{R}}(u)-C_{\mathbb{R}}(u)|\rightarrow 0$. If $\Phi_{\mathbb{R}}(v)$ does not exist ($=-\infty$), then $|\Phi_{\mathbb{R}}(u)-C_{\mathbb{R}}(u)|\rightarrow 0$ can be true (For instance, choose $f(n)=2^n$ in Example \ref{bigO}). Hence $|\Phi_{\mathbb{R}}(u)-C_{\mathbb{R}}(u)|\rightarrow 0$ does not guarantee that $\Phi_{\mathbb{R}}(v)$ exists. In contrast with this fact in $\mathbb{R}$, in $\mathbb{Z}_2$ (respectively in $\mathbb{Z}_3$) Lemma \ref{2y3} holds.

\begin{lem}\label{2y3}
Let $u\rightarrow v$. Then
\begin{equation}\label{23}
\Big|\Phi(u)-C(u)\Big|_2\rightarrow 0,\qquad \Big|T_u(0)-C(u)\Big|_3\rightarrow 0 .
\end{equation}
\end{lem}
\noindent{\it Proof}. $\Big| \Phi(u)-C(u) \Big|_2=\Big| \frac{-\varphi(u)\cdot 2^\ell}{3^h(2^\ell-3^h)} \Big|_2\leq \frac{1}{2^\ell}$, and $\Big| T_u(0)-C(u) \Big|_3=\Big| \frac{-\varphi(u)\cdot 3^h}{2^\ell(2^\ell-3^h)} \Big|_3= \frac{1}{3^h}$. But $\ell\rightarrow\infty$ implies $h\rightarrow\infty$, since $v$ is aperiodic.\hfill $\Box$

\begin{lem}\label{CT}
Let $v$ be an aperiodic word such that $\overline{\lim}\;\big(\frac{h}{\ell}\big)\leq\frac{\ln(2)}{\ln(3)}$, and the elements of set $\Theta$ satisfy the conditions (\ref{below}) of Lemma \ref{=}.\\
\noindent If the limit $C_{\mathbb{R}}(v,\mathcal{L}):=\lim\;(C_{\mathbb{R}}(u^{(\ell)}))_{\ell\in\mathcal{L}}$ exists, where $u^{(\ell)}$ are prefixes of $v$ whose lengths $\ell$ take progressively the values of some strictly monotone increasing sequence $\mathcal{L}$ of positive integers, then $\lim\;(C_{\mathbb{R}}(u^{(\ell)}))_{\ell\in\mathcal{L}}=\lim\;\big(T_{u^{(\ell)}}(0)\big)_{\ell\in\mathcal{L}}$.
\end{lem}
\noindent{\it Proof}. In $\mathbb{R}$ we have  $|T_u(0)-C(u)|=|\frac{-\varphi(u)(3^{h})}{2^{\ell}(2^{\ell}-3^h)}|=|C(u)\cdot \frac{3^h}{2^{\ell}}|$. For the selected $u^{(\ell)}$ the limit of $C(u)$ exists, and $\frac{3^h}{2^{\ell}}\rightarrow 0$ by Lemma \ref{leq}.\hfill $\Box$ 

\begin{lem}\label{uinfty}
Let $v$ be an aperiodic word. If limit $\Phi_{\mathbb{R}}(v)$ exists, then $\Phi_{\mathbb{R}}(u^{\infty})=C_{\mathbb{R}}(u)$ holds for all, except for a finite number of prefix $u$.
\end{lem}
\noindent{\it Proof}. Let $\ell$ be the length and $h$ the height of $u$.   By Lemma \ref{necessary} there is a positive integer $L$ such that $\frac{2^{\ell}}{3^h}<1$ for all $\ell >L$. Then
$\Phi_{\mathbb{R}}(u^{\infty})=\Phi_{\mathbb{R}}(u)+\frac{2^{\ell}}{3^h}\Phi_{\mathbb{R}}(u)+\frac{2^{2\ell}}{3^{2h}}\Phi_{\mathbb{R}}(u)+\frac{2^{3\ell}}{3^{3h}}\Phi_{\mathbb{R}}(u)\ldots=\Phi_{\mathbb{R}}(u)\frac{1}{1-\frac{2^{\ell}}{3^h}}=-\frac{\varphi(u)}{3^h}\cdot \frac{3^h}{3^h-2^{\ell}}=C_{\mathbb{R}}(u)$ for all $u$ larger than $L$.\hfill $\Box$

\begin{lem}\label{trajC}
Let $v$ be an aperiodic word having the limit $\zeta=\Phi_{\mathbb{R}}(v)$. The trajectory of $\;C_{\mathbb{R}}(u)$ and the trajectory of $\;\Phi_{\mathbb{R}}(u)$ converge pointwise to the same limit trajectory $\mathcal{T}_v(\zeta)$ for $u\rightarrow v$.
\end{lem}
\noindent{\it Proof}. Let $\ell_0$ be any fixed positive integer. The ($\ell_0+1$)'th term of the trajectory of $\zeta=\Phi_{\mathbb{R}}(v)=C_{\mathbb{R}}(v)$ is given by $\Phi_{\mathbb{R}}(v_{\ell_0}v_{\ell_0+1}\cdots)=T_{\texttt{u}_0}(\zeta)$, where $\texttt{u}_0:=v_0v_1\cdots v_{\ell_0-1}$ is the corresponding prefix of $v$ with height $h_0$. For any given $\varepsilon>0$ we define $\delta:=\frac{2^{\ell_0}}{3^{h_0}}\varepsilon$. There is a positive integer $L$ such that for all prefix $u$ with length $\ell>L$ we have $|\zeta-C_{\mathbb{R}}(u)|<\delta$ by Lemma \ref{Fi-C}. Then $\zeta=C_{\mathbb{R}}(u)+\xi$ for some $\xi\in\mathbb{R}$ with $|\xi|<\delta$. By Lemma \ref{x+y} we have $T_{\texttt{u}_0}(\zeta)=T_{\texttt{u}_0}(C_{\mathbb{R}}(u)+\xi)=T_{\texttt{u}_0}(C_{\mathbb{R}}(u))+\frac{3^{h_0}}{2^{\ell_0}}\xi$. Thus $|T_{\texttt{u}_0}(\zeta)-T_{\texttt{u}_0}(C_{\mathbb{R}}(u))|=\frac{3^{h_0}}{2^{\ell_0}}|\xi|<\frac{3^{h_0}}{2^{\ell_0}}\delta=\varepsilon$ for all $u$ larger than $L$.

\noindent A similar argument holds for $\Phi_{\mathbb{R}}(u)$. With the same $\varepsilon$ and $\delta$, there is a positive integer $L'$ such that for every prefix $u'$ with length $\ell>L'$ we have $|\zeta-\Phi_{\mathbb{R}}(u')|<\delta$ and, by Lemma \ref{x+y} $|T_{\texttt{u}_0}(\zeta)-T_{\texttt{u}_0}(\Phi_{\mathbb{R}}(u'))|<\varepsilon$ for all $u'$ larger than $L'$. We thus conclude hence that for all prefixes $u$ larger than $\max(L,L')$, points $T_{\texttt{u}_0}(C_{\mathbb{R}}(u))$ and $T_{\texttt{u}_0}(\Phi_{\mathbb{R}}(u))$ lie inside an $\varepsilon$-neighborhood of $T_{\texttt{u}_0}(\zeta)$.\hfill $\Box$

\section{Upper and lower limits of Sturmian trajectories}

\begin{lem}\label{uplo1lim}
Let $1>\alpha>\frac{\ln(2)}{\ln(3)}$ be irrational. Then $\Phi_{\mathbb{R}}(1c_{\alpha})$ is the upper limit, and $\Phi_{\mathbb{R}}(0c_{\alpha})=3\Phi_{\mathbb{R}}(1c_{\alpha})+1$ is the lower limit of trajectory $\mathcal{T}(\Phi_{\mathbb{R}}(1c_{\alpha}))$.
\end{lem}
\noindent{\it Proof}. Let $\alpha=[0;a_1,a_2,\ldots]$ be the simple continued fraction expansion. The \textit{characteristic} word $c_{\alpha}$ is the limit of the \textit{standard sequence} $(s_n)$, which is recursively defined by $s_{-1}:=1$, $s_0:=0$, and $s_n=s_{n-1}^{d_n}s_{n-2}$ for $n>0$, where $d_1:=a_1-1$ and $d_k:=a_k$ for $k\neq 1$ (\cite{Lot:2002}). The finite words $s_n$ have suffix $10$ for $n>0$ even and suffix $01$ for $n>0$ odd.

\noindent [upper limit]. Consider the sequence $(1s_2,1s_4,1s_6,\ldots)=(1s_{2i}):=\big(u^{(i)}\big)_{i=1}^{\infty}$. We have 
\begin{displaymath}
1s_{2i+2}=1s_{2i+1}^{d_{2i+2}}s_{2i}=1\big(s_{2i}^{d_{2i+1}}s_{2i-1} \big)^{d_{2i+2}}s_{2i}=1s_{2i}\big(\cdots \big)1s_{2i}\;,
\end{displaymath}
because the last bit of $s_{2i-1}$ is 1. Thus $1s_{2i}=u^{(i)}$ is both a suffix and prefix of $1s_{2i+2}=u^{(i+1)}$. The $u^{(i)}$'s taken as prefix of $u^{(i+1)}$ are prefix of $1c_{\alpha}$. Hence $\Phi_{\mathbb{R}}(u^{(i)})\rightarrow\Phi_{\mathbb{R}}(1c_{\alpha})$ for $i\rightarrow\infty$. The $u^{(i)}$'s taken as suffix of $u^{(i+1)}$ are factors of $1c_{\alpha}=v_0v_1v_2\cdots$, starting at indexes $\ell_{1}'<\ell_{2}'<\cdots<\ell_{i}'<\cdots$. Let $\ell_i$ be the length and $h_i$ the height  of $u^{(i)}$. Then we have $\Phi_{\mathbb{R}}(v_{\ell_{i}'}v_{\ell_{i}'+1}v_{\ell_{i}'+2}\cdots) = \Phi_{\mathbb{R}}(u^{(i)})+\frac{2^{\ell_{i}}}{3^{h_{i}}}\Phi_{\mathbb{R}}(v_{\ell_{i}'+\ell_i}v_{\ell_{i}'+\ell_i+1}v_{\ell_{i}'+\ell_i+2}\cdots)$. From $\Phi_{\mathbb{R}}(0c_{\alpha})\leq\Phi_{\mathbb{R}}(v_{\ell_{i}'+\ell_i}v_{\ell_{i}'+\ell_i+1}v_{\ell_{i}'+\ell_i+2}\cdots)\leq \Phi_{\mathbb{R}}(1c_{\alpha})$ and $\frac{2^{\ell_{i}}}{3^{h_{i}}}\rightarrow 0$ for $i\rightarrow\infty$ it follows that $\Phi_{\mathbb{R}}(v_{\ell_{i}'}v_{\ell_{i}'+1}v_{\ell_{i}'+2}\cdots)-\Phi_{\mathbb{R}}(u^{(i)})\rightarrow 0$ for $i\rightarrow\infty$. But $\Phi_{\mathbb{R}}(u^{(i)})\rightarrow\Phi_{\mathbb{R}}(1c_{\alpha})$. Therefore, elements $\Phi_{\mathbb{R}}(v_{\ell_{i}'}v_{\ell_{i}'+1}v_{\ell_{i}'+2}\cdots)$ of the trajectory $\mathcal{T}(\Phi_{\mathbb{R}}(1c_{\alpha}))$ converge to $\Phi_{\mathbb{R}}(1c_{\alpha})$. 

\noindent [lower limit]. Consider the sequence $(0s_3,0s_5,0s_7,\ldots)=(0s_{2i+1}):=\big(u^{(i)}\big)_{i=1}^{\infty}$. Then we have 
\begin{displaymath}
0s_{2i+1}=0s_{2i}^{d_{2i+1}}s_{2i-1}=0\big(s_{2i-1}^{d_{2i}}s_{2i-2} \big)^{d_{2i+1}}s_{2i-1}=0s_{2i-1}\big(\cdots \big)0s_{2i-1}\;,
\end{displaymath}
because the last bit of $s_{2i-2}$ is 0. Thus $0s_{2i-1}=u^{(i-1)}$ is both a suffix and prefix of $0s_{2i+1}=u^{(i)}$. The same arguments as above for $\Phi_{\mathbb{R}}(1c_{\alpha})$ prove the claim for $\Phi_{\mathbb{R}}(0c_{\alpha})$.\\
Since $1c_{\alpha}=2^{d_0}+2^{d_1}+2^{d_2}+\cdots$, 
$0c_{\alpha}=2^{d_0'}+ 2^{d_1'}+ 2^{d_2'}+\cdots$ and $d_0'=d_1, d_1'=d_2, d_2'=d_3,  \cdots$, it holds that $\Phi_{\mathbb{R}}(0c_{\alpha})=3\Phi_{\mathbb{R}}(1c_{\alpha})+1$.\hfill $\Box$

\begin{lem}\label{+1/2}
Let $0<\alpha<\frac{\ln(2)}{\ln(3)}$ be irrational. Then $\Phi^*_{\mathbb{R}}(0c_{\alpha})=3\Phi^*_{\mathbb{R}}(1c_{\alpha})+\frac{1}{2}$.
\end{lem}
\noindent{\it Proof}. Let $0<\alpha<1$ be irrational with convergents $\frac{p_k}{q_k}$ ($k=0,1,2,\ldots$) of the simple continued fraction expansion of $\alpha$. The Christoffel words are defined by
\begin{displaymath}
\overline{m}_{p_k/q_k}:=\bigg(\bigg\lceil(j+1)\frac{p_k}{q_k}\bigg\rceil-\bigg\lceil j\frac{p_k}{q_k}\bigg\rceil\bigg)_{j=0}^{q_k-1}=1z_k0.
\end{displaymath}
The height of $\overline{m}_{p_k/q_k}$ is $p_k$ and its length $q_k$. The \textit{central word} $z_k$ is a palindrome (\cite{Lot:2002}).\\
Since $1z_k0$ and $0z_k1$ have the same length $\ell$ and height $h$, it holds that $\varphi(0z_k1)=3\varphi(1z_k0)+2^{\ell-1}-3^h$. We divide both sides by $2^{\ell}-3^h$ and let $\ell\rightarrow\infty$ (\cite{Lopez:Lop}, Section 4).\hfill $\Box$

\begin{lem}\label{Tu0}
Let $v$ be an aperiodic word. If $\;\overline{\lim}\big(\frac{h}{\ell}\big)<\frac{\ln(2)}{\ln(3)}$, then for any $x\in\mathbb{R}$ the pseudo trajectory $\mathcal{T}_v(x)$ converges pointwise to the pseudo trajectory $\mathcal{T}_v(0)$.
\end{lem}
\noindent{\it Proof}. Let $u\in\mathbb{U}$ be a prefix of $v$ with length $\ell$ and height $h$. We have $T_u(x)=T_u(0)+\frac{3^h}{2^\ell} x$ for any $x\in\mathbb{R}$ by Lemma \ref{Rout}, and $\frac{3^h}{2^{\ell}}\rightarrow 0$ for $\ell\rightarrow\infty$ by Lemma \ref{2/3}. Thus $|T_u(x)-T_u(0)|\rightarrow 0$ for $u\rightarrow v$.\hfill $\Box$ 

\begin{lem}\label{limits}
Let $0<\alpha<\frac{\ln(2)}{\ln(3)}$ be irrational. Then $\Phi^*_{\mathbb{R}}(1c_{\alpha})$ is the lower limit,  and $\Phi^*_{\mathbb{R}}(0c_{\alpha})=3\Phi^*_{\mathbb{R}}(1c_{\alpha})+\frac{1}{2}$ is the upper limit of pseudo trajectory $\mathcal{T}_{1c_{\alpha}}(x)$ for any given $x\in\mathbb{R}$ over the Sturmian word $1c_{\alpha}$.
\end{lem}
\noindent{\it Proof (Sketch)}. Let $\mathcal{S}_{\ell,h}$ be the set of all words with fixed length $\ell$ and fixed height $h$, and let $\mathcal{R}$ be the set of left-rotate permutations of $s=s_0s_1\cdots s_{\ell-1}$ generated by $\rho:=(s_0s_1\cdots s_{\ell-1})\rightarrow(s_1s_2\cdots s_{\ell-1}s_0)$. Halbeisen and Hungerb\H{u}hler  have shown that for any prefix $u\in\mathbb{U}$ of $1c_{\alpha}$ with length $\ell$ and height $h$, $\varphi(u)=\max_{s\in\mathcal{S}_{\ell,h}}\{\min_{t\in\mathcal{R}}\{\varphi(t)\} \}$ holds (\cite{Halb:1997}, Lemma 5).\footnote{Therein $\varphi(u)$ is called $M_{\ell,n}$.} In a similar way it can be shown that for any prefix $u\in\mathbb{U}$ of $0c_{\alpha}$ with length $\ell$ and height $h$, $\varphi(u)=\min_{s\in\mathcal{S}_{\ell,h}}\{\max_{t\in\mathcal{R}}\{\varphi(t)\} \}$ holds.

\noindent [lower limit]. The Christoffel words $\overline{m}_{p_k/q_k}=1z_k0$ are prefixes of $1c_{\alpha}$ for $k$ odd. The orbit of
\begin{displaymath}
\frac{\varphi(\overline{m}_{p_k/q_k})}{2^{q_k}-3^{p_k}}=\Phi^*_{\mathbb{R}}(m_{p_k/q_k})\quad\textrm{(k odd)}\footnote{The word $m_{p_k/q_k}$ is purely periodic with the Christoffel word $\overline{m}_{p_k/q_k}$ as its period.} 
\end{displaymath}
is the finite cycle
\begin{displaymath}
\bigg\{\Phi^*_{\mathbb{R}}(m_{p_k/q_k}), T\bigg(\Phi^*_{\mathbb{R}}(m_{p_k/q_k})\bigg), T^2\bigg(\Phi^*_{\mathbb{R}}(m_{p_k/q_k})\bigg), \ldots , T^{q_k-1}\bigg(\Phi^*_{\mathbb{R}}(m_{p_k/q_k})\bigg)\bigg\},
\end{displaymath}
wherein $\Phi^*_{\mathbb{R}}(m_{p_k/q_k})$ is the smallest element by the result of Halbeisen and Hungerb\H{u}hler. We have $m_{p_k/q_k}\rightarrow 1c_{\alpha}$ for $k$ odd (\cite{Lot:2002}, Section 1.2.2., Definition on page 7), and $\Phi^*_{\mathbb{R}}(m_{p_k/q_k})\rightarrow\Phi^*_{\mathbb{R}}(1c_{\alpha})$ by the continuity of function $F^*(x)=\Phi^*_{\mathbb{R}}(m_x)$ at irrational $\alpha$ (\cite{Lopez:Lop}, Section 4, Definition 28). We hence conclude that $\Phi^*_{\mathbb{R}}(1c_{\alpha})$ is the lower limit of $\mathcal{T}_{1c_{\alpha}}(\Phi^*_{\mathbb{R}}(1c_{\alpha}))$ and also of any other $\mathcal{T}_{1c_{\alpha}}(x)$ by Lemma \ref{Tu0}.

\noindent [upper limit]. The words $\overline{m}\;'_{p_k/q_k}=0z_k1$ are prefixes of $0c_{\alpha}$ for $k$ even. The orbit of $\Phi^*_{\mathbb{R}}(m'_{p_k/q_k})$ is a finite cycle wherein $\Phi^*_{\mathbb{R}}(m'_{p_k/q_k})$ is the largest element. We have $m'_{p_k/q_k}\rightarrow 0c_{\alpha}$ for $k$ even, and it can be shown that $\Phi^*_{\mathbb{R}}(m'_{p_k/q_k})\rightarrow\Phi^*_{\mathbb{R}}(0c_{\alpha})$. Therefore, $\Phi^*_{\mathbb{R}}(0c_{\alpha})$ is the upper limit of $\mathcal{T}_{0c_{\alpha}}(\Phi^*_{\mathbb{R}}(0c_{\alpha}))$ and also of any other $\mathcal{T}_{0c_{\alpha}}(x)$ by Lemma \ref{Tu0}. For any given $x\in\mathbb{R}$ we have
\begin{displaymath}
\mathcal{T}_{1c_{\alpha}}(x)=\mathcal{T}_{c_{\alpha}}(T_1(x))=\mathcal{T}_{c_{\alpha}}\bigg(\frac{3x+1}{2}\bigg)\quad\textrm{and}\quad\mathcal{T}_{0c_{\alpha}}(x)=\mathcal{T}_{c_{\alpha}}(T_0(x))=\mathcal{T}_{c_{\alpha}}\bigg(\frac{x}{2}\bigg).
\end{displaymath}
Both trajectories converge pointwise to $\mathcal{T}_{c_{\alpha}}(0)$ by Lemma \ref{Tu0}. Thus they have the same upper limit. Finally, $\Phi^*_{\mathbb{R}}(0c_{\alpha})=3\Phi^*_{\mathbb{R}}(1c_{\alpha})+\frac{1}{2}$ by Lemma \ref{+1/2}.\hfill $\Box$

\section{Terms of \textbf{$-\Phi_{\mathbb{R}}$}\textbf{$(1c_{\alpha})$} for \textbf{$\alpha=\ln(2)/\ln(3)$}}\label{promedio}

\begin{lem}\label{mu}
Let $\alpha=\frac{\ln(2)}{\ln(3)}$. The arithmetic mean $\mu$ of the terms in $-\Phi_{\mathbb{R}}(1c_{\alpha})=\sum_{i=1}^{\infty}\frac{2^{d_{i-1}}}{3^i}:=\sum_{i=1}^{\infty}t_i$ is given by 
\begin{displaymath}
\mu=\lim_{m\rightarrow\infty}\frac{1}{m}\sum_{i=1}^{m}t_i=\frac{1}{\alpha}\int_1^{1+\alpha}\frac{1}{3^x}\;dx=\frac{1}{6\ln(2)}\approx 0.240449
\end{displaymath}
\end{lem}
\noindent{\it Proof}. The terms of $-\Phi_{\mathbb{R}}(1c_{\alpha})$ are given by $t_i=\frac{2^{\lfloor\frac{i-1}{\alpha}\rfloor}}{3^{i}}$ (\cite{Lopez:Lop}, Lemma 11). The map $x\mapsto 2^{\lfloor\frac{x-1}{\alpha}\rfloor}$ for $x\geq 1$ yields a staircase with steps $2^{n-1}$ over the semi-open intervals $[1+(n-1)\alpha, 1+n\alpha)$ where $n$ is a positive integer. The function $[1,\infty)\rightarrow\mathbb{R}$ defined by $\tau(x):=\frac{2^{\lfloor\frac{x-1}{\alpha}\rfloor}}{3^{x}}$ has the following properties:
\begin{itemize}
	\item The graph of $\tau$ is a set of disjoint, congruent and parallel exponential branches: $\{\frac{2^{n-1}}{3^x}\;|\;x\in [1+(n-1)\alpha, 1+n\alpha),\; n=1,2,3,\ldots\}$.
\item $\tau$ has a period $\alpha$. In fact, $\tau(x+(n-1)\alpha)=\frac{2^{\lfloor(x+(n-1)\alpha-1)\frac{1}{\alpha}\rfloor}}{3^{x+(n-1)\alpha}}=\frac{2^{n-1}}{3^{(n-1)\alpha}}\tau(x)=\tau(x)$, since $3^{\alpha}=2$.
\item The range of $\tau$ is $(\frac{1}{6}, \frac{1}{3}]$, since $\tau(1+(n-1)\alpha)=\tau(1)=\frac{1}{3}$ and $\lim_{x\rightarrow(1+n\alpha)^-}\tau(x)=\frac{1}{6}$.
\item $\tau$ is discontinuous exactly at $\{x=1+n\alpha\;|\;n=1,2,3,\ldots\}$.
\end{itemize} 
Each interval $[1+(n-1)\alpha, 1+n\alpha)$ contains at most one of the integers $i=1,2,3,\ldots$. If $1+(n-1)\alpha\leq i< 1+n\alpha$ for some $i$, then $n-1=\lfloor(i-1)\frac{1}{\alpha}\rfloor=d_{i-1}$. Thus for $x=i$ we have $\frac{2^{n-1}}{3^x}=\frac{2^{d_{i-1}}}{3^i}=t_i$, and $(i,t_i)$ is a point of the branch over the interval $[1+\lfloor\frac{i-1}{\alpha}\rfloor\alpha, \;1+\lfloor\frac{i-1}{\alpha}\rfloor\alpha+\alpha\big)$.\\
We project now the set $\{(i,t_i)\;|\;i\in \mathbb{N}\}$ horizontally onto the first branch $\{\frac{1}{3^x}\;|\;x\in [1, 1+\alpha)\}$. The projection $(i,t_i)\longmapsto\big(i-\lfloor(i-1)\frac{1}{\alpha}\rfloor\alpha,\;t_i\big)$ is injective, because $i-\lfloor(i-1)\frac{1}{\alpha}\rfloor\alpha=i'-\lfloor(i'-1)\frac{1}{\alpha}\rfloor\alpha$ implies that the integer $i-i'$ is a multiple of the irrational $\alpha$.\\
We have $i-\lfloor(i-1)\frac{1}{\alpha}\rfloor\alpha=1+\alpha\big((i-1)\frac{1}{\alpha}-\lfloor(i-1)\frac{1}{\alpha}\rfloor\big)$. But $\big((i-1)\frac{1}{\alpha}-\lfloor(i-1)\frac{1}{\alpha}\rfloor \big)_{i\in \mathbb{N}}$ is equidistributed in $[0,1)$ and therefore, $\big(i-\lfloor(i-1)\frac{1}{\alpha}\rfloor\alpha\big)_{i\in\mathbb{N}}$ is equidistributed in the interval $1+\alpha\cdot[0,1)=[1,1+\alpha)$ and hence equidistributed in $[1,1+\alpha]$. For any R-integrable function $f$ over $[1,1+\alpha]$ and especially for $\tau_1(x):=\frac{1}{3^x}$, it holds that
\begin{equation}\label{equi}
\lim_{m\rightarrow\infty}\frac{1}{m}\sum_{i=1}^m f\bigg(i-\bigg\lfloor(i-1)\frac{1}{\alpha}\bigg\rfloor\alpha\bigg)=\frac{1}{\alpha}\int_{1}^{1+\alpha}f(x)\;dx.\\[-10pt]
\end{equation}
\hfill $\Box$

\begin{lem}\label{mug}
Let $\alpha=\frac{\ln(2)}{\ln(3)}$. The geometric mean $\mu_g$ of the terms in $-\Phi_{\mathbb{R}}(1c_{\alpha})=\sum_{i=1}^{\infty}\frac{2^{d_i-1}}{3^i}:=\sum_{i=1}^{\infty}t_i$ is given by 
\begin{displaymath}
\mu_g=\lim_{m\rightarrow\infty}\bigg(\prod_{i=1}^{m}t_i\bigg)^{\frac{1}{m}}=3^{-{\alpha}\int_1^{1+\alpha}x\;dx}=\frac{1}{6}\sqrt{2}\approx 0.23570
\end{displaymath}
\end{lem}
\noindent{\it Proof}. Just as in the proof of Lemma \ref{mu} we have that $\big(i-\lfloor(i-1)\frac{1}{\alpha}\rfloor\alpha\big)_{i\in\mathbb{N}}$ is equidistributed in the interval $[1,1+\alpha]$. In relation (\ref{equi}) we choose the R-integrable function $f(x)=x$ over $[1,1+\alpha]$. Hence
\begin{displaymath}
\lim_{m\rightarrow\infty}\frac{1}{m}\sum_{i=1}^m \bigg(i-\bigg\lfloor(i-1)\frac{1}{\alpha}\bigg\rfloor\alpha\bigg)=\frac{1}{\alpha}\int_{1}^{1+\alpha}x\;dx=1+\frac{\alpha}{2}.
\end{displaymath}
Since $\tau_1(x)=\frac{1}{3^x}$ and $\tau_1\big(i-\lfloor(i-1)\frac{1}{\alpha}\rfloor\alpha \big)=t_i$, we have $i-\lfloor(i-1)\frac{1}{\alpha}\rfloor\alpha=\frac{\ln(t_i)}{-\ln(3)}$. Therefore, 
\begin{displaymath}
\lim_{m\rightarrow\infty}\frac{1}{m}\sum_{i=1}^m \ln(t_i)=-\ln(3)\bigg(1+\frac{\alpha}{2}\bigg).
\end{displaymath}
Finally,
\begin{displaymath}
\lim_{m\rightarrow\infty}\bigg(\prod_{i=1}^{m}t_i\bigg)^{\frac{1}{m}}=e^{\lim_{m\rightarrow\infty}\frac{1}{m}\sum_{i=1}^m \ln(t_i)}=e^{-\ln(3)(1+\frac{\alpha}{2})}=\frac{1}{3^{1+\frac{\alpha}{2}}}=\frac{1}{6}\sqrt{2}.\\[-15pt]
\end{displaymath}
\hfill $\Box$

\section{Pseudo Trajectories over \textbf{$1c_\alpha$} for \textbf{$\alpha=\ln(2)/\ln(3)$}}\label{pseudo}

Throughout this section, let $\alpha=\frac{\ln(2)}{\ln(3)}=[a_0;a_1,a_2,\ldots,a_k,\ldots]$ be the simple continued fraction expansion of $\alpha$ and $\big(\frac{p_k}{q_k}\big)_{k=0}^{\infty}$ its convergents. We have $\alpha=[0;1,1,1,2,2,3,1,5,2,23,2,\ldots]$ and $\big(\frac{p_k}{q_k}\big)_{k=0}^{\infty}=\big(\frac{0}{1}, \frac{1}{1}, \frac{1}{2}, \frac{2}{3}, \frac{5}{8}, \frac{12}{19}, \frac{41}{65}, \frac{53}{84}, \frac{306}{485}, \frac{15601}{24727},\ldots\big)$. The Christoffel words $w_k:=1z_k0$ for small $k$ appear in Table \ref{tab:ChristoffelWords}.
\begin{table}[htb]
	\centering
		\begin{tabular}{lllll}
		$b_{k}$ & $\frac{p_k}{q_k}$ &  $a_k$ &&   $w_k$ \\
		$b_2=\frac{2^1}{3^0}$ & $\frac{0}{1}$ &  $a_0=0$ &   $w_0=0$ & $0$ \\  
		$b_3=\frac{3^1}{2^1}$ & $\frac{1}{1}$ &  $a_1=1$ &   $w_1=1$ & $1$ \\
		$b_4=\frac{2^2}{3^1}$ & $\frac{1}{2}$ &  $a_2=1$ &   $w_2=w_1^{a_2}w_0$ & $1\;0$ \\
		$b_5=\frac{3^2}{2^3}$ & $\frac{2}{3}$ &  $a_3=1$ &   $w_3=w_1w_2^{a_3}$ & $1\;10$ \\
		$b_6=\frac{2^8}{3^5}$ & $\frac{5}{8}$ &  $a_4=2$ &   $w_4=w_3^{a_4}w_2$ & $110\; 110\; 10$ \\
		$b_7=\frac{3^{12}}{2^{19}}$ & $\frac{12}{19}$ &  $a_5=2$ &   $w_5=w_3w_4^{a_5}$ & $110\;(11011010)^2$ \\
		$b_8=\frac{2^{65}}{3^{41}}$ & $\frac{41}{65}$ &  $a_6=3$ &   $w_6=w_5^{a_6}w_4$ & $(110\; 11011010\; 11011010)^3\; 110\; 110\; 10$ \\
		$b_9=\frac{3^{53}}{2^{84}}$ & $\frac{53}{84}$ &  $a_7=1$ &   $w_7=w_5w_6^{a_7}$ & $110\;(11011010)^2\;(110 11011010 11011010)^3 110 110 10$ \\
		\end{tabular}
	\caption{Christoffel words}
	\label{tab:ChristoffelWords}
\end{table}

The sequence $\big(b_{k}\big)_{k=2}^{\infty}$ defined by $b_{k}=\frac{3^{p_{k-2}}}{2^{q_{k-2}}}$ for $k$ odd and $b_{k}=\frac{2^{q_{k-2}}}{3^{p_{k-2}}}$ for $k$ even decreases strictly monotone to $1$. From $p_k-a_kp_{k-1}=p_{k-2}$ and $q_k-a_kq_{k-1}=q_{k-2}$ the relation
\begin{equation}\label{b_k=}
b_k=b_{k+1}^{a_k}b_{k+2}\quad \textrm{ follows},
\end{equation}
and hence $b_{k+2}<b_{k+1}<b_{k+1}b_{k+2}\leq b_k$. Furthermore we have
\begin{displaymath}
-\Phi_{\mathbb{R}}(w_k)=\sum_{i=0}^{p_k-1}\frac{2^{\lfloor\frac{i}{\alpha}\rfloor}}{3^{1+i}}=\sum_{i=1}^{p_k}\frac{2^{\lfloor\frac{i-1}{\alpha}\rfloor}}{3^{i}}:=\sum_{i=1}^{p_k}x_i.
\end{displaymath}
Index $i$ enumerates the $p_k$ terms of $-\Phi_{\mathbb{R}}(w_k)$ in the order they appear in the sum or equivalently, giving the position of the $1$'s in $w_k$. In the next Lemma we sort the $x_i$'s by their value into ascending order to represent them as numbers of the semi-open interval $(\frac{1}{6},\frac{1}{3}]$. The fact that this sorting can be done with a permutation modulo $p_k$ shows the random-like distribution of the $x_i$'s in $(\frac{1}{6},\frac{1}{3}]$, even for small $k$.

\begin{lem}\label{ij}
Let $-\Phi_{\mathbb{R}}(w_k)=\sum_{i=1}^{p_k}x_i$ $(k>1)$. The permutation defined by
\begin{displaymath}
\sigma_k(i):=p_k-\bigg(\frac{(-1)^{k}}{p_{k-1}}(i-1)\bmod{p_k}\bigg)\quad \textrm{or equivalently,}\quad\sigma_k(i):=\frac{(-1)^{k-1}}{p_{k-1}}(i-1)\bmod{p^+_k}
\end{displaymath}
where the sign $+$ means that modulo $p_k$ must be reduced to a number of the set $\{1,2,\ldots,p_k\}$, yields the partition $Y_k:=\big(\frac{1}{6}<y_1<y_2<\cdots <y_j<\cdots <y_{p_k}=\frac{1}{3}\big)$ with $j=\sigma_k(i)$ and $y_{\sigma_k(i)}=x_i$. For convenience, we assign $y_0:=\frac{1}{6}$ as an auxiliary value. 
\end{lem}
\noindent{\it Proof}. The proof is based on a recursion over $k$. We show the first steps in Table \ref{tab:Permutation}. Instead of the $x_i$'s we have written only the $i$'s. In the row for $k=4$, numbers $4,2,5,3,1$ mean that $\frac{1}{6}<x_4<x_2<x_5<x_3<x_1=\frac{1}{3}$. Factors $b_5b_6,b_5,b_5b_6,b_5,b_5$ mean that $\frac{1}{6}b_5b_6=x_4$, $\;x_4b_5=x_2$, $\;x_2b_5b_6=x_5$, $\;x_5b_5=x_3$, $\;x_3b_5=x_1=\frac{1}{3}$.

For $k=4$, $w_4=(w_3)^2w_2=(110) (110) 10$, and the $1$'s correspond to $\frac{2^0}{3^1}$, $\frac{2^1}{3^2}$, $\frac{2^3}{3^3}$, $\frac{2^4}{3^4}$, $\frac{2^6}{3^5}$. By moving the first $(110)$ onto the second $(110)$, the corresponding $\big(\frac{2^0}{3^1}, \frac{2^1}{3^2}\big)$ are multiplied by $\frac{2^3}{3^2}=\frac{1}{b_5}<1$, wherein exponent $3$ is the length and exponent $2$ is the height of $(110)$. Thus $x_3<x_1$ and $x_4<x_2$. Therefore, $x_3b_5=x_1$ and $x_4b_5=x_2$. Moreover, $x_3\frac{2^3}{3^2}=x_3\frac{1}{b_5}=x_5$. So $x_5b_5=x_3$. However, in the row for $k=3$ we have $x_2b_4b_5=x_1$, and $b_4b_5=(b_5)^2b_6b_5$ holds by Relation (\ref{b_k=}). Hence $x_2b_6b_5=x_2(b_5b_6)=x_5$. From  $\frac{1}{6}b_4=x_2$, $b_4=(b_5)^2b_6$ and $x_4b_5=x_2$ it follows that $\frac{1}{6}b_5b_6=x_4$. 

For $k=4$ we have $\frac{2^3}{3^2}=\frac{1}{b_5}=\frac{1}{b_{k+1}}=\frac{2^{q_{k-1}}}{3^{p_{k-1}}}$. Also $x_i=\frac{2^{\lfloor\frac{i-1}{\alpha}\rfloor}}{3^{i}}$. Then
\begin{displaymath}
x_i\frac{1}{b_5}=\frac{2^{\lfloor\frac{i-1}{\alpha}\rfloor}}{3^{i}}\cdot\frac{2^{q_{k-1}}}{3^{p_{k-1}}}=\frac{2^{\lfloor\frac{i-1}{\alpha}\rfloor+q_{k-1}}}{3^{i+p_{k-1}}}=x_{i+p_{k-1}};
\end{displaymath}
\begin{displaymath}
x_i\frac{1}{b_5b_6}=\frac{2^{\lfloor\frac{i-1}{\alpha}\rfloor}}{3^{i}}\cdot\frac{2^{q_{k-1}}}{3^{p_{k-1}}}\cdot\frac{3^{p_k}}{2^{q_k}}=\frac{2^{\lfloor\frac{i-1}{\alpha}\rfloor+q_{k-1}-q_k}}{3^{i+p_{k-1}-p_k}}=x_{i+p_{k-1}-p_k}.
\end{displaymath}
Hence we conclude that in $(4, 2, 5, 3, 1)$, the neighbor of $i$ at the left side is $i+2 = i+p_{3} \bmod{p^+_4}$. 

\begin{table}[htb]
	\centering
		\begin{tabular}{c@{ }c@{ }c@{ }c@{ }c@{ }c@{ }c@{ }c@{ }c@{ }c@{ }c@{ }c@{ }c@{ }c@{ }c@{ }c@{ }c@{ }c@{ }c@{ }c@{ }c@{ }c@{ }c@{ }c@{ }c}
	$w_2=(w_1)^1w_0$&&&&& &&&&&&&&&&&&&&&&&&&$1$\\
	$\frac{1}{2}$&&&&& &&&&&&&&&&&&&&&&&&$b_3b_4$&\\
	\\
	$w_3=w_1(w_2)^{1}$&&&&& &&&&&$2$&&&&&&&&&&&&&&$1$\\
	$\frac{2}{3}$&&&&&&&&  & {$b_4$} &&&&&&&&{$b_4b_5$}&\\
	\\
	$w_4=(w_3)^{2}w_2$&&&&& &$4$&&&&$2$&&&&&&$5$&&&&$3$&&&&$1$\\
	$\frac{5}{8}$&&&&&$b_5b_6$  & \multicolumn{4}{c}{$\;\;b_5$} &&\multicolumn{5}{c}{$b_5b_6$}&&\multicolumn{3}{c}{$\;\;b_5$}&&\multicolumn{3}{c}{$b_5$}&\\
	\\
	$w_5=w_3(w_4)^{2}$&&$6$&&$11$& &$4$&&$9$&&$2$&&$7$&&$12$&&$5$&&$10$&&$3$&&$8$&&$1$\\	
				 $\frac{12}{19}$&$b_6$&&$b_6$&&$b_6b_7$&&$b_6$&&$b_6b_7$&&$b_6$&&$b_6$&&$b_6b_7$&&$b_6$&&$b_6b_7$&&$b_6$&&$b_6b_7$&\\	
					\end{tabular}
		\caption{Permutation $\sigma_k(i)$}
	\label{tab:Permutation}
\end{table}

For $k=5$, $w_5=w_3(w_4)^2=110 (11011010) (11011010)$. By moving the first $(11011010)$ onto the second $(11011010)$, the corresponding $\big(\frac{2^3}{3^3}, \frac{2^4}{3^4}, \frac{2^6}{3^5}, \frac{2^7}{3^6}, \frac{2^9}{3^7}\big)$=$(x_3, x_4, x_5, x_6, x_7)$ are multiplied by $\frac{2^8}{3^5}=b_6>1$. Elements $x_3, x_4, x_5$ already have been placed in the row for $k=4$. Thus $x_3b_6=x_8$, $x_4b_6=x_9$ and $x_5b_6=x_{10}$. However, $x_7\cdot\frac{1}{b_6}=x_2$ and $x_7b_6=x_{12}$. Also, $x_6\cdot\frac{1}{b_6}=\frac{1}{6}$ and $x_6b_6=x_{11}$. But $\frac{1}{6}$ is smaller than any $x_i$, and $x_1=\frac{1}{3}$ is larger than any other $x_i$. Therefore, it is convenient to define that $x_1$ maps onto $x_6$ by $b_6$. 

In the row for $k=4$ we have $x_2b_5b_6=x_5$, and it holds that $b_5b_6=(b_6)^2b_7b_6$ by Relation (\ref{b_k=}). Thus $x_2(b_6)^2b_7b_6=x_5$. But $x_2b_6=x_7$, $x_7b_6=x_{12}$, so that $x_{12}b_6b_7=x_5$. The same arguments hold for $\frac{1}{6}b_5b_6=x_4$. Moreover, we have $x_4b_5=x_2$ and it holds that $x_4(b_6)^2b_7=x_2$ by Relation (\ref{b_k=}). But $x_4b_6=x_9$, so that $x_9b_6b_7=x_2$.

Hence we conclude that in $(6, 11, 4, 9, 2, 7, 12, 5, 10, 3, 8, 1)$, the neighbor of $i$ on the right side is $i+5 = i+p_{4} \bmod{p^+_5}$. The neighbor of $i$ on the left side is $i-5=i-p_4\bmod{p^+_5}$.

The neighbor on the left side of $i$ is $i+p_{k-1} \bmod{p^+_k}$ for $k$ even and $i-p_{k-1} \bmod{p^+_k}$ for $k$ odd. The leftmost $i$ corresponding to $j=1$ in a row of Table \ref{tab:Permutation} is $p_k-1$ steps away from $i=1$. Thus $i_{left}=1+(-1)^kp_{k-1}(p_k-1)\bmod{p^+_k}$. Then $i=1+(-1)^kp_{k-1}(p_k-j)\bmod{p^+_k}$ for $j\in\{1, 2,\ldots, p_k\}$. Since $\gcd(p_{k-1}, p_k)=1$, the congruence $i\equiv 1+(-1)^kp_{k-1}(p_k-j)\pmod{p_k}$ can be divided by $p_{k-1}$ to find $j=\sigma_k(i)$.\hfill $\Box$

It follows from Lemma \ref{ij} that in Table \ref{tab:Permutation} the leftmost factors are $b_{k+1}b_{k+2}$ for $k$ even and $b_{k+1}$ for $k$ odd. The rightmost factors are $b_{k+1}$ for $k$ even and $b_{k+1}b_{k+2}$ for $k$ odd. Note that for $k=2$ we have $b_3b_4=b_2$ being the leftmost and rightmost factor.

\begin{lem}\label{bb=2}
It holds that $\big(b_{k+1}\big)^{p_k} \big(b_{k+2}\big)^{p_{k-1}}=2$ for the factors of Table \ref{tab:Permutation}.
\end{lem}
\noindent{\it Proof}. If $k$ even, then $\big(b_{k+1}\big)^{p_k} \big(b_{k+2}\big)^{p_{k-1}}=\frac{3^{p_{k-1}p_k}}{2^{q_{k-1}p_k}}\cdot\frac{2^{q_kp_{k-1}}}{3^{p_kp_{k-1}}}=2^{q_kp_{k-1}-q_{k-1}p_k}=2^{(-1)^k}=2$. Similar for $k$ odd.\hfill $\Box$

\begin{lem}\label{Jo}
Let $-\Phi_{\mathbb{R}}(w_k)=\sum_{i=1}^{p_k}x_i$ $(k>1)$. The terms $x_i$ sorted by their value into ascending order are given by the sequence
\begin{displaymath}
\big(y_j\big)_{j=1}^{p_k}=\bigg(\frac{2^{(-1+(-1)^{k-1}jq_{k-1})\bmod{q_k}}}{3\;^{p_k\;-\;(-(1+(-1)^{k-1}jp_{k-1})\bmod{p_k})}}\bigg)_{j=1}^{p_k}= \bigg(\frac{2^{(-1+(-1)^{k-1}jq_{k-1})\bmod{q_k}}}{3^{(1+(-1)^{k-1}jp_{k-1})\bmod{p^+_k}}}\bigg)_{j=1}^{p_k}.
\end{displaymath}
\end{lem}
\noindent{\it Proof}. We have $y_0=\frac{1}{6}=\frac{2^{-1}}{3^1}$. Then for $k$ odd $y_1=\frac{2^{-1}}{3^1}b_{k+1}=\frac{2^{-1}}{3^1}\cdot\frac{2^{q_{k-1}}}{3^{p_{k-1}}}=\frac{2^{-1+q_{k-1}}}{3^{1+p_{k-1}}}$, and for $k$ even  
$y_1=\frac{2^{-1}}{3^1}b_{k+1}b_{k+2}= \frac{2^{-1}}{3^1}\cdot\frac{3^{p_{k-1}}}{2^{q_{k-1}}}\cdot\frac{2^{q_{k}}}{3^{p_{k}}}=\frac{2^{-1-q_{k-1}+q_k}}{3^{1-p_{k-1}+p_k}}=\frac{2^{(-1-q_{k-1})\bmod{q_k}}}{3^{(1-p_{k-1})\bmod{p^+_k}}}$.\hfill $\Box$

\begin{lem}\label{yby}
Let $k>1$ be odd, $Y_k:=\big(\frac{1}{6}=y_0<y_1<y_2<\cdots <y_j<\cdots <y_{p_k}=\frac{1}{3}\big)$, and $b^{p_k}=2$ such that $b\in\mathbb{R}^+$. 
\begin{displaymath}
\textrm{If } 1\leq j\leq p_k-1,\textrm{ then }y_{j}<b^jy_0<y_{j+1}.
\end{displaymath}
\end{lem}
\noindent{\it Proof}. Note that $b^{p_k}y_0=2y_0=y_{p_k}=\frac{1}{3}$. By Lemma \ref{bb=2} we have $b^{p_k}=2=b_{k+1}^{p_k}b_{k+2}^{p_{k-1}}=b_{k+1}^{p_k-p_{k-1}}\big(b_{k+1}b_{k+2}\big)^{p_{k-1}}$ and therefore, $(p_k-p_{k-1})\ln(b_{k+1})+p_{k-1}\ln(b_{k+1}b_{k+2})=p_k\ln(b)$. In this equation we replace either $\ln(b_{k+1}b_{k+2})$ by $\ln(b_{k+1})$ or, $\ln(b_{k+1})$ by $\ln(b_{k+1}b_{k+2})$ and find that $b_{k+1}<b<b_{k+1}b_{k+2}$. Since $k$ is odd $b_{k+1}y_0=y_1<by_0$, and because of $b<b_{k+1}b_{k+2}<b_{k+1}^2<b_{k+1}^2b_{k+2}$, we have $by_0<y_2$. Hence $y_1<by_0<y_2$.

\noindent Since $k$ is odd, factors $b_{k+1}b_{k+2}$ arise for those $j=\sigma_k(i)$ where $i+p_{k-1}$ must be reduced modulo $p_k$. Therefore,
\begin{displaymath}
y_j=y_0\big(b_{k+1}\big)^j\big(b_{k+2}\big)^{\big\lfloor\frac{j\cdot p_{k-1}}{p_k}\big\rfloor} < y_0\big(b_{k+1}\big)^j\big(b_{k+2}\big)^{\frac{j\cdot p_{k-1}}{p_k}} = b^jy_0.
\end{displaymath}
Then
\begin{displaymath}
\frac{b^jy_0}{y_j}=\big(b_{k+2}\big)^{\big(\frac{p_{k-1}}{p_k}\cdot j - \big\lfloor\frac{p_{k-1}}{p_k}\cdot j\big\rfloor\big)}<b_{k+2}<\frac{y_{j+1}}{y_j}.
\end{displaymath}
Hence $y_j<b^jy_0<y_{j+1}$.\hfill $\Box$

\begin{lem}\label{bln2}
Let $b^{p_k}=2$ be such that $b\in\mathbb{R}^+$ and $k>1$. Then $p_k\frac{b-1}{b}<\ln(2)<p_k(b-1)$.
\end{lem}
\noindent{\it Proof}. By defining $\beta:[\frac{1}{6},\frac{1}{3}]\rightarrow [3,6]$ by $\beta(x):=\frac{1}{x}$, then $\int_{\frac{1}{6}}^{\frac{1}{3}}\frac{1}{x}\;dx=\ln(2)$. We consider the partition $B_k:=\big(\frac{1}{6}=y_0<by_0<\cdots <b^jy_0<\cdots <b^{p_k-1}y_0<2y_0=\frac{1}{3}\big)$ of the interval $[\frac{1}{6},\frac{1}{3}]$. The lower Riemann sum is $p_k\frac{b-1}{b}$ and the upper Riemann sum $p_k(b-1)$.\hfill $\Box$

\begin{lem}\label{pkmu}
Let $\mu=\frac{1}{6\ln(2)}$ be the arithmetic mean of Lemma \ref{mu} and $k>1$. Then  
\begin{displaymath}
\sum_{j=0}^{p_k-1}b^jy_0<p_k\mu<\sum_{j=1}^{p_k}b^jy_0=\sum_{j=0}^{p_k-1}b^jy_0+\frac{1}{6}.
\end{displaymath}
\end{lem}
\noindent{\it Proof}. From Lemma \ref{bln2} it follows that $\frac{1}{6}\cdot\frac{1}{p_k(b-1)}<\frac{1}{6}\cdot\frac{1}{\ln(2)}<\frac{1}{6}\cdot\frac{b}{p_k(b-1)}$. Hence $\frac{1}{6(b-1)}<p_k\mu<\frac{b}{6(b-1)}$. We have $\sum_{j=0}^{p_k-1}b^jy_0=y_0\frac{b^{p_k}-1}{b-1}=\frac{1}{6(b-1)}$ and $\sum_{j=1}^{p_k}b^jy_0=\frac{b}{6(b-1)}=\frac{1}{6(b-1)}+\frac{1}{6}$.\hfill $\Box$ 

\begin{lem}\label{y-pkmu}
Let $k>1$ be odd. Then $\big|-\Phi_{\mathbb{R}}(w_k)-p_k\mu\big|<\frac{1}{6}$.
\end{lem}
\noindent{\it Proof}. By Lemma \ref{yby} we have $0<\sum_{j=1}^{p_k}b^jy_0-\sum_{j=1}^{p_k}y_j <\frac{1}{6}$, since $\sum_{j=1}^{p_k}\big(b^jy_0-y_j\big)$ is the sum of $p_k$ disjoint sub-intervals, and this sum is less than the length of $(\frac{1}{6},\frac{1}{3}]$. We combine this inequality with Lemma \ref{pkmu}:
\begin{displaymath}
-\frac{1}{6}<\sum_{j=1}^{p_k}b^jy_0-\frac{1}{6}-\sum_{j=1}^{p_k}y_j=\sum_{j=0}^{p_k-1}b^jy_0-\sum_{j=1}^{p_k}y_j<p_k\mu-\sum_{j=1}^{p_k}y_j\;;
\end{displaymath}
\begin{displaymath}
p_k\mu-\sum_{j=1}^{p_k}y_j<\sum_{j=1}^{p_k}b^jy_0-\sum_{j=1}^{p_k}y_j<\frac{1}{6}.
\end{displaymath}
So we have $-\frac{1}{6}<p_k\mu-\sum_{j=1}^{p_k}y_j<\frac{1}{6}$. By Lemma \ref{Jo} we have $-\Phi_{\mathbb{R}}(w_k)=\sum_{i=1}^{p_k}x_i=\sum_{j=1}^{p_k}y_j$. Thus $-\frac{1}{6}<-\Phi_{\mathbb{R}}(w_k)-p_k\mu<\frac{1}{6}$.\hfill $\Box$ 

In a similar way it can be proved that Lemma \ref{y-pkmu} also holds for $k>1$ even. Moreover, in a bit lengthy proof it can be shown that we have
\begin{displaymath}
0<-\Phi_{\mathbb{R}}(w_k)-p_k\mu<\frac{1}{6}\quad\textrm{ ($k>1$ odd)}, \quad\quad\frac{1}{12}<-\Phi_{\mathbb{R}}(w_k)-p_k\mu<\frac{1}{4}\quad\textrm{ ($k>1$ even)}.
\end{displaymath}

\begin{lem}\label{ol}
Let $w$ be any factor of $1c_{\alpha}$ with length $\ell$. Then $-\Phi_{\mathbb{R}}(w)=\ell\mu+\emph{\textbf{o}}\big(\ell\big)$.
\end{lem}
\noindent{\it Proof}. The word $1c_{\alpha}$ is built up by Christoffel words $w_k$. Lemma \ref{y-pkmu} holds even for small $k$, so that the error of the approximation by $\ell\mu$ grows very slowly for an increasing $\ell$,  while the lengths $q_k$ of $w_k$ increase very fast.\hfill $\Box$   

\begin{lem}\label{ln2/ln3}
Let $\frac{p_k}{q_k}$ ($k=0,1,2,\ldots$) be the convergents of the simple 
continued fraction expansion of $\alpha=\frac{\ln(2)}{\ln(3)}$, $\mathcal{T}_{
1c_{\alpha}}(x)$the pseudo trajectory of any $x\in\mathbb{R}$ over $1c_{\alpha
}$ and $u(q_k)$ a prefix of $1c_{\alpha}$ with length $q_k$. Then
\begin{displaymath}
\lim_{k\rightarrow\infty}\frac{T_{u(q_k)}(x)}{q_k}=\frac{1}{6\ln(3)}
\quad\textrm{(for k odd)}\quad\textrm{and}\quad
\lim_{k\rightarrow\infty}\frac{T_{u(q_k)}(x)}{q_k}=\frac{1}{2\ln(3)}
\quad\textrm{(for k even)}.
\end{displaymath}
\end{lem}
\noindent{\it Proof}. [$k$ odd]. In this case $u(q_k)$ is a Christoffel word. 
We have $\big|-\Phi_{\mathbb{R}}(u(q_k))-p_k\mu\big|<\frac{1}{6}$ by Lemma 
\ref{y-pkmu}. Then $p_k\mu-\frac{1}{6}<\frac{2^{q_k}}{3^{p_k}}T_{u(q_k)}(0)<p_
k\mu+\frac{1}{6}$ by (\ref{rel}). Thus $\frac{{3^{p_k}}}{2^{q_k}}\big(p_k\mu-
\frac{1}{6}\big)<T_{u(q_k)}(0)<\frac{{3^{p_k}}}{2^{q_k}}\big(p_k\mu+\frac{1}{6
}\big)$ and therefore,
\begin{equation}\label{3mu}
\frac{{3^{p_k}}}{2^{q_k}}\bigg(\frac{p_k}{q_k}\mu-\frac{1}{6q_k}\bigg)< \frac{
T_{u(q_k)}(0)}{q_k} <\frac{{3^{p_k}}}{2^{q_k}}\bigg(\frac{p_k}{q_k}\mu+\frac{1
}{6q_k}\bigg).
\end{equation}
If $k\rightarrow\infty$, then $\frac{{3^{p_k}}}{2^{q_k}}\rightarrow 1$, $\frac
{p_k}{q_k}\rightarrow\frac{\ln(2)}{\ln(3)}$ and $\frac{1}{6q_k}\rightarrow 0$
. But $\frac{\ln(2)}{\ln(3)}\mu=\frac{1}{6\ln(3)}$. 

\noindent [$k$ even]. For $k$ even the prefix $u(q_k)$ of $1c_{\alpha}$ is not a Christoffel word $1z_k0$. We have now $u(q_k)=1z_k1$, and the height is $p_k+1$. We substitute $p_k$ by $p_k+1$ in (\ref{3mu}). If $k\rightarrow\infty$, then $\frac{{3^{p_k+1}}}{2^{q_k}}\rightarrow 3$, $\frac{p_k+1}{q_k}\rightarrow\frac{\ln(2)}{\ln(3)}$ and $\frac{1}{6q_k}\rightarrow 0$.

\noindent Finally, we apply Lemma \ref{Rout}. For any $x\in\mathbb{R}$ it holds that $T_{u(q_k)}(x)=T_{u(q_k)}(0)+\frac{3^{p_k}}{2^{q_k}} x$ if $k$ is odd, and $T_{u(q_k)}(x)=T_{u(q_k)}(0)+\frac{3^{p_k+1}}{2^{q_k}} x$ if $k$ is even. Thus $T_{u(q_k)}(x)\rightarrow T_{u(q_k)}(0)+x$ for $k$ odd, and $T_{u(q_k)}(x)\rightarrow T_{u(q_k)}(0)+3x$ for k even (see Example \ref{drei}).\hfill $\Box$

\section{A sufficient condition for the existence of $\mathbf{\Phi_{\mathbb{R}}(1c_v)}$}
\begin{lem}\label{suff}
Let $v$ be an aperiodic word such that 
$\underline{\lim}\; \big(\frac{h}{\ell}\big)_{\ell=1}^{\infty}=\frac{\ln(2)}{\ln(3)}=\alpha$, and let $n_{\ell}=h(\ell)-\lceil \ell\alpha\rceil$ be the relative height of $v$ over the Sturmian word $1c_{\alpha}$ at position $\ell$, such that $\big( \frac{2^{\ell}}{3^{\lceil\ell\alpha\rceil+n_{\ell}}}\big)_{\ell=1}^{\infty}\rightarrow 0$. Further, $\ell_j:=\max\{\ell\;|\;n_{\ell}=j\}$ for $j=0,1,2,\ldots$.

\begin{displaymath}
\textrm{If }\quad\overline{\lim}\;\bigg(\frac{\ell_{j+2}-\ell_{j+1}}{\ell_{j+1}-\ell_j}\bigg)_{j=1}^{\infty}< 3\quad\textrm{, then limits } \Phi_{\mathbb{R}}(1c_v)\leq\Phi_{\mathbb{R}}(v)\textrm{ exist}.
\end{displaymath}
\begin{displaymath}
\textrm{If }\quad\overline{\lim}\;\bigg(\frac{\ell_{j+2}-\ell_{j+1}}{\ell_{j+1}-\ell_j}\bigg)_{j=1}^{\infty}> 3\quad\textrm{, then limit } \Phi_{\mathbb{R}}(1c_v) \textrm{ does not exist.} 
\end{displaymath}
\end{lem}
\noindent{\it Proof}. The $\ell_j$ are well defined by Lemma \ref{max}. We apply the ratio test on the series (\ref{serie1cv}) in Lemma \ref{Phi1cv}.
 Since either $\lceil\ell_{j+1}\alpha\rceil - \lceil\ell_{j}\alpha\rceil=\lceil(\ell_{j+1}-\ell_j)\alpha\rceil$, or $\lceil\ell_{j+1}\alpha\rceil - \lceil\ell_{j}\alpha\rceil=\lceil(\ell_{j+1}-\ell_j)\alpha\rceil -1$, it holds for the ratio that
\begin{displaymath}
\rho_j=\frac{1}{3}\cdot \frac{2^{\ell_{j+1}-\ell_j}}{3^{\lceil\ell_{j+1}\alpha\rceil-\lceil\ell_j\alpha\rceil}}\cdot\frac{\Phi_{\mathbb{R}}(s'^{(j+1)})}{\Phi_{\mathbb{R}}(s'^{(j)})} = c_j\cdot \frac{2^{\ell_{j+1}-\ell_j}}{3^{\lceil(\ell_{j+1}-\ell_j)\alpha\rceil}}\cdot\frac{\Phi_{\mathbb{R}}(s'^{(j+1)})}{\Phi_{\mathbb{R}}(s'^{(j)})},\; \textrm{where }c_j=\frac{1}{3} \textrm{ or }c_j=1.  
\end{displaymath}
We apply Lemma \ref{ellj} and Lemma \ref{ol}. Then we have
\begin{displaymath}
\rho_j=c_j\cdot \frac{2^{\ell_{j+1}-\ell_j}}{3^{\lceil(\ell_{j+1}-\ell_j)\alpha\rceil}}\cdot \frac{(\ell_{j+2}-\ell_{j+1})\mu+\emph{\textbf{o}}\big(\ell_{j+2}-\ell_{j+1}\big)}{(\ell_{j+1}-\ell_{j})\mu+\emph{\textbf{o}}\big(\ell_{j+1}-\ell_{j}\big)},
\end{displaymath}
wherein $\frac{1}{3}<\frac{2^{\ell_{j+1}-\ell_j}}{3^{\lceil(\ell_{j+1}-\ell_j)\alpha\rceil}}<1$ by Lemma \ref{ln2ln3}, and wherein the last factor
\begin{displaymath}
\frac{(\ell_{j+2}-\ell_{j+1})\mu+\emph{\textbf{o}}\big(\ell_{j+2}-\ell_{j+1}\big)}   {(\ell_{j+1}-\ell_{j})\mu+\emph{\textbf{o}}\big(\ell_{j+1}-\ell_{j}\big)} = \frac{\frac{\ell_{j+2}-\ell_{j+1}}{\ell_{j+1}-\ell_j} +\frac{ \emph{\textbf{o}}\big(\ell_{j+2}-\ell_{j+1}\big)}{(\ell_{j+1}-\ell_{j})\mu}}   {1+\frac{\emph{\textbf{o}}\big(\ell_{j+1}-\ell_{j}\big)}{(\ell_{j+1}-\ell_{j})\mu}}
\end{displaymath}
approaches $\frac{\ell_{j+2}-\ell_{j+1}}{\ell_{j+1}-\ell_j}$ for an increasing $j$.

If $\overline{\lim}\;(\rho_j)<1$, the series (\ref{serie1cv}) converges, and if $\overline{\lim}\;(\rho_j)>1$, the series diverges. If we have infinitely often $c_j=\frac{1}{3}$ (that is, $\lceil\ell_{j+1}\alpha\rceil - \lceil\ell_{j}\alpha\rceil=\lceil(\ell_{j+1}-\ell_j)\alpha\rceil$ for infinitely many $j$), then the Lemma as claimed holds. Note that limit $\Phi_{\mathbb{R}}(v)$ exists by Lemma \ref{iff}. 

The assumption that there are only finitely many $c_j=\frac{1}{3}$ is in contradiction with Lemma \ref{piecewise}, namely, it yields divergence for $\overline{\lim}\big(\frac{\ell_{j+2}-\ell_{j+1}}{\ell_{j+1}-\ell_j}\big)>1$ while Lemma \ref{piecewise} states the convergence for $\overline{\lim}\;\big(\frac{h_{j+1}}{h_j}\big)<\frac{3}{2}$ of the corresponding heights.\hfill $\Box$

\section{Proof of Result \ref{funf}}\label{geom}

\begin{figure}[htbp]
\begin{center}
\epsfxsize=3.5in
\epsfbox{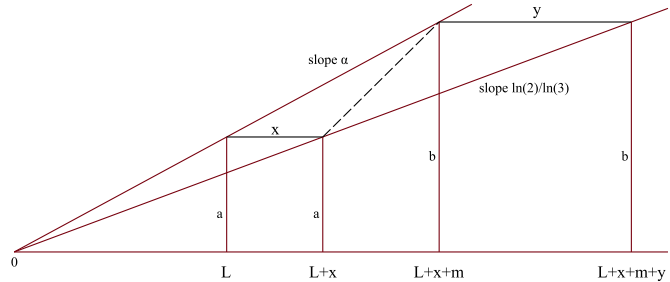}
\end{center}
\caption{Proof of Result \ref{funf}}
\label{propo}
\end{figure}

In Figure \ref{propo} where $L$ is the length of a prefix of $v$, we see that $a:=L\alpha=(L+x)\frac{\ln(2)}{\ln(3)}$. Thus $x=L\big(\alpha\frac{\ln(3)}{\ln(2)}-1\big)$. With $x$, $y$, and the dashed line we symbolize the geometric representation of the digits. The maximum of consecutive $0$'s not reaching the slope $\frac{\ln(2)}{\ln(3)}$ is given by $\lfloor x \rfloor$. These 0's are followed by $m$ 1's and $\lfloor y\rfloor$ 0's. It holds that
\begin{displaymath}
\frac{y}{x}=\frac{L+x+m+y}{L+x}=\frac{b}{a}=\frac{L+x+m}{L}.
\end{displaymath}
In the latter expression we substitute $m=\lfloor\frac{(L+x)\alpha-L\alpha}{1-\alpha}\rfloor=\lfloor\frac{x\alpha}{1-\alpha}\rfloor$ and $x=L\big(\alpha\frac{\ln(3)}{\ln(2)}-1\big)$. 
\begin{displaymath}
\frac{y}{x}=\frac{L\alpha\frac{\ln(3)}{\ln(2)}+\bigg\lfloor L(\alpha\frac{\ln(3)}{\ln(2)}-1)\frac{\alpha}{1-\alpha}\bigg\rfloor}{L} = \frac{L\alpha\frac{\ln(3)}{\ln(2)}+ L(\alpha\frac{\ln(3)}{\ln(2)}-1)\frac{\alpha}{1-\alpha}}{L}-\frac{\xi}{L}\quad (\textrm{for }0<\xi<1).
\end{displaymath}
We simplify and obtain $\frac{y}{x}=\frac{\alpha(\frac{\ln(3)}{\ln(2)}-1)}{1-\alpha}-\frac{\xi}{L}$. If $x\rightarrow\infty$, then $L\rightarrow\infty$. Thus $\lim_{x\rightarrow\infty}(\frac{y}{x})=\frac{\alpha(\frac{\ln(3)}{\ln(2)}-1)}{1-\alpha}$. Finally, $\lim_{x\rightarrow\infty}(\frac{\lfloor y\rfloor}{\lfloor x\rfloor})=\lim_{x\rightarrow\infty}(\frac{y}{x})$.\hfill $\Box$


\section{Lemmas concerning Example \ref{On}}\label{end}
In the following we calculate the exponents of 1's and 0's in the aperiodic word $v$ of Example \ref{On}. We have $v=$\\
$10\; 1^40\; 1^30\; 1^30\; 1^{4}0\; 1^{3}0\;\;\; 1^{4}0^2\; 1^{5}0^2\; 1^{5}0^2\; 1^{5}0^2\; 1^{5}0^2\; 1^{5}0^2\;\;\; 1^{5}0^{3}\; 1^{7}0^{3}\;1^{7}0^{3}\; 1^{7}0^{3}\; 1^{7}0^{3}\; 1^{7}0^{3}\;\;\; 1^{7}0^{4}\cdots$.
Each turn $n=1, 2, \ldots$ of the loop yields a factor $1^m0^{f(n)}$. 
Exponent $m$ depends on $n$. So we write $m_n$ (for instance $m_5=4$). Then $m_n=\lfloor\frac{(L_{n-1}+n)\alpha-H_{n-1}}{1-\alpha}\rfloor$ where $H_{n-1}$ is the height of a prefix with length $L_{n-1}$ having $n-1$ factors (for instance $H_4=1+4+3+3=11$ and $L_{4}=15$). Note that $H_0=L_0=0$. It holds that $L_{n-1}=Z_{n-1}+H_{n-1}$ where $Z_{n-1}$ is the number of 0's in the first $n-1$ factors.

\begin{lem}\label{expo}
Let $v$ be the aperiodic word of Example $\ref{On}$. Then
\begin{displaymath}
H_n=\bigg\lfloor\frac{(Z_{n-1}+n)\alpha}{1-\alpha}\bigg\rfloor,\quad m_{n+1}= \bigg\lfloor\frac{(Z_{n}+n+1)\alpha}{1-\alpha}\bigg\rfloor - \bigg\lfloor\frac{(Z_{n-1}+n)\alpha}{1-\alpha}\bigg\rfloor,
\end{displaymath}
\begin{displaymath}
Z_n=\frac{\lfloor\frac{n}{6}\rfloor(\lfloor\frac{n}{6}\rfloor+1)}{2}\cdot 6 + \bigg(\bigg\lfloor\frac{n}{6}\bigg\rfloor+1\bigg)\cdot (n\bmod{6})\;.
\end{displaymath}
\end{lem}
\noindent{\it Proof}.
\begin{displaymath}
m_n=\bigg\lfloor\frac{(L_{n-1}+n)\alpha-H_{n-1}}{1-\alpha}\bigg\rfloor = \bigg\lfloor\frac{(Z_{n-1}+H_{n-1}+n)\alpha-H_{n-1}}{1-\alpha}\bigg\rfloor = \bigg\lfloor\frac{(Z_{n-1}+n)\alpha}{1-\alpha}\bigg\rfloor-H_{n-1}.
\end{displaymath}
Obviously, $H_1=m_1$. If $n>0$, then $H_{n-1}=\sum_{i=1}^{n-1}m_i$. Therefore,
\begin{displaymath}
m_{n+1}=\bigg\lfloor\frac{(Z_{n}+n+1)\alpha}{1-\alpha}\bigg\rfloor-H_{n}=\bigg\lfloor\frac{(Z_{n}+n+1)\alpha}{1-\alpha}\bigg\rfloor-\sum_{i=1}^{n}m_i
\end{displaymath}
\begin{displaymath}
=\bigg\lfloor\frac{(Z_{n}+n+1)\alpha}{1-\alpha}\bigg\rfloor-\sum_{i=1}^{n-1}m_i-\bigg(\bigg\lfloor\frac{(Z_{n-1}+n)\alpha}{1-\alpha}\bigg\rfloor-H_{n-1}\bigg) = \bigg\lfloor\frac{(Z_{n}+n+1)\alpha}{1-\alpha}\bigg\rfloor - \bigg\lfloor\frac{(Z_{n-1}+n)\alpha}{1-\alpha}\bigg\rfloor
\end{displaymath}
\begin{displaymath}
\textrm{Hence }H_n=\bigg\lfloor\frac{(Z_{n-1}+n)\alpha}{1-\alpha}\bigg\rfloor\quad\textrm{and}\quad m_{n+1}= \bigg\lfloor\frac{(Z_{n}+n+1)\alpha}{1-\alpha}\bigg\rfloor - \bigg\lfloor\frac{(Z_{n-1}+n)\alpha}{1-\alpha}\bigg\rfloor.
\end{displaymath}
Since $f(n)=n-\lfloor\frac{5}{6}n\rfloor=(1,1,1,1,1,1,\;2,2,2,2,2,2,\ldots$), we can write $v=B_1B_2B_3\cdots B_i\cdots$ as a  product of blocks $B_i$ built up by 6 factors which have exactly $i$ 0's (for instance $B_2=1^{4}0^2\; 1^{5}0^2\; 1^{5}0^2\; 1^{5}0^2\; 1^{5}0^2\; 1^{5}0^2$). Thus $B_i=1^{m_{6i-5}}0^i\;1^{m_{6i-4}}0^i\;\cdots 1^{m_{6i}}0^i$.

We count the number $Z_n$ of 0's in the first $n$ factors. From $n=\lfloor\frac{n}{6}\rfloor 6+(n\bmod{6})$ and $1+2+3+\cdots +\lfloor\frac{n}{6}\rfloor =\frac{\lfloor\frac{n}{6}\rfloor(\lfloor\frac{n}{6}\rfloor+1)}{2}$, it follows that
\begin{displaymath}
Z_n=\frac{\lfloor\frac{n}{6}\rfloor(\lfloor\frac{n}{6}\rfloor+1)}{2}\cdot 6 + \bigg(\bigg\lfloor\frac{n}{6}\bigg\rfloor+1\bigg)\cdot (n\bmod{6})\;.\\[-15pt]
\end{displaymath}
\hfill $\Box$

Function $m_n$ increases, because $Z_n$ is strictly monotone increasing. But, it holds that
\begin{lem}\label{-10123}
\begin{displaymath}
-1\leq m_{n+1}-m_n\leq 3\;.\\[-15pt]   
\end{displaymath}
\end{lem}
\noindent{\it Proof}.
\begin{displaymath}
m_{n+1}-m_n=\bigg (\frac{(Z_{n}+n+1)\alpha}{1-\alpha}-\xi_1\bigg ) - 2\bigg(\frac{(Z_{n-1}+n)\alpha}{1-\alpha}-\xi_2\bigg ) +\bigg (\frac{(Z_{n-2}+n-1)\alpha}{1-\alpha}-\xi_3\bigg ), 
\end{displaymath}
where $0<\xi_k<1$. Thus $m_{n+1}-m_n=\frac{\alpha}{1-\alpha}\bigg(Z_{n}+n+1-2(Z_{n-1}+n)+Z_{n-2}+n-1\bigg)-\xi_1+2\xi_2-\xi_3=\frac{\alpha}{1-\alpha}(Z_{n}-Z_{n-1}+Z_{n-2}-Z_{n-1})-\xi_1+2\xi_2-\xi_3$, and $-2<-\xi_1+2\xi_2-\xi_3<2$.

Exponent $m_{n+1}$ depends on $Z_{n-2}$, $Z_{n-1}$, and $Z_{n}$. Factor $1^{m_{n+1}}0^i$ belongs to block $B_i$ with $i=\lfloor\frac{n}{6}\rfloor +1$. We have three cases: factor $1^{m_{n+1}}0^i$ is the first, the second, or another term of $B_i$.
\begin{center}
\begin{tabular}{rllccrllccrll}
\multicolumn{3}{l}{term 1}&&&\multicolumn{3}{l}{term 2}&&&\multicolumn{3}{l}{$3\leq$ term $\leq 6$}\\
$Z_n-Z_{n-1}$&=&$i-1$&&&$Z_n-Z_{n-1}$&=&$i$&&&$Z_n-Z_{n-1}$&=&$i$\\
$Z_{n-2}-Z_{n-1}$&=&$-i+1$&&&$Z_{n-2}-Z_{n-1}$&=&$-i+1$&&&$Z_{n-2}-Z_{n-1}$&=&$-i$\\
\end{tabular}
\end{center}

Hence, if $1^{m_{n+1}}0^i$ is not the second term, then $-2<1.709\cdot 0-\xi_1+2\xi_2-\xi_3<2$, so that $m_{n+1}-m_n \in\{-1,0,1\}$. However, if it is the second term, then  $1.709-2<1.709\cdot 1-\xi_1+2\xi_2-\xi_3<1.709+2$, so that $m_{n+1}-m_n \in\{0,1,2,3\}$.\hfill $\Box$

\begin{lem}\label{quot}
The quotient of the heights of consecutive factors converges to $1$. The same holds for the quotient of lengths.
\end{lem}
\noindent{\it Proof}. Since $m_n-1\leq m_{n+1}$, we have (for $n\geq 2)$
\begin{displaymath}
1-\frac{1}{m_n}\leq\frac{m_{n+1}}{m_n} = \frac{ \bigg\lfloor\frac{(Z_{n}+n+1)\alpha}{1-\alpha}\bigg\rfloor - \bigg\lfloor\frac{(Z_{n-1}+n)\alpha}{1-\alpha}\bigg\rfloor }{ \bigg\lfloor\frac{(Z_{n-1}+n)\alpha}{1-\alpha}\bigg\rfloor - \bigg\lfloor\frac{(Z_{n-2}+n-1)\alpha}{1-\alpha}\bigg\rfloor } < \frac{ (Z_n-Z_{n-1}+1)+1 }{(Z_{n-1}-Z_{n-2}+1)-1 }=\frac{i+a}{i+b}=\frac{1+\frac{a}{i}}{1+\frac{b}{i}}\;.
\end{displaymath}
Depending on the position of $1^{m_{n+1}}0^i$ in block $B_i$, we have $a\in\{1,2\}$ and $b\in\{-1,0\}$ as possible values for $a$ and $b$.\hfill $\Box$

\begin{figure}[htbp]
\begin{center}
\epsfxsize=6in
\epsfbox{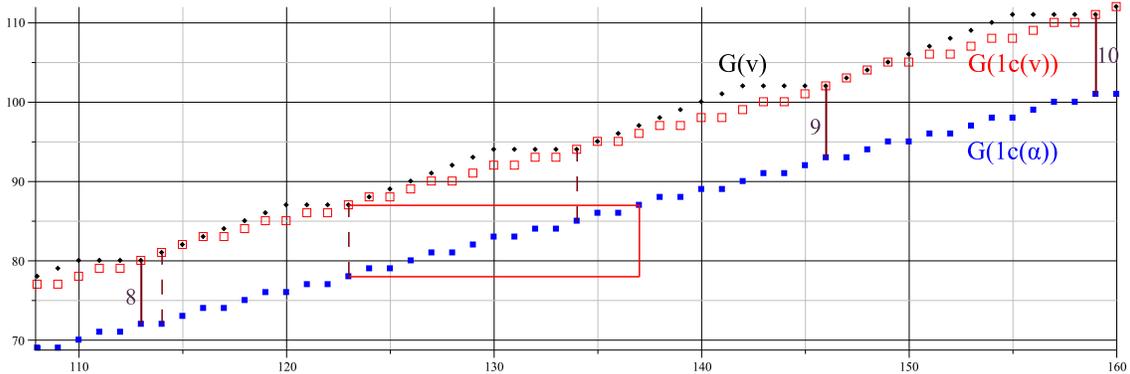}
\end{center}
\caption{Proof of Lemma \ref{max3}}
\label{detalle}
\end{figure}

\begin{lem}\label{max3}
Let $v$ be the aperiodic word of Example $\ref{On}$ and $\ell_j=\max\{\ell\;|\;n_{\ell}=j\}$ as in Lemma \ref{max} and Lemma \ref{iff}. In interval $[1+\ell_j,\ell_{j+1}]$ lie at most three complete factors of the form $1^m0^{f(n)}$. If the end of the third factor is at $L_n$, then $L_n\leq\ell_{j+1}\leq L_n+2$.
\end{lem}
\noindent{\it Proof}. We explain the idea of the proof in Figure \ref{detalle}. We see that $\ell_8=113$, $\ell_9=146$ and $\ell_{10}=159$. The end of block $B_3=\cdots\cdot\cdot 1^70^3$ and the beginning of block $B_4=1^70^4\;1^80^4\;1^90^4\cdots$ is at $L=123$ (Note that the first digit of $B_4$ corresponds to $L=124$). Thus $\ell_9-\ell_8=33$, and the length of the three factors $1^70^3\;1^70^4\;1^80^4$ is $33$ as well.

Associated to $1^70^3\;1^70^4\;1^80^4$ is the Sturmian word $s'^{(8)}=1s_{114}s_{115}\cdots s_{145}$ (see Lemma \ref{iff}) which lies between $\ell_8+1$ and $\ell_9$ nine units above $1c_{\alpha}$. By construction of $v$ in each turn of the loop, the slope is re-adjusted, so that a new factor approximates $s'^{(8)}$ in the 9'th floor as closely as possible and therefore, $G(v)$ and $G(s'^{(8)})$ coincide in points $(123,87)$, $(134,94)$, $(146,102)$. 

Factor $1^70^3$ of $B_3$ has number $n=18$. Thus $f(18)=18-15$. If we adjoin $18$ zeroes instead of only $3$, starting at (120,87), then the slope would be slightly less than $\ln(2)/\ln(3)$, and the last $0$ 
would give point $(138,87)$ of $G(0c_{\alpha})$. This is so, because $m$ and $0c_{\alpha}$ are calculated with the function floor. Hence the last but one $0$ yields a point of $G(1c_{\alpha})$ as indicated in Figure \ref{detalle} by a horizontal line segment. 

We conclude that the height of $1c_v$ at L=123 can be calculated in two different ways with the same result: $\lceil(123+14)\alpha\rceil=\lceil 123\alpha\rceil+9\;$. Note that $14=\lfloor\frac{5}{6}\cdot 18\rfloor-1$. Since $L_n$ is the length of a prefix of $v$ ending with the $n'th$ factor, we see in Figure \ref{detalle} that $\lceil(L_{18}+14)\alpha\rceil=\lceil L_{18}\alpha\rceil+9,\; \lceil(L_{19}+14)\alpha\rceil=\lceil L_{19}\alpha\rceil+9,\; \lceil(L_{20}+15)\alpha\rceil=\lceil L_{20}\alpha\rceil+9$.

The equation $\lceil(L_{n}+x)\alpha\rceil=\lceil L_{n}\alpha\rceil+9$ has only two integer solutions, namely $\lfloor 9/\alpha\rfloor=14$ and $\lceil 9/\alpha\rceil=15$. But in sequence $(\lfloor\frac{5}{6}\cdot n\rfloor-1)_{n=1}^{\infty}=(-1,0,1,2,3,4,4,5,6,7,8,9,9,10\ldots)$ the unique positive integers which appear $2$ times have the form $4+5k$ (for instance 14). Hence distance $9$ cannot occur more than $3$ times as an end of consecutive factors. The same holds for other distances.
 
Let $D_n$ be the distance between $1c_{\alpha}$ and $v$ (resp. $1c_v$) at the end of the $n$'th factor. Then
\begin{displaymath}
D_n=\bigg\lceil\bigg(L_{n}+\bigg\lfloor\frac{5}{6}\cdot n\bigg\rfloor-1\bigg)\alpha\bigg\rceil - \lceil L_{n}\alpha\rceil,
\end{displaymath}
where $L_n=Z_n+H_n$. We can check that in the sequence of distances no number is repeated more than 3 times:\\
$(-1, 0, 1, 1, 2, 2, 3, 3, 4, 4, 5, 5, 5, 6, 7, 8, 8, 9, 9, 9, 10, 11, 11, 12, 12, 13, 13, 14, 15, 15, 15,\ldots)$.

\begin{figure}[htbp]
\begin{center}
\epsfxsize=4.5in
\epsfbox{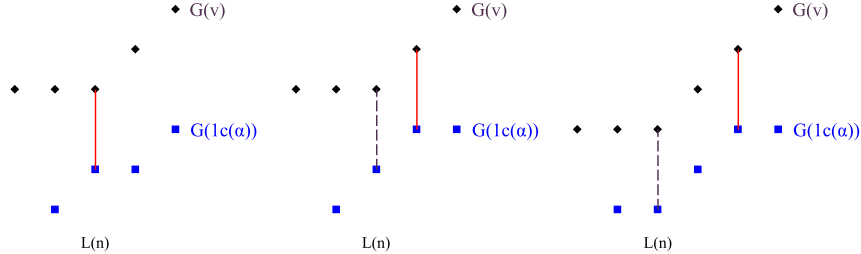}
\end{center}
\caption{$L_n\leq\ell_{j+1}\leq L_n+2$ }
\label{3c}
\end{figure}

In Figure \ref{detalle} distances $n_{\ell}=8$, $n_{\ell}=9$ and $n_{\ell}=10$ appear exactly at the end of any factor $1^m0^{f(n)}$. However, any other $n_\ell$'s can fall one or two places beyond the end of a factor. Figure \ref{3c} shows the $3$ possibilities. Since $1c_{\alpha}$ does not have more than two consecutive $1$'s, the $n_{\ell}$ fall at most two places beyond the end of a factor.

\hfill $\Box$ 

\begin{lem}\label{h/h}
Let $v$ be the aperiodic word of Example $\ref{On}$. Then $\overline{\lim}_{j\rightarrow\infty}(\frac{\ell_{j+2}-\ell_{j+1}}{\ell_{j+1}-\ell_j})\leq 3$ and 
$\overline{\lim}_{j\rightarrow\infty}(\frac{h'_{j+1}}{h'_{j}})\leq 3$, where $h'_j=\lceil (\ell_{j+1}-\ell_j)\alpha\rceil$ are the corresponding heights.
\end{lem}
\noindent{\it Proof}. In sequence $D_n$ we look for the largest $n$ where the distance $n_{\ell_j}$ for a fixed $j$ appears at $L_n$. We obtain\\ 
\noindent $n_{\mathrm{\max}}\in \{2,4,6,8,10,13,14,15,17,20,21,23,25,27,28,31,33,34,37,38,40,43,44,46,47,\ldots\}$.\\  
For instance, the ultimate distance $1$ occurs at $L_4$, the ultimate distance $2$ at $L_6$, and the ultimate distance $4$ at $L_{10}$. The first differences of the foregoing sequence are\\
$(2,2,2,2,3,1,1,2,3,1,2,2,2,1,3,2,1,3,1,2,3,1,2,1,3,1,2,3,1,2,3,1,2,2,2,1,3,2,1,\ldots)$.\\
They yield the number of factors in interval $(\ell_j,\ell_{j+1}]$. For instance, in $(\ell_0,\ell_{_1}]$ we have two factors and in $(\ell_4,\ell_{_5}]$ three factors.

The quotient $\frac{L_{n+2}-L_{n+1}}{L_{n+1}-L_n}$ of the lengths of consecutive factors converges to $1$ by Lemma \ref{quot}. Therefore, the quotient $\frac{L_{n+3}-L_{n}}{L_{n}-L_{n-1}}$ of the total length of three factors preceded by the length of one factor converges to $3$, provided that there are infinitely many pairs $(1,3)$ in the sequence above. From $L_{n}\leq\ell_{j+1}\leq L_{n}+2$ we conclude that $L_{n+3}-L_{n}\leq\ell_{j+2}-\ell_{j+1}\leq L_{n+3}-L_{n}+2$ and $L_{n}-L_{n-1}\leq\ell_{j+1}-\ell_{j}\leq L_{n}-L_{n-1}+2$. If there are infinitely many pairs $(1,3)$, then $\lim_{j\rightarrow\infty}\frac{\ell_{j+2}-\ell_{j+1}}{\ell_{j+1}-\ell_j} = \lim_{n\rightarrow\infty}\frac{L_{n+3}-L_{n}+a_n}{L_{n}-L_{n-1}+b_n} = \lim_{n\rightarrow\infty}\frac{L_{n+3}-L_{n}}{L_{n}-L_{n-1}}=3$, where $a_n,b_n\in\{0,1,2\}$. 
Similar arguments hold for any other pair in the sequence above. There is at least one pair occurring infinitely many times. The set of possible limit-points is $\{\frac{1}{3}, \frac{1}{2}, \frac{2}{3}, 1, \frac{3}{2}, 2, 3\}$. Experimentation shows that all these and no other limit-points actually exist.

Let $1c_v$ be the associated word to $v$. The quotient of lengths (resp. heights) of consecutive factors of $v$ converges to $1$, and $\ell_{j+1}-\ell_j$ is  not larger than the length of $3$ factors. The Sturmian words $s'^{(j)}$ lie between $1+\ell_{j}$ and $\ell_{j+1}$ and have the same length (and the same height) $\ell_{j+1}-\ell_j$ as the corresponding factors of $v$. Therefore,  $\overline{\lim}_{j\rightarrow\infty}(\frac{\ell_{j+2}-\ell_{j+1}}{\ell_{j+1}-\ell_j})\leq 3$. Let $h'_j=\lceil (\ell_{j+1}-\ell_j)\alpha\rceil$ be the height of $s'^{(j)}$. Then  $\overline{\lim}_{j\rightarrow\infty}(\frac{h'_{j+1}}{h'_{j}})\leq 3$.\hfill $\Box$

\begin{lem}\label{exist}
Let $v$ be the aperiodic word of Example \ref{On}. Limit $\Phi_{\mathbb{R}}(v)$ exists.
\end{lem}
\noindent{\it Proof}. Lemma \ref{suff} allows no decision about the existence of the limit, because in Example \ref{On}  $\;\overline{\lim}\;\big(\frac{\ell_{j+2}-\ell_{j+1}}{\ell_{j+1}-\ell_j}\big)_{j=1}^{\infty}\leq 3$ holds.

In Example \ref{On} we replace the function $f(n)$ by $f'(n):=n-\lfloor\frac{1}{4}n\rfloor$ and get an aperiodic word $v'$ with a faster increasing number of consecutive $0$'s than in $v$.\\
$v'=10\;1^40^2\;1^50^3\;1^70^3\;1^60^4\;1^90^5\;1^{10}0^6\;1^{12}0^6\;1^{12}0^7\;1^{14}0^8\;1^{15}0^9\;1^{17}0^9\;1^{17}0^{10}\;1^{17}0^{11}\;1^{19}0^{12}\;\cdots\; $.\\
\noindent We count the number $Z'_n$ of $0$'s contained in the first $n$ factors $1^{m'}0^{f'(n)}$ of $v'$.
\begin{displaymath}
Z'_n=\frac{n(n+1)}{2}-\frac{\big(\big\lfloor\frac{n}{4}\big\rfloor-1\big)\big\lfloor\frac{n}{4}\big\rfloor}{2}\cdot 4-\bigg\lfloor\frac{n}{4}\bigg\rfloor\bigg((n \bmod{4})+1\bigg).               
\end{displaymath}
Distance $D'_n$ between $1c_{\alpha}$ and $v'$ at the end of the $n$'th factor is given by
\begin{displaymath}
D'_n=\bigg\lceil\bigg(L'_{n}+\bigg\lfloor\frac{1}{4}\cdot n\bigg\rfloor-1\bigg)\alpha\bigg\rceil - \lceil L'_{n}\alpha\rceil,
\end{displaymath}
where $L'_n=Z'_n+H'_n$. In the sequence\\
\noindent $\scriptstyle{(D'_n)_{n=1}^{\infty}=(-1, -1, -1, 0, 0, 0, 0, 1, 1, 1, 0, 1, 1, 1, 1, 2, 2, 2, 2, 2, 2, 3, 2, 3, 3, 3, 3, 4, 4, 4, 4, 4, 4, 4, 4, 5, 5, 5, 5, 6, 6, 6, 6, 6, 6, 6, 6, 7, 7, 7, 7, 8, 7, 7, 7, 8, 8, 8, 8,\ldots)}$,\\
no distance is repeated more than $8$ times. The ultimate distance $0$ occurs at $L'_{11}$, the ultimate $1$ at $L'_{15}$, and the ultimate $2$ at $L'_{23}$, etc. In fact,\\
$n'_{\mathrm{\max}}\in\{11,15,23,27,35,39,47,55,59,67,71,79,84,91,99,103,111,119,123,131,135, \ldots \}$.\\
We calculate the first differences: \\
$\scriptstyle{(4, 8, 4, 8, 4, 8, 8, 4, 8, 4, 8, 5, 7, 8, 4, 8, 8, 4, 8, 4, 8, 5, 7, 6, 6, 7, 5, 8, 7, 5, 7, 5, 8, 8, 4, 4, 8, 6, 6, 8, 4, 8, 8, 4, 8, 4, 8, 8, 4, 8, 4, 8, 8, 4, 7, 5, 8, 4, 8, 8, 4, 8, 4, 8, 8, 4, 8, 4,\ldots})$.\\
They yield the number of factors in the interval $(\ell'_j,\ell'_{j+1}]$. It follows that the set of possible limit-points is $\{\frac{1}{2},\frac{5}{8},\frac{5}{7},\frac{3}{4},\frac{6}{7},\frac{7}{8},1,\frac{8}{7},\frac{7}{6},\frac{4}{3},\frac{7}{5},\frac{8}{5},2\}$. Hence, $\overline{\lim}_{j\rightarrow\infty}(\frac{\ell'_{j+2}-\ell'_{j+1}}{\ell'_{j+1}-\ell'_j})\leq 2$. The limit $\Phi_{\mathbb{R}}(v')$ exists by Lemma \ref{suff}, and we have $\Phi_{\mathbb{R}}(v')\approx -8.38314$.

Let $G(v)$ and $G(v')$ be the geometric representations of $v$ and $v'$. They coincide in the first $7$ points. There exists a positive integer $L$, namely $L=14$, such that the points of $G(v')$ fall below $G(v)$ for all $\ell\geq L$, because the prefixes of $v'$ occupy many places with $0$'s while $v$ takes advantage filling up these places with $1$'s. Thus, the suffix of $v'$ starting at $L=13$ can be lifted onto the corresponding suffix of $v$ so that $\Phi_{\mathbb{R}}(v')$ is a lower bound of $\Phi_{\mathbb{R}}(v)$.\hfill $\Box$

\noindent {\bf Result \ref{sieben}, Sketch of proof.}

\noindent Replacing the function $f'(n)$ of the foregoing proof by $f''(n):=n-\lfloor\frac{1}{q}n\rfloor$ for a large positive integer $q$ we get an aperiodic word $v''$ with a very fast increasing number of consecutive $0$'s such that $\Phi_{\mathbb{R}}(v'')$ is a lower bound of $\Phi_{\mathbb{R}}(v)$ for any $v$ considered in Result \ref{sieben}. For $q\geq 2$ we have
\begin{displaymath}
Z''_n=\frac{n(n+1)}{2}-\frac{\big(\big\lfloor\frac{n}{q}\big\rfloor-1\big)\big\lfloor\frac{n}{q}\big\rfloor}{2}\cdot q-\bigg\lfloor\frac{n}{q}\bigg\rfloor\bigg((n \bmod{q})+1\bigg).               
\end{displaymath}
The distance $D''_n$ between $1c_{\alpha}$ and $v''$ at the end of the $n$'th factor is given by
\begin{displaymath}
D''_n=\bigg\lceil\bigg(L''_{n}+\bigg\lfloor\frac{1}{q}\cdot n\bigg\rfloor-1\bigg)\alpha\bigg\rceil - \lceil L''_{n}\alpha\rceil,
\end{displaymath}
where $L''_n=Z''_n+H''_n$. We can verify that in the sequence $(D''_n)_{n=1}^{\infty}$ there is no distance repeating more than $2q+1$ times and that $\overline{\lim}\leq 2$ for any $q\geq 2$.\hfill $\Box$

\noindent {\bf Evidence that $\Phi_{\mathbb{R}}(v)$ of Example \ref{On} is irrational.}

\begin{figure}[htb]
\begin{center}
\epsfxsize=3in
\epsfbox{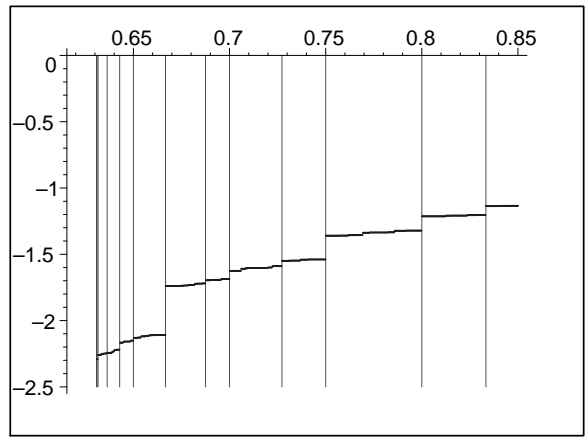}
\end{center}
\caption{Devil's staircase}
\label{diablo1}
\end{figure}

\noindent In the pseudo-code we have $\alpha$ as parameter. Let $\alpha$ be any real number (irrational or rational) such that $\frac{\ln(2)}{\ln(3)}<\alpha<1$. The limit $\Phi_{\mathbb{R}}(v_{\alpha})$
with such a parameter exists by Lemma \ref{limarriba}. Moreover by Lemma \ref{orbitarriba}, it is an irrational number even in case of a rational $\alpha$ since $v_{\alpha}$ is aperiodic and $\Phi_{\mathbb{R}}(v_{\alpha})\notin\mathbb{Q}_{odd}$. The plot of the function $\alpha\rightarrow \Phi_{\mathbb{R}}(v_{\alpha})$ in Figure \ref{diablo1} reveals gaps at the points with rational $\alpha$, but $\frac{\ln(2)}{\ln(3)}$ is irrational. The biggest gaps occur at the rational $\frac{p}{q}$ with small $p$ and $q$. Hence we may expect that $\lim_{\alpha\rightarrow \ln(2)/\ln(3)}(\Phi_{\mathbb{R}}(v_{\alpha}))=\Phi_{\mathbb{R}}(v)$.
 
We note that the gaps arrise in a natural way, because the aperiodic words $v_{\alpha}$ behave likely as Sturmian words. Actually, the expression $m_{n}(\alpha)= \big\lfloor (Z_{n-1}+n)\frac{\alpha}{1-\alpha}\big\rfloor - \big\lfloor (Z_{n-2}+n-1)\frac{\alpha}{1-\alpha}\big\rfloor$ has aspect of a mechanical word with slope $\frac{\alpha}{1-\alpha}$, but the coefficients $Z_{n-1}+n$ and $Z_{n-2}+n-1$ are not consecutive integers.


\end{document}